\documentclass[12 pt]{article}

\usepackage{amscd, amssymb, amsthm, amsmath,eucal, verbatim}
\usepackage[dvips]{graphics}

\newcommand{\bpf}{\noindent {\em Proof.} }
\newcommand{\epf}{\qed \vspace{+10pt}}

\newcommand{\proj}{\mathbf P}

\newcommand{\rarr}{\rightarrow}
\newcommand{\oh}{{\mathcal{O}}}
\newcommand{\com}{\mathbb{C}}
\newcommand{\Q}{\mathbb{Q}}
\newcommand{\A}{A_\Gamma}
\newcommand{\Z}{\mathbb{Z}}
\newcommand{\R}{\mathbb{R}}
\newcommand{\mgn}{\overline{M}_{g,n}}
\newcommand{\M}{\overline{M}}
\newcommand{\p}{\partial}
\newcommand{\hodge}{{\mathbb{E}}}
\newcommand{\llr}{{\mathbb{L}}}

\newcommand{\br}{{br}}

\def\scup{\mathbin{\text{\scriptsize$\cup$}}}
\def\scap{\mathbin{\text{\scriptsize$\cap$}}}

\newcommand{\E}{\mathsf{E}}
\newcommand{\V}{\mathsf{V}}
\newcommand{\T}{\mathsf{T}}
\newcommand{\bG}{\mathsf{G}}
\newcommand{\bH}{\mathsf{H}}
\newcommand{\bGt}{\mathsf{G}^3}
\newcommand{\bGth}{\mathsf{G}^{\ge 3}}
\newcommand{\mm}{\mathsf{m}}
\newcommand{\la}{\left\langle}
\newcommand{\ra}{\right\rangle} 
\newcommand{\Su}{\Sigma} 
\newcommand{\bem}{{\boldsymbol\emptyset}} 

\DeclareMathOperator{\val}{val} 
\DeclareMathOperator{\tk}{tk} 
\DeclareMathOperator{\Prob}{Prob} 
\DeclareMathOperator{\per}{per} 
\DeclareMathOperator{\Cell}{Cell} 
\DeclareMathOperator{\hmt}{hmt} 
\DeclareMathOperator{\und}{und}
\DeclareMathOperator{\Asm}{Asm} 
\DeclareMathOperator{\asm}{asm} 
\DeclareMathOperator{\Aut}{Aut} 
\DeclareMathOperator{\tr}{tr} 
\DeclareMathOperator{\Map}{Map}
\DeclareMathOperator{\map}{map}
\DeclareMathOperator{\Ai}{Ai}

\newcommand{\KAi}{\mathsf{K}_{\textup{Ai}}}

\newtheorem{pr}{Proposition}[section]
\newtheorem{lm}[pr]{Lemma}
\newtheorem{tm}{Theorem}
\newtheorem{cor}[pr]{Corollary}

\theoremstyle{definition}

\newtheorem{definition}[pr]{Definition}

\numberwithin{equation}{section}

\begin{document}
\title{Gromov-Witten theory, Hurwitz numbers,
and Matrix models, I}
\author{A.~Okounkov and R.~Pandharipande}
\date{15 March 2001}
\maketitle

\tableofcontents

\part{Overview}

\section{Introduction}
\subsection{Gromov-Witten theory, matrix models, and integrable
hierarchies}

Our goal here is to present a new path connecting
the intersection theory of the  moduli space 
$\overline{M}_{g,n}$
of stable
curves  
to the theory of matrix models.
The relationship
between these subjects was first discovered by E. Witten
in 1990 
through a study of 2 dimensional quantum
gravity \cite{W}. The path integral in the quantum gravity theory on
a genus $g$ topological surface $\Sigma_g$  
admits two natural interpretations.
First, the free energy of the theory may be expressed as a generating
series of tautological   
intersections products in $\overline{M}_{g,n}$.
A second approach via approximations by singular metrics on 
$\Sigma_g$ is connected to the asymptotic expansions
of Hermitian
matrix integrals --- the theory of matrix models.
Witten conjectured that the
Korteweg-de Vries  equations 
(known to control the associated matrix models)
governed the intersection theory of $\overline{M}_{g,n}$.
These KdV equations completely determine 
the tautological intersections in $\overline{M}_{g,n}$. 
As there was no previous mathematical approach to these
intersections, the relationship
to matrix models and integrable systems came as a beautiful surprise.

In 1992, M.~Kontsevich provided a mathematical
connection between the intersection theory of $\overline{M}_{g,n}$
and matrix models. The connection required two main steps.
First, Kontsevich constructed a combinatorial model for the intersection
theory of $\overline{M}_{g,n}$ via
a topological stratification 
of the moduli space defined by
Strebel differentials \cite{K1}.
 The combinatorial model
expresses the tautological intersections 
as sums over trivalent
graphs on $\Sigma_g$. 
Further details of Kontsevich's construction,
some quite subtle, are discussed in \cite{Lo}. 

Kontsevich's second step
was to interpret the trivalent graph  summation 
as a Feynman diagram expansion for a new 
matrix integral (Kontsevich's matrix model). The  
KdV equations were then deduced from the analysis of the matrix
integral. The details of
this second step
are discussed in several papers, see for example \cite{dF,DFGZ,DIZ,IZ}.

Witten's conjecture (Kontsevich's theorem) is remarkable from several
perspectives --- it is certainly among
the deepest known properties of the moduli
space of curves. Once the connection to matrix models is
made, combinatorial techniques and ideas from the
theory of integrable systems may be used study the
free energy $F$ and the partition function $Z=e^F$.
For example,
Witten's conjecture
may be reformulated in terms of Virasoro constraints:
the KdV equations for $F$ are equivalent to
the annihilation of  $Z$ 
by a specific set of differential operators which form a 
representation of (a part of) the Virasoro algebra.

The moduli of stable curves $\overline{M}_{g,n}$
may be naturally viewed in
the richer context of the moduli of stable maps 
$\overline{M}_{g,n}(X)$ from curves to target varieties $X$. 
Gromov-Witten theory is the study of 
tautological intersections in $\overline{M}_{g,n}(X)$.
The development of Gromov-Witten theory was motivated
by Gromov's work on the moduli of pseudo-holomorphic maps in
symplectic geometry and Witten's study of
2 dimensional gravity \cite{Gro,W}.   
It is expected that the intersection theory
of $\overline{M}_{g,n}(X)$ will again be governed by
matrix models and their 
associated integrable hierarchies.

In particular, 
the Gromov-Witten theory of the target $X=\proj^1$
has been intensively studied by
the physicists T.~Eguchi, K.~Hori, C.-S.~Xiong,
Y.~Yamada, and S.-K.~Yang.
A conjectural formal matrix model for $\proj^1$
has led to  a precise prediction for  Gromov-Witten
theory analogous to Witten's conjecture:
intersections in $\overline{M}_{g,n}(\proj^1)$
are governed by the Toda equations  (see \cite{EgY, Ge1, P1}).

For arbitrary $X$, it remains unclear at present what the
corresponding matrix model or the 
integrable hierarchy should be. However, there 
exists a precise conjecture for the associated
Virasoro constraints. This was formulated in  1997
for an arbitrary nonsingular projective 
target variety $X$
 by Eguchi, Hori, and Xiong (using also ideas of S. Katz) \cite{EgHX}.
This Virasoro conjecture generalizes the Virasoro formulation 
of Witten's conjecture and is one of the most fundamental
open questions in Gromov-Witten theory.

\subsection{Hurwitz numbers}
The goal of the present paper is to provide a new and complete
proof of Kontsevich's combinatorial formula for 
intersections in $\overline{M}_{g,n}$. Our approach  
uses
a connection between intersections in $\overline{M}_{g,n}$
and the enumeration of branched coverings of $\proj^1$ --- Strebel
differentials play no role. In fact, two models for
the intersection theory of $\overline{M}_{g,n}$ are naturally
found from our perspective: Kontsevich's model and an alternate
model called 
the \emph{edge-of-the-spectrum matrix model}.
The relation between the latter matrix model and $\overline{M}_{g,n}$
was recognized in \cite{O2} and then used in \cite{O1}. 

Concretely, we consider the enumeration problem of Hurwitz covers
of  $\proj^1$. 
Let $\mu$ be a partition of $d$ of length $l$.
Let $H_{g,\mu}$ be  the {\em Hurwitz number}: the
number  of
genus $g$ degree $d$ covers of $\proj^1$ 
with profile $\mu$ over $\infty$ and
simple ramification over a fixed set of finite points.
The path from the intersection theory of the moduli space of curves to 
matrix models developed here uses two approaches to
the Hurwitz numbers.

First, the
numbers $H_{g,\mu}$ may be expressed in terms of
tautological intersection products in $\overline{M}_{g,l}$.
The $l$-point generating series for intersections
then arises naturally via the large $N$ asymptotics of 
$H_{g,N\mu}$.

The relationship between the numbers $H_{g,\mu}$ and
the intersection theory of $\overline{M}_{g,l}$ was independently
discovered in \cite{FanP} (for $\mu= 1^d$) and
\cite{ELSV} (for all $\mu$).
The method of \cite{FanP} is a direct 
calculation in the Gromov-Witten theory of $\proj^1$.
The Hurwitz numbers arise by definition as
intersections in $\overline{M}_g(\proj^1)$. The
virtual localization formula of \cite{GrP} precisely
relates these intersections to $\overline{M}_{g,l}$.
The study of $H_{g,\mu}$ for general $\mu$ within
the Gromov-Witten framework was completed in \cite{GrV}.
The
method of \cite{ELSV} follows a different path ---
the result is obtained by an analysis of
a twisted Segre class
construction for cones over $\overline{M}_{g,l}$.

Second, the Hurwitz numbers may be approached 
via graph enumeration. 
The large $N$ asymptotics of $H_{g,N\mu}$ is then related
to the sum over trivalent graphs arising in Kontsevich's  model.
This asymptotic analysis involves  probabilistic
techniques, in particular, a study of random 
trees is required.

\subsection{Plan of the paper}
The Hurwitz path from
the intersection theory of $\overline{M}_{g,n}$ to
matrix models 
draws motivations and techniques from several distinct
areas of mathematics. A parallel goal of the paper is
to provide an exposition of the circle of ideas
involving Gromov-Witten theory, Hurwitz numbers, 
and random graphs. 

The paper consists of three parts. The first 
part covers the background material and
explains the general strategy of the proof. 
We start with a review of
Witten's conjecture and
Kontsevich's 
combinatorial model for tautological
intersections in Section \ref{kcmod}.

The Hurwitz numbers, which are the main focus 
of the paper, are  discussed in Section \ref{shur}.
Three characterizations 
of $H_{g,\mu}$ are given in Section \ref{hur}. The relationship
between the Hurwitz numbers and the intersection theory of  moduli space
is introduced in Sections \ref{hodhur}-\ref{hurasym}.
A summary of the asymptotic study of $H_{g,\mu}$
via graph enumeration is given
in Section \ref{grasym}. 

Section \ref{smm}, concluding Part I of the paper,
 is devoted
to a brief discussion of the edge-of-the-spectrum
matrix model and Kontsevich's matrix 
model. We also
discuss there another connection between Hurwitz
numbers and integrable hierarchies via the Toda
equations. We will return to the material 
of Section \ref{smm} in the sequel  \cite{OP} to
this paper which will contain a more in depth
discussion of matrix models and integrable
hierarchies arising in Gromov-Witten
theory. 

Part II of the paper, consisting
of Sections \ref{st1}-\ref{fin1}, contains
a survey of the proof in Gromov-Witten theory of 
the formula for $H_{g,\mu}$ in the intersection theory of
$\overline{M}_{g,l}$.
Our exposition follows \cite{FanP,GrV}. 
An effort is made here to balance the geometrical ideas
with the tools needed from Gromov-Witten theory: branch morphisms,
virtual classes,
and the virtual localization formula.

In Part III of the paper, we investigate
the asymptotics of Hurwitz numbers using the methods of \cite{O2}.
Results from the theory of random trees,
summarized in Section \ref{st2}, play a significant role
in this asymptotic analysis. In the end, Kontsevich's 
combinatorial model is
precisely recovered from the asymptotics of the Hurwitz numbers. 

Finally, there are two appendices.
The classical
recursions for the Hurwitz numbers are recalled in Appendix \ref{degdeg}.
These recursive formulas are obtained by 
studying the degenerations of covers as a finite branch point is
moved to $\infty$. The degeneration formulas provide
an elementary, if not very efficient, method of computing $H_{g, \mu}$.
A short table of the 
values of the various integrals discussed in the paper
is given in Appendix \ref{tabtab}. The tables cover the cases
of $g\leq 2$ and $d\leq 4$.

\subsection{Acknowledgments}
We thank J.~Bryan,
C.~Faber,
B.~Fantechi, E.~Getzler, A.~Givental,
T.~Graber, E.~Ionel, Y.~Ruan, M.~Shapiro, R.~Vakil, and C.-S.~Xiong
for many discussions 
about Hurwitz numbers and
Gromov-Witten theory.  We thank Jim Pitman for his aid 
with the literature on random trees.
A.~O.\ was partially supported by DMS-9801466 and a Sloan foundation
fellowship. 
R.~P.\ was
partially supported by
DMS-0071473 and fellowships from the Sloan and Packard foundations.

\section{Kontsevich's combinatorial model for the intersection theory of $\mgn$}
\label{kcmod}

\subsection{Intersection theory of $\mgn$ and KdV}
\label{wcon}
The intersection theory
of $\overline{M}_{g,n}$ must be studied
in the category of Deligne-Mumford stacks (or alternatively, in
the orbifold category) to correctly handle the automorphism groups
of the pointed curves.
$\overline{M}_{g,n}$ is a complete, irreducible,
nonsingular Deligne-Mumford stack 
of complex dimension $3g-3+n$.
Intersection theory for $\overline{M}_{g,n}$ was first
developed in \cite{Mu} (see also \cite{Vi}).

We will require the tautological $\psi$ classes in
$H^2(\overline{M}_{g,n},{\mathbb{Q}})$.
For each marking $i$, there exists a 
canonical line bundle $\llr_i$ on 
$\overline{M}_{g,n}$ determined by the following prescription:
the fiber of $\llr_i$ at the stable pointed 
curve $(C,x_1, \ldots,x_n)$ is the cotangent space
$T^*_{C}(x_i)$
of $C$ at $x_i$. We note while $\llr_i$
is a {\em stack} line bundle, $\llr_i$ only determines a 
$\mathbb{Q}$-divisor on the coarse moduli space.
Let $\psi_i$ 
denote the first Chern class of
$\llr_i$.

Witten's conjecture concerns the 
complete set of evaluations of intersections of the $\psi$ classes:
\begin{equation}
\label{wittint}
\int_{\overline{M}_{g,n}} \psi_1^{k_1} \cdots \psi_n^{k_n}.
\end{equation}
The symmetric group $S_n$ acts naturally on  $\overline{M}_{g,n}$
by permuting the markings. Since the $\psi$ classes are
permuted by this $S_n$ action, the integral (\ref{wittint}) is
unchanged by a permutation of the exponents $k_i$.
A concise notation for these
intersections which exploits the
$S_n$ symmetry is given by:
\begin{equation}
\label{products}
\langle\tau_{k_1} \cdots \tau_{k_n}\rangle_g = \int
_{\overline {M}_{g,n}} \psi_1^{k_1} \cdots \psi_n^{k_n}.
\end{equation}
Such products are well defined when
the $k_i$ are non-negative integers and the dimension condition
$3g-3+n-\sum k_i=0$ holds.
In all other cases, $\langle\prod_{i=1}^{n} \tau_{k_i}\rangle_g$ is
defined to be zero. The empty product $\langle1\rangle_1$ is also 
set to zero.

The simplest integral is 
$$\langle\tau_0^3\rangle_0 = \int_{\overline{M}_{0,3}} 
\psi_1^0\psi_2^0 \psi_3^0 = 1,$$
since $\overline{M}_{0,3}$ is a point.
In fact, the genus 0 integrals are determined by the closed
form \cite{W}:
\begin{equation}
\label{genz}
\langle \tau_{k_1} \cdots \tau_{k_n} \rangle_0 =
\binom{n-3}{k_1, \ldots, k_n}.
\end{equation}
The first elliptic integral is $\langle \tau_1 \rangle_1=1/24$
which may be computed, for example, by studying a pencil
of cubic plane curves.

A fundamental property of the integrals (\ref{products}) is the
{\em string equation}: for $2g-2+n>0$,
$$\langle\tau_0 \prod_{i=1}^{n} \tau_{k_i}\rangle_g =
\sum_{j=1}^{n} \langle\tau_{k_j-1} \prod_{i\neq j} \tau_{k_i}\rangle_g.$$
Equation (\ref{genz}) easily follows from the string
equation and the evaluation $\langle \tau_0^3 \rangle_0=1$.
A second property is the {\em dilaton equation}:
for $2g-2+n>0$,
$$\langle\tau_1 \prod_{i=1}^{n} \tau_{k_i}\rangle_g =
(2g-2+n) \ \langle\prod_{i=1}^{n} \tau_{k_i}\rangle_g.$$
The string equation, dilaton equation, and the evaluation
$\langle\tau_1\rangle_1=1/24$ determine all the integrals (\ref{products})
in genus 1.

Both the string and dilaton equations are derived from a comparison
result describing the behavior of the $\psi$ classes under
pull-back via the  map $$\pi: \overline{M}_{g,n+1} \rightarrow
\overline{M}_{g,n}$$ forgetting the
last point. 
Let $i\in \{1,\ldots,n\}$. The basic formula is:
\begin{equation}
\label{pback}
\psi_i= \pi^*(\psi_i) + [D_i]
\end{equation}
where $D_i$ is the boundary divisor in $\overline{M}_{g,n}$ with
genus splitting $g+0$ and 
marking splitting $\{1,\ldots,\hat{i}, \ldots,n\} \cup \{i, n+1\}$.
That is, the general point of $D_i$ corresponds
to a reducible curve $C= C_1 \scup C_2$ connected by a single node
satisfying:
\begin{enumerate}
\item[(i)]
$C_1$ is nonsingular of  genus $g$
\item[(ii)]
$C_2$ is nonsingular of genus 0.
\item[(iii)] The markings $\{1,...,n\} \setminus \{i \}$
 lie on $C_1$ and the remaining marking $\{i,n+1\}$
lie on $C_2$.
\end{enumerate}
The relation (\ref{pback}) implies the string
and dilaton equations by a direct geometric argument (see, for example,
 \cite{W}).

The KdV equations are differential equations satisfied
by a generating series of the $\psi$ intersections.
Let $t$  denote the set of variables $\{t_i\}_{i=0}^\infty$. 
Let $\gamma=\sum_{i=0}^{\infty} t_i \tau_i$ be the formal sum.
Consider the formal
generating function for the integrals (\ref{products}):
$$F_g(t)= \sum_{n=0}^{\infty} 
\frac{\langle\gamma^n\rangle_g}{n!}.$$
The expression $\langle\gamma^n\rangle_g$  is defined by monomial expansion
and multilinearity in the variables $t_i$. More concretely,
\begin{eqnarray*}
F_g(t) & = & \sum_{n\geq 1} \frac{1}{n!} \sum_{k_1,\ldots,k_n} \langle
\tau_{k_1} \cdots \tau_{k_n} \rangle_g  \ t_{k_1} \cdots t_{k_n} \\ 
& = & \sum_{\{n_i\}} 
\langle
\tau_0^{n_0} \tau_1^{n_1} \tau_2^{n_2} \cdots\rangle_g \prod_{i=0}^{\infty} 
\frac{t_i^{n_i}}{n_i!},
\end{eqnarray*}
where the last sum is over all sequences of nonnegative integers $\{n_i\}$
with finitely many nonzero terms. 
Let $F$ denote the full generating function:
$$F= \sum_{g=0}^{\infty} F_g.$$
The genus subscript $g$ of a non-vanishing bracket
 $\langle \tau_{k_1} \ldots \tau_{k_n} \rangle_g$
is determined by the dimension condition 
$3g-3+n -\sum_{i=1}^n k_i =0.$
Hence, $F$ is a faithful generating series of all the $\psi$ intersections
in $\overline{M}_{g,n}$.

We will use the following notation for the derivatives of $F$:
\begin{equation}
\label{pproductss}
\langle\langle\tau_{k_1} \tau_{k_2} \cdots \tau_{k_n}\rangle\rangle = 
\frac{\p}{\p t_{k_1}} \frac{\p}{\p t_{k_1}} \cdots \frac{\p}{\p t_{k_1}} F.
\end{equation}
Note $\langle\langle\tau_{k_1} \tau_{k_2} 
\cdots \tau_{k_n}\rangle\rangle|_{t_i=0}= \langle\tau_{k_1}
\tau_{k_2} \cdots \tau_{k_n}\rangle$.

$F$
was conjectured by Witten to equal
the free energy in $2$ dimensional quantum gravity
and therefore to satisfy the KdV hierarchy.
The classical 
KdV equation (first studied in the $19^{th}$ century
to describe shallow water waves) is:
\begin{equation}
\label{kkdv}
\frac{\p U}{\p t_1}= U 
\frac{\p U}{\p t_0}
+ \frac{1}{12} \frac{\p^3 U}{\p t^3_0}.
\end{equation}
Witten conjectured $U=\langle \langle t_0 t_0 \rangle \rangle$
satisfies (\ref{kkdv}). The KdV hierarchy for $F$ may be written
in the following form (equation (\ref{kkdv}) is recovered
in case $n=1$).

\vspace{+10pt}
\noindent{\bf Witten's Conjecture.} For all
$n\geq 1$,
\begin{equation}
\label{kdv}
(2n+1)\langle\langle\tau_n \tau_0^2 \rangle\rangle =
\end{equation}
$$\langle\langle\tau_{n-1} \tau_0\rangle\rangle
\langle\langle\tau_0^3\rangle\rangle+ 
2\langle\langle
\tau_{n-1}\tau_0^2\rangle\rangle\langle\langle\tau_0^2\rangle\rangle+
\frac{1}{4}\langle\langle\tau_{n-1} \tau_0^4\rangle\rangle.$$
\vspace{+10pt}

As an example, consider equation (\ref{kdv}) for $n=3$ evaluated at
$t_i=0$. We obtain:
$$7\langle\tau_3 \tau_0^2\rangle_1= 
\langle\tau_2 
\tau_0\rangle_1 \langle\tau_0^3\rangle_0 +{\frac{1}{4}} \langle\tau_2
\tau_0^4\rangle_0.$$
Use of the string equation yields:
$$7\langle\tau_1\rangle_1= 
\langle\tau_1\rangle_1 \langle\tau_0^3\rangle_0
 + \frac{1}{4} \langle\tau_0^3\rangle_0.$$ 
Hence, we conclude $\langle\tau_1\rangle_1=1/24$.
The KdV equations (\ref{kdv}) and 
the string equation together
determine all the integrals (\ref{products}) from 
$\langle \tau_0^3 \rangle_0 =1$.
Therefore, $F$ is 
uniquely determined by Witten's conjecture.

\subsection{Kontsevich's combinatorial model}
\label{kmod}
We now explain the model found by Kontsevich
for the generating series:
\begin{equation}
\label{kgen}
K_{g}(s_1, \ldots, s_n) =
\sum_{\sum k_i =3g-3+n}
\langle \tau_{k_1} \cdots \tau_{k_n} \rangle_{g}
\prod_{i=1}^n \frac{(2k_i-1)!!}{s_i^{2k_i+1}}
\end{equation}
of $\psi$ intersections in $\overline{M}_{g,n}$.

Let $\Sigma_g$ be an oriented topological surface of genus $g$.
A \emph{map} $G$ on $\Sigma_g$ is a triple $(V,E,\phi)$ 
satisfying the following conditions:
\begin{enumerate}
\item[(i)]
$V \subset \Sigma_g$ is a finite set of vertices,
\item[(ii)]
$E$ is finite set of edges:  
\begin{enumerate}
\item[$\bullet$] each edge is a simple path in $\Sigma_g$
connecting two vertices of $V$,
\item[$\bullet$] self-edges at vertices are permitted,
\item[$\bullet$] distinct edge paths intersect only in vertices,
\end{enumerate}
\item[(iii)] the graph $G$ is connected,
\item[(iv)] the complement of the union of the edges in $\Sigma_g$
is a disjoint union of topological disks, called the \emph{cells}
of $G$,
\item[(v)] $\phi$ is a bijection of the set $\Cell(G)$ of cells
with $\{1,\ldots, |\Cell(G)|\}$.
\end{enumerate}

The origin of the term ``map'' is the following: one can visualize
the cells of a map $G$ as different countries into which $G$ 
divides the surface $\Sigma_g$.  

The {\em valence} of a vertex $v$ is given by the number of half-edges
incident to $v$. A map $G$ is called {\em trivalent} if
every vertex has valence exactly 3. The map
$G$ is called {\em stable} if 
$$
2g-2+|\Cell(G)|>0\,.
$$
Two maps  $G$ and $G'$ on $\Sigma_g$ are isomorphic if there is
an orientation preserving homeomorphism of $\Sigma_g$
which maps $G$ to $G'$ and respects $\phi$.
The automorphism group $\text{Aut}(G)$
is the finite group of symmetries of $(V,E,\phi)$
induced by orientation preserving homeomorphisms of $\Sigma_g$
that map $G$ to $G$
and respect the marking $\phi$.

\begin{figure}[!hbt]
\centering
\scalebox{.7}{\includegraphics{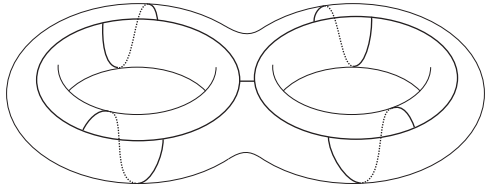}}
\caption{A trivalent map on a genus $2$ surface}
\label{fig0}
\end{figure}

Let $\bG_{g,n}$ denote the set of isomorphism classes
of maps on $\Sigma_g$ with $n$ cells and let
$\bGt_{g,n}\subset \bG_{g,n}$ denote the subset of
trivalent maps.
The trivalent condition and
the Euler characteristic constraint on $G\in \bGt_{g,n}$
imply:
\begin{gather}
  |V|= \frac{2}{3} |E|\,, \\
|V|= 2(2g-2+n) \,,
\end{gather}
$|E|$, $|V|$ denote the cardinality of $E$ and $V$
respectively.
It is then easy to see that $\bGt_{g,n}$ is a finite set.
An example of an element of $\bGt_{2,3}$ is shown in
Figure \ref{fig0}.

Let $g\geq 0$ and $n$ be fixed in the stable range
$2g-2+n>0$.
Let the variables $s_1, \ldots, s_n$  
correspond to the markings of $G\in \bGt_{g,n}$.
Each edge $e\in E$ of $G$ borders two cells. Let
$i$ and $j$ be the labels assigned by $\phi$  to these
cells. If both sides of $e$ border the {\em same} cell, then $i=j$.
We denote
$$
\widetilde{s}(e)= s_i+ s_j \,.
$$

The fundamental result proven by Kontsevich is
the following formula for $K_{g}$ in terms
of combinatorics of trivalent maps:

\begin{tm} \label{kon}
$K_{g}$ is obtained by a sum over trivalent maps:
\begin{equation}
\label{dfgg}
K_{g}(s_1,\ldots,s_n) = \sum_{G\in \bGt_{g,n}}
\frac{2^{2g-2+n}}{|\text{\em Aut}(G)|} \prod_{e\in E} \frac{1}{\widetilde{s}(e)}.
\end{equation}
\end{tm}

Kontsevich's proof requires a topological decomposition
of $\overline{M}_{g,n}$ obtained via the theory of
Strebel differentials (see \cite{K1}, Appendix B). 
Aspects of the boundary behavior of this geometry are quite subtle.
A discussion
can also be found in \cite{Lo}.

\section{Hurwitz numbers}\label{shur}

\subsection{Three definitions of Hurwitz numbers}
\label{hur}
Three equivalent definitions of the Hurwitz numbers are 
discussed
in this section.
Definitions \ref{Hur1} and \ref{Hur2}  will be used to provide a new proof
of Theorem 1 connecting $\psi$ intersections in
$\overline{M}_{g,n}$ to Kontsevich's combinatorial model.
Definition \ref{Hur3} relates the Hurwitz numbers to the
combinatorics of the symmetric group and arises in the
connection between  Hurwitz
numbers and  the Toda equations in the  Gromov-Witten theory of $\proj^1$.

\subsubsection{Enumeration of branched coverings} 

We start with the definition of the Hurwitz numbers $H_{g,\mu}$
via covers of $\proj^1$.
Let $g\geq 0$ and let $\mu$ be a non-empty partition.
 Let $|\mu|$ denote the
sum of the parts of $\mu$, and let $\ell(\mu)$
denote the length of $\mu$.
A Hurwitz cover of $\proj^1$ of
 genus $g$ and ramification $\mu$ over $\infty$ is
a morphism
$$\pi: C \rarr \proj^1$$ 
satisfying the following properties:
\begin{enumerate}
\item[(i)] $C$ is a nonsingular, irreducible
 genus $g$ curve,
\item[(ii)] the divisor $\pi^{-1} (\infty)\subset C$ has profile equal
            to the partition $\mu$,
\item[(iii)] the map $\pi$ is simply ramified over 
${\mathbf{A}}^1=\proj^1 \setminus
\infty$.
\end{enumerate}
Note that condition (ii) implies 
$$
\deg\pi= |\mu|\,.
$$

By the Riemann-Hurwitz formula,
the number of simple ramification points of $\pi$ over  ${\mathbf{A}}^1$
is: $$r(g,\mu)=2g-2+|\mu|+ \ell(\mu).$$
Let $U_r$ denote a fixed set of $r=r(g,\mu)$ distinct points in 
${\mathbf{A}}^1$ --- it will be convenient for us to
take $U_r$ equal to the set of $r^{th}$ roots of unity
in $\com = {\mathbf{A}}^1$. We will require the simple ramification
points of $\pi$
to lie over $U_r$.

Two covers 
$$\pi: C \rarr \proj^1, \ \pi': C' \rarr \proj^1$$ are isomorphic if
there exits an isomorphism of curves $\phi: C \rarr C'$ satisfying
$\pi'\circ \phi= \pi$. Each  cover $\pi$ has an naturally
associated automorphism
group $\text{Aut}(\pi)$.

\begin{definition}\label{Hur1} 
$H_{g,\mu}$ is a weighted count of the distinct
Hurwitz covers $\pi$ of genus $g$ with ramification $\mu$
over $\infty$ and simple ramification over $U_r$.
Each such cover is weighted by $1/|\text{Aut}(\pi)|$.
\end{definition}

\subsubsection{Enumeration of branching graphs} 

The Hurwitz numbers admit a second definition via 
enumeration of graphs, see for example \cite{Ar}.
 Let $g\geq 0$ and $\mu$ be fixed. Let 
$r=r(g,\mu)$
and $U_r=\{\zeta_1, \ldots, \zeta_r\}$ be the
set of $r^{th}$ roots of unity as above.

A \emph{branching graph} $H$ on an oriented topological surface
$\Sigma_g$ consists of the data
$(V,E,\gamma: E\rarr U_r)$ satisfying the following 
conditions:
\begin{enumerate}
\item[(i)]
the vertex set $V \subset \Sigma_g$ consists of $|\mu|$ distinct points,
\item[(ii)]
the set $E$ consists of $r$ edges:
\begin{enumerate}
\item[$\bullet$] each edge is a simple path in $\Sigma_g$
connecting two vertices of $V$,
\item[$\bullet$] self-edges at vertices are not permitted,
\item[$\bullet$] distinct edge paths intersect only in vertices,
\end{enumerate}
\item[(iii)] the graph $H$ is connected,
\item[(iv)] the  function
$\gamma$ is a bijection,
\item[(v)] at each vertex $v\in V$, the cyclic order
of the edge markings (with respect to the orientation of
$\Sigma_g$) agrees with the cyclic order of the roots of unity
(with respect to the standard orientation of $\com$),
\item[(vi)] the complement of the union of the edges is a
disjoint union of $l=\ell(\mu)$ topological disks
$D_1,\ldots,D_l$.
\end{enumerate}
Let $D_i$ be a cell bounded by the sequence of
edges $e_{12},\ldots, e_{s1}$ 
of the graph $H$. Assume the edge circuit is {\em clockwise} 
with respect to the orientation of $D_i$ restricted
from $\Sigma_g$. Then, to each pair of edges
$e_{k-1,k}$, $e_{k,k+1}$ there is an associated positive 
angle given by:
$$
\measuredangle(e_{k-1,k},e_{k,k+1})=
\arg\left( \frac{\gamma(e_{k-1,k})}{\gamma(e_{k,k+1})} \right) \in (0,2\pi].
$$
The sum of these angles along the boundary of $D_i$ 
is a multiple of $2\pi$. In other words, the 
following {\it perimeter} of the cell $D_i$
$$
\per(D_i) = \frac1{2\pi} \sum_{k=1}^s \measuredangle(e_{k-1,k},e_{k,k+1})
$$
is a positive integer.

The cyclic ordering condition (v) implies
that $\sum_{i=i}^l \per(D_i) = |\mu|$.
The last condition in the definition of a branching 
graph is:

\begin{enumerate}
\item[(vii)] The partition $\mu$ equals $(\per(D_1),\ldots,\per(D_l))$.
\end{enumerate}

Two 
branching graphs $H$ and $H'$ on $\Sigma_g$ 
are isomorphic if there exists an orientation
preserving homeomorphism of $\Sigma_g$ which
maps $H$ to $H'$ and respects the edge markings.
The automorphism group $\text{Aut}(H)$
is the finite group of symmetries of $(V,E)$
induced by orientation preserving homeomorphisms of $\Sigma_g$
which map $H$ to $H$
and respect the edge markings.

Let $\bH_{g,\mu}$ denote the set of isomorphism 
classes of genus $g$ branching graphs with
perimeter $\mu$. 
The second definition of the Hurwitz numbers is
by an enumeration of graphs:

\begin{definition}\label{Hur2} 
$H_{g,\mu}$ equals a weighted count of
the branching graphs $H$ in $\bH_{g,\mu}$, where
each graph $H$ is weighted by $1/|\text{Aut}(H)|$.
\end{definition}

Definition \ref{Hur2} can be seen to agree with Definition \ref{Hur1} by
a direct association of a branching graph to
each Hurwitz cover with ramification $\mu$.
Let $\pi:C \rarr \proj^1$ be a Hurwitz cover of genus $g$
with ramification $\mu$ over infinity and simple
ramification over $U_r$.
First, 
observe that $\pi$ is unramified over the open unit
disk at the origin:
$$
B\subset\com = {\mathbf{A}}^1\,.
$$
Therefore, $\pi^{-1}(B)$ is the disjoint union
of exactly $|\mu|$ open disks $$B_1, \ldots, B_{|\mu|} \subset C.$$
Let $\overline{B}_i$ and $\p B_i= \overline{B}_i \setminus B_i$
denote the closure and the boundary of $B_i$ respectively. 

Let $q$ be an intersection point of two
different closed disks $\overline{B}_i$ and $\overline{B}_j$.
Then $q$ must be a ramification point of $\pi$ and hence
$\pi(q)\in U_r$.
In fact, as $\pi$ is simply ramified over $U_r$,
every element $\zeta\in U_r$ must lie under exactly
one intersection of different closed disks.
Therefore, there are exactly $r$  intersection points
of pairs of closed disks $Q=\{q_1, \ldots, q_r\}$, 
in bijective correspondence $\pi$ with the set $U_r$. 

Define a branching graph $H=(V,E,\gamma: E \rarr U_r)$
 on the Riemann surface $C$
by the following data:
\begin{enumerate}
\item[(a)] $V=\pi^{-1}(0)$,
\item[(b)] the edge set  $E$ corresponds to the
intersection set $Q$,
\item[(c)] the function $\gamma: E \rarr U_r$
is defined by the projection $\pi:Q \rarr U_r$.
\end{enumerate}
The edges $E$ are constructed as follows. Suppose
$$
q = \overline{B}_i \scap \overline{B}_{j} 
$$
and $\zeta=\pi(q)$. Let $[0,\zeta]$ be the
segment connecting $0$ to $\zeta$ in $\mathbf{A}^1$. 
The edge associated to $q$ is defined to be
the unique component of $\pi^{-1}([0,\zeta])$ that
connects the centers of $B_i$ and $B_{j}$.
The required
conditions (i)-(vii) of a branching graph are
easily checked. 

Conversely, every branching graph on $\Sigma_g$ with
perimeter $\mu$ 
corresponds to  a Hurwitz cover with ramification $\mu$
which can be obtained 
by reversing the above construction.
The automorphism groups
of the Hurwitz cover and of the branching graph
coincide under this identification.
We therefore conclude that Definitions \ref{Hur1} and \ref{Hur2} agree.

\begin{figure}[!hbt]
\centering
\scalebox{.7}{\includegraphics{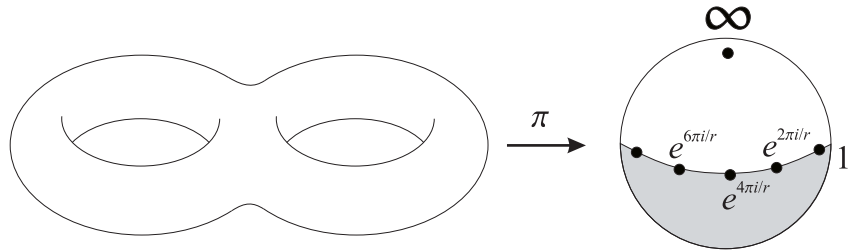}}
\caption{A covering $\pi$ ramified over $\infty$ and roots of unity}
\label{fig1}
\end{figure}

\begin{figure}[!hbt]
\centering
\scalebox{.5}{\includegraphics{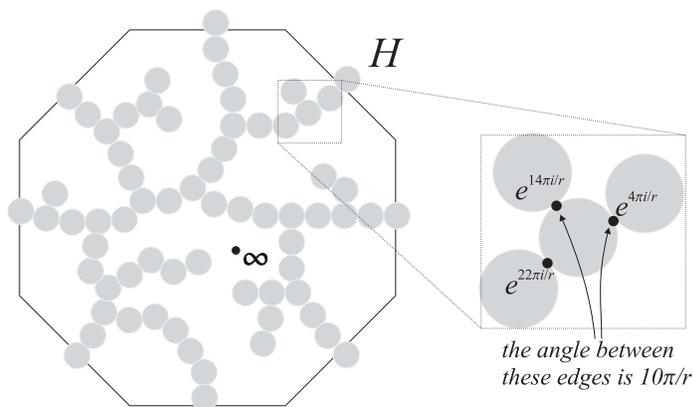}}
\caption{Preimage on $\Sigma_2$ of the unit circle under the map $\pi$}
\label{fig2}
\end{figure}

Figures \ref{fig1} and \ref{fig2} should help visualize the 
relationship between Definitions  \ref{Hur1} and \ref{Hur2}.
Suppose we have a covering $\pi$ of $\proj^1$ which satisfies the
conditions of  Definition  \ref{Hur1}, such as the one shown 
schematically in Figure \ref{fig1}.  

In Figure \ref{fig2}, we see
the preimage of the unit circle $B$ under $\pi$ consists
of $\deg\pi$ disks which meet at the ramification points of $\pi$. Such
points correspond bijectively under $\pi$ to the roots of unity. The centers of
the disks
form the  vertices of the branching graph $H$, and the intersection points of the
disks correspond to the edges of $H$. Since the edges of $H$ are
labeled by roots of unity, we can define the angle between two
edges and then the perimeters of the cells of $H$. In Figure \ref{fig2},
most edge labels of $H$ are omitted except on a small part of $H$ which
is magnified. 

\subsubsection{Counting factorizations into transpositions}

A third approach to the
Hurwitz numbers
via the combinatorics of the symmetric
group $S_{|\mu|}$ also plays a role in Gromov-Witten theory. 
A Hurwitz cover of genus $g$ with ramification $\mu$
over $\infty$ and simple ramification over $U_r$
can be associated to an ordered sequence of
transpositions $(\gamma_1, \ldots, \gamma_r)$ of $S_{|\mu|}$ 
satisfying the following two
properties:
\begin{enumerate}
\item[(a)] $\gamma_1, \ldots, \gamma_r$ generate $S_{|\mu|}$,
\item[(b)] the product  $\gamma_1 \gamma_2 \cdots \gamma_r$
has cycle structure $\mu$.
\end{enumerate}
The associated Hurwitz cover is found by the following
topological construction.

The fundamental group $\pi_1({\mathbf{A}^1} \setminus U_r)$
is freely generated by the loops around the points $U_r$.
Let $$\widetilde{ {\mathbf{A}^1} \setminus U_r}$$ denote the
universal cover of ${\mathbf{A}^1} \setminus U_r$.
The sequence $(\gamma_1, \ldots, \gamma_r)$
defines an action of $\pi_1({\mathbf{A}^1} \setminus U_r)$
on $\{1,2,\dots,|\mu|\}$. This determines  an unramified,
 $|\mu|$-sheeted
covering space 
$$
\pi^0: C^0\rarr {\mathbf{A}^1} \setminus U_r
$$
defined by the mixing construction:
$$
C^0 = \widetilde{ {\mathbf{A}^1} \setminus U_r} 
\times_{\pi_1({\mathbf{A}^1} \setminus U_r)} \{1,2,\ldots 
|\mu|\}.
$$
The covering $C^0$ is connected by condition (a).
$C^0$ is naturally endowed with a complex structure and
may be canonically completed to yield a Hurwitz cover $\pi:C \rarr \proj^1$
of genus $g$ and ramification $\mu$ by condition (b).
All Hurwitz covers of genus $g$ with ramification $\mu$ over
$\infty$ and simple ramification over $U_r$
arise in this way. Therefore, the following
definition of the Hurwitz numbers is equivalent to Definition \ref{Hur1} 

\begin{definition}\label{Hur3}
$H_{g,\mu}$ equals $1/|\mu|!$ times the number of $r$-tuples
of 2-cycles satisfying (a) and (b).  
\end{definition}

Formulas for $H_{g,\mu}$ in terms of
the characters of the symmetric group 
were deduced by Burnside from this perspective.
In fact, Hurwitz's original computations of 
covering numbers were obtained via symmetric group
calculations \cite{Hu}.

\subsection{Hurwitz numbers and the intersection theory of
$\overline{M}_{g,n}$}
\label{hodhur}
The Hurwitz numbers  are naturally expressed
in terms of tautological intersections in $\overline{M}_{g,n}$.
However, we will require here not only the $\psi$ classes
arising in Witten's conjecture, but also the $\lambda$
classes. 
Let the {\em Hodge bundle} $$\hodge \rarr \overline{M}_{g,n}$$
be the rank $g$ vector bundle with fiber
$H^0(C,\omega_C)$ over the moduli point $(C,p_1, \ldots,p_n)$.
The  $\lambda$ classes are the Chern classes of the Hodge bundle:
$$\lambda_i= c_i(\hodge) \in H^{2i}(\overline{M}_{g,n},
\mathbb{Q}).$$
The $\psi$ and $\lambda$ classes are {\em tautological} classes
on the moduli space of curves. 
A foundational treatment
of the tautological intersection
theory of $\overline{M}_{g,n}$ can be found in \cite{Mu} (see \cite{Fa,FaP3}
for a current perspective).

Let $\mu=(\mu_1,\ldots,\mu_l)$ be a non-empty partition with positive
parts. Let $\text{Aut}(\mu)$ denote the permutation
group of symmetries of the parts of $\mu$. 
The Hurwitz numbers $H_{g,\mu}$ are related to
the intersection theory of $\overline{M}_{g,l}$
by the following formula.

\begin{tm}
\label{rrrr} 
Let $2g-2+\ell(\mu)>0$. 
The Hurwitz number $H_{g,\mu}$ satisfies:
\begin{equation}
\label{hodgehurq}
H_{g,\mu}= 
\frac{(2g-2+|\mu|+l)!}{|\text{\em Aut}(\mu)|} 
\prod_{i=1}^l \frac{\mu_i^{\mu_i}}{\mu_i!}
\int_{\overline{M}_{g,l}}
\frac{\sum_{k=0}^g (-1)^k \lambda_k} 
{\prod_{i=1}^l (1-\mu_i \psi_i)}.
\end{equation}
\end{tm}

Theorem \ref{rrrr} was proven by 
T. Ekedahl, Lando, M. Shapiro, and
Vainshtein \cite{ELSV} using a theory
of twisted Segre classes for cone bundles
over $\overline{M}_{g,n}$.
In case $\mu= 1^d$, the case
of trivial ramification over $\infty$, formula (\ref{hodgehurq})
was independently found and proven in \cite{FanP} via 
a direct integration in Gromov-Witten theory.
This approach was later 
refined in \cite{GrV} to yield the formula (\ref{hodgehurq})
for the general partition $\mu$.

The proof of \cite{FanP} begins with an 
integral formula in Gromov-Witten theory for the
Hurwitz numbers.
Let $\overline{M}_{g}(\proj^1,d)$
be the moduli space of stable maps of genus $g$ and degree $d$
to $\proj^1$.
There is branch morphism:
$$ \br:\overline{M}_{g}(\proj^1,d) 
\rarr \text{Sym}^{2g-2+2d}(\proj^1)$$
which assigns to each stable maps $f:C \rarr \proj^1$
the branch divisor in the target \cite{FanP}.
Using Definition 1 of the Hurwitz numbers and properties of
the virtual class, an integral formula 
\begin{equation}
\label{hhhh}
H_{g,1^d} = \int_{[\overline{M}_{g}(\proj^1,d)]^{vir}} \br^*(\xi_p)
\end{equation}
may be obtained.
Here, $\xi_p$ is (the Poincar\'e dual) of the point class of
$\text{Sym}^{2g-2+2d}(\proj^1)$.

Integrals in Gromov-Witten are evaluated
against the virtual fundamental class of the moduli space
of maps $[\overline{M}_{g}(\proj^1,d)]^{vir}$.
The moduli space of maps itself may be quite ill-behaved
as all possible stable maps occur --- including
maps with reducible domains, collapsed components,
and maps defined by special linear series. In general,
$\overline{M}_g(\proj^1,d)$ is reducible and of impure
dimension. However, Gromov-Witten theory is based
on the remarkably uniform behavior of the virtual
class. Integrals against the virtual class are
{\em easier} to understand than general intersections
in the moduli space of maps.

The virtual localization formula of \cite{GrP}
provides a direct approach to the
integral in (\ref{hhhh}).
The moduli space $\overline{M}_{g}(\proj^1,d)$ has
a natural $\com^*$-action induced by the
standard $\com^*$-action on $\proj^1$.
By construction, $\br^*(\xi_p)$ is seen to be
an $\com^*$-equivariant class.
The $\com^*$-fixed loci in $\overline{M}_g(\proj^1,d)$
are well-known to be products of moduli spaces
of pointed curves \cite{K2,GrP}. The localization formula
then precisely relates equivariant integrals against 
$[\overline{M}_{g}(\proj^1,d)]^{vir}$ to 
tautological intersections in the moduli space
of pointed curves. 
Formula (\ref{hodgehurq}) for $\mu = 1^d$
is the result.

In case $\mu$ is arbitrary, the above strategy may be
followed on an appropriate  component  
of the moduli space  $$\overline{M}_g(\proj^1,d(\mu)=|\mu|)$$
via an elegant localization analysis provided in \cite{GrV}.

Sections \ref{st1}-\ref{vll} contains a review of the 
Gromov-Witten theory of $\proj^1$ and the virtual
localization formula.
The proof of Theorem \ref{rrrr} is presented
in Section \ref{fin1} following \cite{FanP,GrV}.

\subsection{Asymptotics of the Hurwitz numbers I: $\psi$
integrals}
\label{hurasym}
Let $\mu$ be a partition with $l$ parts $\mu_1, \ldots,\mu_l$ (assumed here to
be {\em distinct}).
Let $N\mu$ 
denote the partition obtained by scaling
each part of $\mu$ by $N$.
The asymptotics of 
$H_{g,N\mu}$
as $N \rarr \infty$ are easily related
to the $l$-point function in $2$-dimensional quantum gravity
by Theorem \ref{rrrr}. After a Laplace transform,
Kontsevich's series (\ref{kgen}) is found.

The $l$-point function $P_{g}$ is defined
by the following equation (for $2g-2+l>0$):
\begin{equation}
  \label{P_g}
 P_{g}(x_1,\ldots,x_l)= \sum_
{\sum_i k_i=3g-3+l}
\langle\tau_{k_1} \cdots \tau_{k_l}\rangle_{g} 
\prod_{i=1}^l x_i^{k_i}\,. 
\end{equation}
The $l$-point function
$P_{g}$ contains the data of the full set of 
$\psi$ integrals on $\overline{M}_{g,l}$.

Define the function $H_{g}(\mu_1,\ldots,\mu_l)$ as the following limit:
\begin{equation}
\label{limmm}
H_{g}(\mu_1,\ldots,\mu_l) =
\lim_{N\rarr \infty}
\, \frac{1}{N^{3g-3+l/2}}\,
\frac{H_{g,N\mu}} {e^{N|\mu|}\ 
r(g,N\mu)!},
\end{equation}
A direct application of Theorem 2 together with Stirling's formula 
\eqref{Stirl} 
then yields the following
result governing the asymptotics of the Hurwitz numbers.
\begin{pr}
\label{lhur}
 We have:
$$
H_{g}(\mu_1,\ldots,\mu_l)= 
\frac{1}{(2\pi)^{l/2}} \frac{1}{\prod_{i=1}^l \mu_i^{1/2}} \ P_{g}(\mu_1,\ldots,\mu_l)\,.
$$
\end{pr}

Let $\mu$ be a vector with 
distinct, positive, {\em rational}
parts. The asymptotics of $H_{g,N\mu}$ 
are then well-defined over sufficiently divisible $N$, and
Proposition \ref{lhur} remains valid. It is natural to define
$H_g(x_1, \ldots,x_l)$ for all positive real values by
Proposition \ref{lhur}.

Let $LH_{g}$ denote the Laplace transform of the function
$H_{g}$:
\begin{eqnarray*}
LH_{g}(y_1,\ldots,y_l) & = & \int_{x\in {\mathbb{R}}^l_{> 0}} 
e^{-y\cdot x} \frac{1}{(2\pi)^{l/2}} \frac{1}{\prod_{i=1}^l x_i^{1/2}}\
 P_{g}(x)\ dx \\
& = & \sum_{\sum k_i={3g-3+l}}
\langle\tau_{k_1} \cdots \tau_{k_l}\rangle_{g} 
\prod_{i=1}^l \frac{(2k_i-1)!!}{(2y_i)^{k_i+\frac{1}{2}}}
\end{eqnarray*}
The variable substitution $s_i= \sqrt{2y_i}$ relates
the answer to Kontsevich's model.

\begin{tm}
\label{connn}
The Laplace transform of $H_{g}$ in the variables $s_i$
equals Kontsevich's generating series for $\psi$ integrals:
$$
LH_{g}(y_1,\ldots,y_l) = \sum_{\sum k_i={3g-3+l}}
\langle\tau_{k_1} \cdots \tau_{k_l}\rangle_{g} 
\prod_{i=1}^l \frac{(2k_i-1)!!}{s_i^{2k_i+1}}\,, \quad s_i= \sqrt{2y_i} \,.
$$
\end{tm}

We have completed the path from Hurwitz numbers to $\psi$ integrals
via Definition 1 and
Gromov-Witten theory. The result after taking the appropriate
asymptotics and
the Laplace transform is Kontsevich's series (\ref{kgen}).

\subsection{Asymptotics of the Hurwitz numbers II: graph enumeration}
\label{grasym}
Let $\mu$ be a partition with $l$ distinct parts as above.
The asymptotics of the Hurwitz numbers $H_{g,N\mu}$
may be studied alternatively via
Definition \ref{Hur2} and an analysis of graphs.
The result after Laplace transform exactly equals Kontsevich's
sum over trivalent maps on $\Sigma_g$ (\ref{dfgg}).
The two approaches to the asymptotics of the Hurwitz numbers
together yield a new proof of Theorem \ref{kon}.

Let $\bGth_{g,n} \subset \bG_{g,n}$ denote the subset of 
maps with at least trivalent vertices. In case $2g-2+n>0$, there
exists  a natural map
$$
\hmt:   \bG_{g,n} \to \bGth_{g,n}
$$
which we call the \emph{homotopy type map}. It is constructed
as follows. 

First, given a map $G\in \bG_{g,n}$ one repeatedly removes
all univalent vertices  from $G$ together with the incident edges
until there are no more univalent vertices. After that, one
removes all $2$-valent vertices by concatenating their incident
edges. The resulting map is, by definition, $\hmt(G)$. It is
clear that 
$$
|\Cell(G)|= |\Cell(\hmt(G))| \,.
$$
By definition, two maps $G$ and $G'$ on $\Sigma_g$
 have the same homotopy type
if $\hmt(G)=\hmt(G')$. 

\begin{figure}[!hbt]
\centering
\scalebox{.5}{\includegraphics{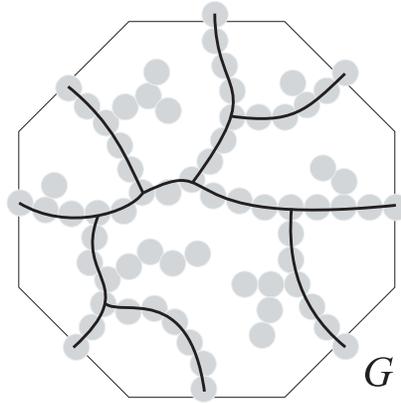}}
\caption{The homotopy type of the graph $H$ from Figure \ref{fig2}}
\label{fig3}
\end{figure}

In case the parts of $\mu$ are distinct,
there is a natural mapping 
$$\und: \bH_{g,\mu} \rarr \bG_{g, \ell(\mu)} $$
from branching graphs to underlying maps which
forgets the edge labels. The composition of $\und$ and $\hmt$
defines homotopy type
and  homotopy
equivalence for branching graphs in $\bH_{g,\mu}$. For example,
the homotopy type $G$ corresponding to the
branching graph $H$ from Figure \ref{fig2} is
shown in Figure \ref{fig3}.
Kontsevich's combinatorial model is naturally found from
the asymptotic enumeration of branching graphs by their homotopy type.

For any $G\in\bGth_{g,l}$, let $H_{G,\mu}$ denote the (weighted)
number of branching graphs $H$ on $\Sigma_g$ of
homotopy type $G$. 
By Definition \eqref{Hur2} of the Hurwitz numbers,
$$
H_{g,N\mu} = \sum_{G\in \bGth_{g,l}}  
H_{G,N\mu},
$$
Since $\bGth_{g,l}$ is a finite set, we have 
\begin{equation}
\label{limmmv}
H_{g}(\mu_1,\ldots,\mu_l) =
\sum_{G\in \bGth_{g,l}}  \lim_{N\to\infty} \,
\frac{1}{N^{3g-3+l/2}} \, 
\frac{ H_{G,N\mu}}{e^{N|\mu|}\ 
r(g,N\mu)!} \,. 
\end{equation}
The contribution of $G$ to (\ref{limmmv}) is determined
by an asymptotic analysis in Section \ref{fin2}.
If $G$ is not trivalent, the contribution vanishes.
For trivalent graphs, the contribution of $G$ to (\ref{limmmv})
is found to equal, after the Laplace transform,
the contribution of $G$ to \eqref{dfgg}. As a consequence,
we obtain the following result: 

\begin{tm}
\label{grasy}
The Laplace transform of $H_{g}$ 
in the variables $s_i$ equals a sum over trivalent graphs:
$$
LH_{g}(y_1,\ldots,y_l) = \sum_{G\in \bGt_{g,l}}
\frac{2^{2g-2+l}}{|\text{\em Aut}(G)|} \prod_{e\in E} \frac{1}{\widetilde{s}(e)}\,,
\quad \quad s_i= \sqrt{2y_i} \,.$$
\end{tm}
\noindent
Theorems \ref{connn} and \ref{grasy} together provide
a new proof of Theorem \ref{kon}.

The analysis of Section \ref{fin2} is based on the study of trees
undertaken in Section \ref{st2}. The (multivalued) inverse
of the homotopy
type map may be viewed as generating
trees over the edges of $G$.
The large $N$ asymptotics of $H_{G,N\mu}$ is  
thus governed by the theory of random trees.

\section{Matrix models and integrable hierarchies}\label{smm} 

In this section we indicate several connections
between the material in this paper and the theory of matrix models
and integrable hierarchies. A more detailed treatment of this fundamentally
important subject will be given in the sequel to this paper \cite{OP}.
Some references to existing literature are given below. 

\subsection{Edge-of-the-spectrum matrix model} 
\subsubsection{Wick's formula} 
Consider the linear space of all $N\times N$ Hermitian
matrices and the Gaussian measure on it with 
density $\displaystyle e^{-\frac12\tr M^2}$. The 
expectations with respect to this measure will be denoted
by
\begin{equation}
  \label{gaussme}
\la f \ra_{N}  =  \frac{\displaystyle \int 
f(M) \, \exp\left(-\tr M^2/2\right) \, dM}
{\displaystyle \int 
\exp\left(-\tr M^2/2\right) \, dM} \,.  
\end{equation}
It is clear that this measure has mean zero and its
covariance matrix is easily found to be
\begin{equation}
  \label{covar1}
  \la M_{ij} M_{kl} \ra_{N}  = 
\begin{cases} 1\,, & (k,l)=(j,i) \,, \\
0 \,, & \textup{otherwise}  \,.
\end{cases}
\end{equation}

Expectations of any monomials in the $M_{ij}$'s can be computed
using Wick's rule: the expectation 
 is a sum over all ways to group the
factors in pairs of the products of the pair covariances.
For example:
\begin{multline*} 
\la M_{ab} M_{cd} M_{ef} M_{gh} \ra = 
\la M_{ab} M_{cd} \ra \la  M_{ef} M_{gh} \ra + \\
\la M_{ab} M_{ef} \ra \la  M_{cd} M_{gh} \ra +
\la M_{ab} M_{gh} \ra \la  M_{cd} M_{ef} \ra  = \\
\delta_{ad} \delta_{bc} \delta_{eh} \delta_{fg}+
\delta_{af} \delta_{be} \delta_{ch} \delta_{dg}+
\delta_{ah} \delta_{bg} \delta_{cf} \delta_{de}
\,.
\end{multline*}
The combinatorics of such expansions can be very
conveniently handled using diagrammatic techniques (a very
accessible introduction to this subject can be found in \cite{Z}). 

\begin{figure}[!hbt]
\centering
\scalebox{.7}{\includegraphics{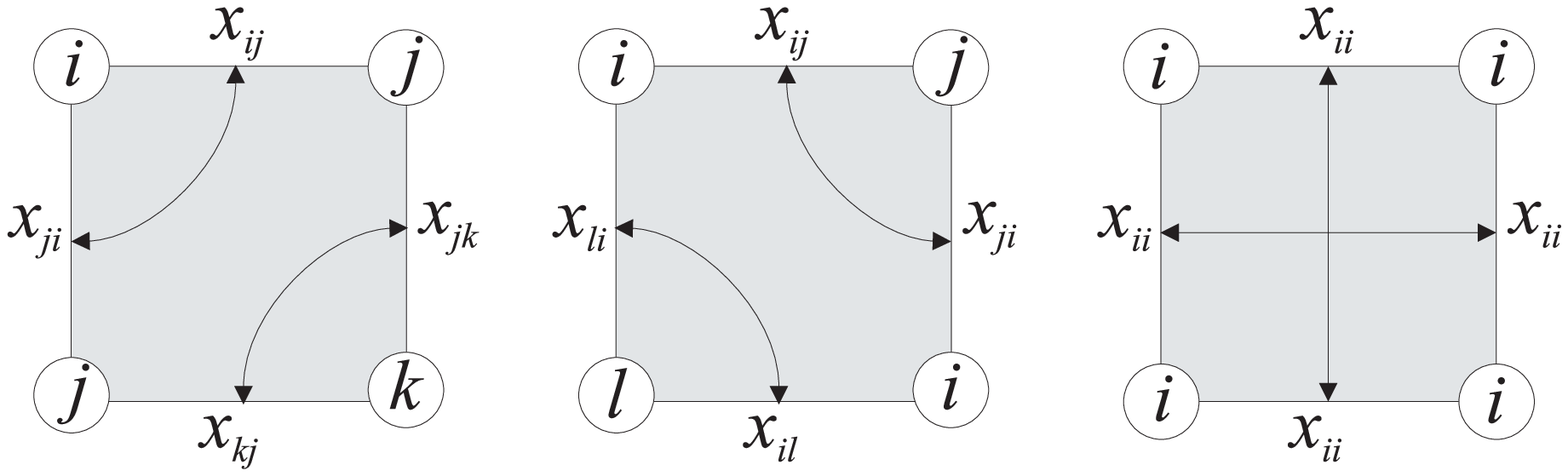}}
\caption{Diagrammatic interpretation of $\la \tr M^4 \ra_{N}$}
\label{fig4}
\end{figure}

For example, the  diagrammatic interpretation of the expectation
$$
\la \tr M^4 \ra_{N} = \sum_{i,j,k,l=1}^N \la M_{ij} M_{jk} M_{kl} 
M_{li} \ra_N  
$$ 
is the following. We place the indices 
$i$, $j$, $k$, $l$ on the vertices of a square and place
the matrix elements $M_{ij}$, $M_{jk}$, $M_{kl}$, $M_{li}$
on the corresponding edges. The pairing in Wick's 
formula can be interpreted as gluing pairs of sides of the
square together. Formula \eqref{covar1} implies  then that the side
identifications 
have to satisfy:
\begin{itemize}
\item[(i)] identified  vertices  carry equal indices,
\item[(ii)] the result is a closed and orientable surface, 
\end{itemize}
see Figure 
\ref{fig4}. Since each combinatorial scheme in
Figure \ref{fig4} contributes a power of $N$ for every
vertex on the resulting surface, we conclude that
$$
\la \tr M^4 \ra_{N} = 2 N^3 + N \,.
$$ 
Similarly, the expectation $\la \tr M^k \ra_{N}$
can be diagrammatically interpreted as counting
surfaces glued out of a $k$-gon. 

\subsubsection{Asymptotics of maps on surfaces}

More generally, an expectation of the form 
$$
\la \prod_{j=1}^l \tr M^{k_i}\ra_N
$$
counts surfaces that one can glue out of a $k_1$-gon,
$k_2$-gon, $\ldots$, and a $k_l$-gon. More specifically, each polygon here 
comes with a choice of a special vertex because the 
monomial 
\begin{equation}
  \label{Mijk}
  M_{i_1 i_2} M_{i_2 i_3} \cdots M_{i_k i_1}
\end{equation}
corresponds to a  $k$-gon  diagram 
with factors $M_{i_r i_{r+1}}$
placed on its edges together with a choice of the
vertex from which we start reading the word \eqref{Mijk}. 

As a matter of fact, we have already encountered such a combinatorial
structure under the name of a "map". Indeed, if a surface $\Su$ is
 glued out of $l$ polygons, then the
boundaries of the polygons form, according to the definition
in Section \ref{kmod}, a map on the surface $\Su$
with $l$ cells. It follows that: 
\begin{equation}
  \label{ESMM}
  \frac1
{N^{|k|/2}}\la \prod_{j=1}^l \tr M^{k_i}\ra_N=  \sum_{\Su} 
N^{\chi(\Su)-l} \, \Map_\Su(k_1,\dots,k_l) \,. 
\end{equation}
Here, $|k|=\sum k_i$. The summation is over all orientable,
but not necessarily connected, homeomorphism classes of
surfaces $\Su$.  $\Map_\Su(k_1,\dots,k_l)$ is the number of maps $G$ on 
$\Su$ satisfying:
\begin{enumerate}
\item[(i)] $G$ is a map on $\Su$ with $l$ cells marked by $1,\dots,l$,
\item[(ii)] the perimeters of cells (in the usual graph metric) are
$k_1,k_2,\dots,k_l$,
\item[(iii)] on the boundary of each cell, one vertex is specified
as the first vertex.
\end{enumerate}
The isomorphisms of such objects are isomorphisms
of the underlying maps which preserve the additional structure. 
The choice of a vertex at the boundary of each cell eliminates
all nontrivial automorphisms. 

As the function  $\Map_\Su(k_1,\dots,k_l)$ vanishes unless $|k|$
is even,  we will assume $|k|$ to be
even. Also, as the enumeration of maps on disconnected
surfaces is easily deduced from the connected case, we
will study the function $\Map_{g}$ enumerating maps on
the genus $g$ connected surface $\Su_g$.

Consider now the limit as  
the $k_i$'s increase to infinity at fixed relative
rates. In other words, introduce an extra parameter $\kappa$
and assume that
$$
\frac{k_i}{\kappa} \to x_i \,, \quad \kappa \to\infty \,.
$$ 
The following limit
\begin{equation}
  \label{map}
  \map_g(x_1,\dots,x_l)= \lim_{\kappa\to\infty}
 \frac{\Map_{g}(k_1,\dots,k_l)}
{2^{|k|} \, \kappa^{3g-3+3l/2}}
\end{equation}
was computed in \cite{O2} and, by comparison
with Kontsevich's combinatorial model, it was observed
that
\begin{equation}
  \label{Pmap}
P_g(x_1,\dots,x_l) = \frac{\pi^{l/2}}{2^g} \, 
\frac{\map_g(2x_1,\dots,2x_l)}{\sqrt{x_1\cdots x_l}} \,,  
\end{equation}
where $P_g$ denotes the $l$-point function defined in Section \ref{hurasym}. 

\begin{figure}[!hbt]
\centering
\scalebox{.4}{\includegraphics{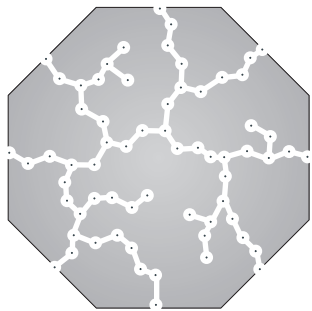}}
\caption{The map on $\Sigma_2$ corresponding to the
graph from Figure \ref{fig2}}
\label{fig6}
\end{figure}

Comparing 
Proposition \ref{lhur} to (\ref{Pmap}), we find the asymptotics of
the enumeration of branching graphs $\bH_{g,\mu}$ and
the asymptotics of map enumeration by $\Map_g(k_1,\dots,k_l)$ 
are closely related.
Each branching graph determines an underlying map
by forgetting edge labels (see Section \ref{grasym}).
The branching graph of Figure \ref{fig2} determines
the map shown in Figure
\ref{fig6}.  
The function from branching graphs to underlying maps
destroys the perimeter data of the branching graph.
However, the asymptotic
behavior of perimeters is governed by 
basic principles which apply for both 
the branching graphs {\em and} the
underlying maps. Borrowing terminology
from statistical physics, the
enumeration of branchings graphs by their perimeters and the 
enumeration of maps
by their perimeters belong to the same {\em universality class}.  
This universality class is quite large (see, for example, \ \cite{Sosh}).
Another classical combinatorial problem in the same 
universality class is the problem of increasing subsequences
in a random permutation, see \cite{O2}. 
The methods that  we use in Sections \ref{st2} and \ref{fin2} to
analyze the asymptotics of the Hurwitz numbers are parallel
to the methods used in \cite{O2} for the asymptotic enumeration of
maps. 

In the case of branching graphs, the asymptotics is related to
the intersection theory of $\mgn$ by
Proposition \ref{lhur}.  
Therefore, a conceptual
explanation of relation \eqref{Pmap} is obtained (complementing the
derivation of \cite{O2}).

\subsubsection{Edge of the spectrum}

The asymptotic function $\map_g$ has a natural extension $\map_\Su$
to disconnected surfaces $\Su$ which satisfies the obvious 
multiplicativity in connected components. Formulas
\eqref{ESMM} and \eqref{map} together imply, provided each
$k_i$ is even, the limit:
\begin{equation}
  \label{asmap}
\la \prod_{j=1}^l \tr 
\left(\frac{M}{2\sqrt N}\right)^{k_i}\ra_N \to \sum_{\Su} 
 \map_\Su(x_1,\dots,x_l) \,,
\end{equation}
as $N\to\infty$ and $k_i\to\infty$ in such a way that 
$$
\frac{k_i}{N^{2/3}} \to x_i \,.
$$
In case some of the $k_i$ are odd, certain distributions
of the $k_i$ between the connected pieces of $\Sigma$ become
prohibited by parity and, consequently, the corresponding terms in 
\eqref{asmap} should be omitted. 

It is well known (see \ \cite{Me}) that, as $N\to\infty$, 
the eigenvalue distribution of the scaled matrix $\dfrac{M}{2\sqrt N}$
converges to the (non-random) semicircle law with density
$$
\frac{2}{\pi} \, \sqrt{1-x^2} \, dx \,, \quad  x\in [-1,1] \,.
$$
It is clear that the eigenvalues near the edges $\pm 1$ of the
spectrum 
make the maximal contribution to the traces of large powers
of $M\big/2\sqrt N$ in \eqref{asmap}. This is why we call the matrix
model \eqref{asmap} the edge-of-the-spectrum matrix model. 

The behavior of eigenvalues near the edges $\pm1$ in the $N\to\infty$
limit is very well studied, see for example \cite{TW}. Let 
$\rho(x_1,\dots,x_l; N)$ denote the $l$-point correlation
function for the eigenvalues of $M\big/2\sqrt N$. 
By definition, $\rho(x_1,\dots,x_l;N) \, \prod dx_i$ is the probability 
of finding an eigenvalue in each of the  infinitesimal intervals
$[x_i,x_i+dx_i]$. These
correlation functions have  the following  $N\to\infty$
asymptotics
\begin{equation}
  \label{KAi}
 N^{-2l/3} \rho\left(1+\frac{x_1}{N^{2/3}}, \dots , 1+\frac{x_l}{N^{2/3}}
\right) \to \det\big[\KAi(x_i,x_j)
\big]_{1\le i,j \le l} \,, 
\end{equation}
where $\KAi$ is the following kernel involving the classical
Airy function 
$$
\KAi(x,y)=\frac{\Ai(2x)\, \Ai'(2y) -\Ai'(2x)\, \Ai(2y)}{x-y}  \,.
$$
The formula \eqref{KAi} together with \eqref{Pmap} results in a 
closed Gaussian integral formula for the $l$-point function
$P_g$, see \cite{O1}. It also shows that the appearance of Airy 
functions in both \eqref{KAi} and \cite{K1} is not a coincidence.

Another application of the edge-of-the-spectrum matrix model
is the following. After Kontsevich's combinatorial formula
\eqref{dfgg} is established, the derivation of
Witten's KdV equations requires an additional analysis. 
Kontsevich's original approach was to
study an associated matrix integral (Kontsevich's matrix
model) which will be discussed in Section \ref{kmod2}.
Alternatively, one can deduce, as was done in \cite{O1}, the KdV equations using
the edge-of-the-spectrum model and the  the work of
Adler, Shiota, and van Moerbeke \cite{ASV}. 

\subsection{Kontsevich's matrix model}
\label{kmod2}

 Let
$\Lambda$ be a diagonal $N\times N$ matrix with
positive real eigenvalues $s_1, \ldots, s_N$.
Instead of the Gaussian measure \eqref{gaussme} one can
consider a more general Gaussian measure on the space
of Hermitian $N\times N$ matrices $M$ with
density $\displaystyle e^{-\frac12\tr \Lambda M^2}$. We 
denote expectations of a function $f(M)$ respect to this measure
by
$$
\la f \ra_{N, \Lambda}  =  \frac{\displaystyle \int 
f(M) \, \exp\left(-\tr \Lambda M^2/2\right) \, dM}
{\displaystyle \int 
\exp\left(-\tr \Lambda M^2/2\right) \, dM}
$$
The covariance matrix of this Gaussian measure is 
easily found to be: 
$$
\la M_{ij} M_{kl} \ra_{N, \Lambda}  = 
\begin{cases} \dfrac{2}{s_i+s_j}\,, & (k,l)=(j,i) \,, \\
0 \,, & \textup{otherwise}  \,.
\end{cases}
$$ 
Expectations of any monomials in the $M_{ij}$ can be again computed
using  Wick's rule.

Kontsevich's matrix integral 
$\Theta_N$ is defined by:
\begin{equation}
\label{mintt}
 \Theta_N(s_1,\ldots, s_N) = \la \exp\left( \frac{i}6 \, \tr M^3 \right)
\ra_{N,\Lambda} \,. 
\end{equation}
Expanding the exponential by Taylor series and 
applying Wick's formula leads to the expansion: 
\begin{equation}
\label{zxc}
\text{log} \
\Theta_N(s_1,\ldots,s_N) =
\sum_{g\geq 0} \sum_{n\geq 1}   
(-2)^{2g-2+n}\sum_{G \in \bGt_{g,n}(N)}
\frac{1}{|\text{ Aut}(G)|} \prod_{e\in E} \frac{1}{\widetilde{s}(e)},
\end{equation}
where
$\bGt_{g,n}(N)$ 
denotes the set of trivalent maps with $n$ marked cells labeled
by a subset of the numbers $\{1,2,\dots,N\}$. The logarithm
function in \eqref{zxc} has the effect of selecting only connected
diagrams.

A change of variables is required  to relate $\Theta_N$
to the free energy $F$ arising in Witten's conjectures. 
Let $t^N$ denote the variable set $\{t^N_i\}_{i=0}^\infty$.
For 
$i \geq 0$, let
\begin{equation}
\label{vvb}
t^N_i = - \sum_{k=1}^N \frac{(2i-1)!!}{s_k^{2i+1}}.
\end{equation}
Substitution into $F$ yields:
\begin{eqnarray*}
F(t^N) & = & \sum_{n\geq 1,\ k_1,\ldots,k_n}
\frac{1}{n!} \langle \tau_{k_1}\cdots \tau_{k_n}\rangle \ t^N_{k_1} \cdots t^N_{k_n} \\
& = &
\sum_{n\geq 1,\ k_1,\ldots,k_n}
\frac{(-1)^n}{n!} \langle \tau_{k_1}\cdots \tau_{k_n}\rangle 
\sum_{1\leq l_1,\ldots,l_n \leq N} \prod_{i=1}^n \frac{(2k_i-1)!! 
}{s_{l_i}^{2k_i-1}} \\
& = & 
\sum_{g\geq 0} \sum_{n\geq 1}   
(-2)^{2g-2+n}\sum_{G \in \bGt_{g,n}(N)}
\frac{1}{|\text{ Aut}(G)|} \prod_{e\in E} \frac{1}{\widetilde{s}(e)},
\end{eqnarray*}
The last equality is a consequence of Theorem \ref{kon}.
Therefore, $$F(t^N)= \text{log}\ \Theta_N(s_1,\ldots,s_N).$$
As $N\rarr \infty$, the change of variables (\ref{vvb}) is faithful to
higher and higher orders. The entire function $F$ may be recovered
in the large $N$ limit.

\begin{tm} $F$ is the large $N$ limit of Kontsevich's matrix
model:
$$F(t)= \lim _{N\rarr \infty} \ \Theta_N(t^N).$$
\end{tm}

Witten's KdV equations for $F$ are proven in \cite{K1} 
from the analysis of Kontsevich's matrix integral. An exposition 
of this analysis can be found in \cite{dF,DFGZ,DIZ}.

\subsection{Matrix models of $2$-dimensional quantum gravity}

In  quantum gravity, one wishes to compute a
Feynman integral of matter fields over all possible topologies and metrics
on a 2-dimensional worldsheet. One
way to make mathematical sense out of such integration is 
to interpret the result as a suitable integral over the moduli
spaces of curves, see \cite{W}. Another approach is to 
discretize the problem: instead of all possible metrics
one can consider, for example, only surfaces
tessellated into standard squares, or into more general 
polygons. In a suitable limit, in which the number of
tiles goes to infinity, one expects to be able to compute physically
significant quantities from this approximations.

Diagrammatic techniques for matrix integrals provide
a very powerful tool for enumerating tessellations and
investigating their asymptotic behavior (see, for example, 
the surveys \cite{dF,DFGZ} as well the original papers  
\cite{BK,Dou,DS,GM1,GM2}. More concretely, consider an integral
over the space of $N\times N$ Hermitian matrices 
of the following  form 
$$
Z(V,N)=\int e^{-N \tr V(M)} \, dM
$$
where 
$$
V(x) = \tfrac12 x^2 +\gamma(x)  \in \R[x] 
$$
is a polynomial (usually assumed to
be even). After an expansion by Wick's formula,  $Z(V,N)$ yields a 
weighted enumeration of 
surfaces tessellated into polygons. The weight involves the
genus of the surface, the automorphisms group of the tessellation, and 
the coefficients of the polynomial $V$ corresponding to the
tiles of the tessellation.

The physically interesting limit (the {\em double
scaling} limit) is obtained when the coefficients of the polynomial 
$V$ approach certain critical values simultaneously as
$N\rarr \infty$.
Formal manipulation with asymptotics of orthogonal polynomials
shows that this limit is governed by the KdV hierarchy, 
see for example \cite{dF,DFGZ} for a survey.  This is precisely what 
led Witten to conjecture that the same hierarchy describes
intersections on the moduli spaces of curves. 

However, rigorous mathematical investigation of the corresponding
double-scaling asymptotics of orthogonal polynomials is a very difficult 
problem.  At present, only the case of even quartic potential $V$ has 
been analyzed completely \cite{BI}. In this respect, the matrix integral
$Z(V,N)$, is a much more problematic object than Kontsevich's
matrix model or the edge of the spectrum  matrix model. 

\subsection{The Toda equation for $\proj^1$}
\label{tttod}

The moduli space $\overline{M}_{g,n}$ may be viewed as
the moduli space of maps to a point.
The Hurwitz path to matrix models is
found in the geometry of maps to $\proj^1$. 
It is perhaps natural then to seek a 
link between the Gromov-Witten theory of 
target varieties $X$ and matrix models via the
geometry of maps to $X\times \proj^1$.
While this direction has promise,
no constructions have yet
been found even for $X=\proj^1$.

Instead,
the study of the Gromov-Witten theory of the
target variety $X=\proj^1$ is {\em again}
linked to the Hurwitz numbers. 
The Toda equation (conjecturally) constrains
the free energy $F$ of $\proj^1$. 
The generating series
$H$ of the Hurwitz numbers has been proven
to satisfy an analogous Toda equation via a representation theoretic analysis
of $H_{g,\mu}$ \cite{O3}.
The functions $F$ and $H$ may be partially
identified through the basic Hurwitz numbers $H_{g,1^d}$ \cite{P1}.
The two Toda equations agree in this region of overlap.


We explain here the basic relationship between 
the Gromov-Witten theory
of $\proj^1$, the Hurwitz numbers, and the Toda equation.
The tautological classes in $H^*(\overline{M}_{g,n}(\proj^1,d), {\mathbb{Q}})$
which we will consider are of two types.
First,
the classes $\psi_i$ are defined on
the moduli space $\overline{M}_{g,n}(\proj^1,d)$  
by the same construction used for $\overline{M}_{g,n}$:  $\psi_i$ is
the Chern class of the $i^{th}$ cotangent line bundle.
The tautological evaluation maps,
$$\text{ev}_j: \overline{M}_{g,n}(\proj^1,d)\rarr \proj^1,$$
defined for each marking $j$ provide a structure not present in the study of 
$\overline{M}_{g,n}$. The second type of tautological class is:
$$\text{ev}_j^*(\omega) \in
H^2(\overline{M}_{g,n}(\proj^1,d), {\mathbb{Q}}),$$
where $\omega\in H^2(\proj^1, {\mathbb{Q}})$ is the point class.
The intersections of products of $\psi_i$ and $\text{ev}_j^*(\omega)$
in $\overline{M}_{g,n}(\proj^1,d)$
are the {\em gravitational descendents} of $\proj^1$.   
The bracket notation for the descendent integrals is:
\begin{equation}
\label{nedd}
\langle  \prod_{i=1}^r \tau_{a_i} \cdot \prod_{j=r+1}^{r+s}
\tau_{b_j}(\omega) \rangle_{g,d}^{\proj^1} =
\int_{[\overline{M}_{g,n}(\proj^1,d)]^{vir}} \prod_{i=1}^r \psi_{i}^{a_i} \cdot
\prod_{j={r+1}}^{r+s} \psi_j^{b_j} \text{ev}_j^*(\omega).
\end{equation}
All integrals in Gromov-Witten theory are evaluated against the
virtual fundamental class.

The free energy $F$ of $\proj^1$ is a complete generating function
of the integrals (\ref{nedd}).
Let the variables $x_i$ and $y_j$ correspond to the classes
$\psi_i$ and $\text{ev}_j^*(\omega)$. Let $x$ and
$y$ denote the sets of variables $\{x_i\}_{i=0}^\infty$ and 
$\{y_i\}_{i=0}^\infty$
respectively.
$F$ is defined by the formula:
\begin{equation}
\label{ggg}
F(\lambda,x,y)= \sum_{g\geq 0} 
\sum_{d\geq 0} \sum_{n\geq 0} \lambda^{2g-2} \frac{
\langle \gamma^n \rangle_{g,d}^{\proj^1}}{n!},
\end{equation}
where $\gamma= \sum_{i\geq 0} x_i \tau_i+y_i\tau_i(\omega)$. 
The bracket in (\ref{ggg}) is viewed as linear in the
variables $x$ and $y$.

The (conjectural) Toda
equation for $F$ may be written in the following
form:
\begin{equation}
\label{todaa}
\exp\Big( F(x_0+\lambda) + F(x_0-\lambda)
- 2 F \Big) = \lambda^2  F_{y_0 y_0},
\end{equation}
where $F(x_0\pm \lambda)= F(\lambda, x_0\pm \lambda, x_1,x_2,\ldots,
y_0,y_1,y_2,\ldots).$
Equation (\ref{todaa}) has its origins in 
the study of matrix models believed to be related
to the Gromov-Witten theory of $\proj^1$ \cite{EgY}.
Proofs of the genus 0 and 1 implications of the
Toda equation can be found in \cite{P1}.
The  Toda equation (\ref{todaa}) determines $F$
from degree $d=0$ descendent invariants of $\proj^1$.

A very similar Toda equation holds for the generating function of
the Hurwitz numbers $H_{g,\mu}$. 
Let $p$ denotes the variable set $\{p_i\}_{i=1}^\infty$.
For each partition $\mu$ of $d$ with parts $\mu_1,\ldots,\mu_l$,
let $$p_\mu = p_{\mu_1} \cdots p_{\mu_l}.$$
Define the Hurwitz generating function $H$ by:
$$H(\lambda,y_0, p) =
\sum_{g\geq 0} \sum_{d > 0} \sum_{\mu \vdash d} \lambda^{2g-2}e^{dy_0}
\frac{H_{g,\mu}}{(2g-2+d+l)!} \ p_\mu.$$
The definition of the Hurwitz numbers 
via the symmetric group in Section \ref{hur}
may be used to prove a Toda equation for $H$. More precisely, the function 
$H$ is linked to the Toda lattice 
hierarchy of Takasaki and Ueno in representation theory \cite{O3}.
One specialization of this hierarchy is the following: 
\begin{pr}
\label{ror}
$H$ satisfies the Toda equation:
\begin{equation}
\label{todah}
\exp\Big( H(y_0+\lambda) + H(y_0-\lambda)
- 2 H \Big) = \lambda^2 e^{-y_0} H_{p_1 y_0}.
\end{equation}
\end{pr}
\noindent

The Toda equations for the free energy $F$ and the 
Hurwitz function $H$
are connected through a partial identification
of these two series. Perhaps the Toda equation
for $F$  could be proven by
a better understanding of this relationship.

Let
$H_{g,d}$ be the Hurwitz number of genus $g$,
degree $d$, simply ramified covers of
$\proj^1$. By definition, $H_{g,d}$ equals $H_{g,1^d}$.
The generating function $\tilde{H}$ of the Hurwitz numbers $H_{g,d}$
is obtained by a restriction of $H$:
\begin{eqnarray*}
\tilde{H}(\lambda, y_0) & = & 
\sum_{g\geq 0} \sum_{d>0} \lambda^{2g-2} e^{dy_0} \frac{H_{g,d}}{(2g-2+2d)!} \\
& = & H(\lambda,y_0, p_1=1, p_{i\geq 2}=0).
\end{eqnarray*}

The Hurwitz numbers $H_{g,d}$ arise in Gromov-Witten theory
as 
descendent integrals of $\proj^1$
 \cite{P1}.
\begin{pr}
\label{yay}
For all $g \geq 0$ and $d>0$,
$$H_{g,d} = \langle \tau_1(\omega)^{2g+2d-2} \rangle^{\proj^1} _{g,d}.$$
\end{pr}
\noindent
The generating function $\tilde{H}$ is therefore obtained 
by a restriction of $F$:
\begin{eqnarray*}
\tilde{H}(\lambda, y_0) & = &  
\sum_{g\geq 0} \sum_{d>0} 
\lambda^{2g-2} e^{dy_0} 
\frac{\langle \tau_1(\omega)^{2g+2d-2} \rangle^{\proj^1} _{g,d}
}{(2g-2+2d)!} \\
& = &
 F(\lambda,x_{i\geq 0}=0,y_0, y_1=1,y_{i\geq 2}=0). 
\end{eqnarray*}

There are two natural Toda equations for $\tilde{H}$ obtained from
the Toda equations for $F$ and $H$ respectively.
\begin{tm}
\label{cvb}
The two Toda equations (\ref{todaa}) and (\ref{todah}) 
specialize to a unique Toda equation for $\tilde{H}$:
\begin{equation}
\label{todahs}
\exp\Big( \tilde{H}(y_0+\lambda) + \tilde{H}(y_0-\lambda)
- 2 \tilde{H} \Big) = \lambda^2 e^{-y_0} \tilde{H}_{y_0 y_0}.
\end{equation}
\end{tm}
\noindent Theorem \ref{cvb} provides strong evidence
for the (conjectural) Toda equation for $F$.

The Toda equation for the Hurwitz series $H$ was found in
the search for a proof of prediction (\ref{todahs}) of
the Toda equation for $\proj^1$. 
One may reasonably hope
the connection between the Toda equations for $H$ and
$F$ is stronger than Theorem \ref{cvb}. However,
a direct extension of Proposition \ref{yay} relating 
all the Hurwitz numbers $H_{g,\mu}$ to descendents 
has not been discovered. 
The natural context for the Toda equation in \cite{O3}
suggests the larger class of {\em double Hurwitz numbers} 
may be related fundamentally to the Gromov-Witten theory
of $\proj^1$.


\part{Hurwitz numbers in Gromov-Witten theory} 

\section{\bf Gromov-Witten theory of $\proj^1$}
\label{st1}
\subsection{Stable maps}
Let $X$ be a nonsingular projective variety.
A path integral over the space of differential maps
$\pi:\Sigma_g\rarr X$ naturally arises in the topological 
gravity theory with target $X$.
A stationary phase
analysis then yields the following
string theoretic result: the path integral localizes to
the space of holomorphic maps from Riemann surfaces to $X$ \cite{W}.
The path integral therefore should be equivalent to 
classical integration over a space
of holomorphic maps.

The moduli of maps may be studied in algebraic geometry
by the equivalence of the holomorphic and algebraic
categories in complex dimension 1.
However, the moduli space $M_{g,n}(X,\beta)$ of $n$-pointed algebraic 
maps $\pi:(C, p_1,\ldots,p_n)\rarr X$ satisfying
\begin{enumerate}
\item[(i)] $C$ is a nonsingular curve of genus $g$,
\item[(ii)] $p_1,\ldots,p_n\in C$ are distinct points,
\item[(ii)]$\pi_*[C] = \beta \in H_2(X, {\mathbb{Z}})$,
\end{enumerate}
is not compact.
For example, the domain may degenerate to a nodal curve,
the points may meet, 
or the map itself may acquire a singularity. The compactification
$$M_{g,n}(X,\beta)\subset \overline{M}_{g,n}(X,\beta)$$ by stable maps 
plays a central role
in Gromov-Witten theory --- it is conjectured to be the
correct compactification for calculating the path integral
of the gravity theory.

The moduli space of stable maps
$\overline{M}_{g,n}(X,\beta)$
parameterizes 
$n$-pointed algebraic maps
$$\pi: (C, p_1,\ldots,p_n)\rarr X$$ satisfying:
\begin{enumerate}
\item[(i)]
$C$ is a compact, connected, reduced, (at worst) nodal curve
of arithmetic genus $g$,
\item[(ii)] $p_1, \ldots, p_n \in C$ are distinct and lie in 
           the nonsingular locus,
\item[(iii)] $\pi_*[C]=\beta$,
\item[(iv)]  $\pi$ has no infinitesimal automorphisms.
\end{enumerate}
A {\em special} point of the domain $C$ is a marked point $p_i$ or
a nodal point. An {\em infinitesimal automorphism} of a map $\pi$ is a tangent field $v$
of the
domain $C$ which vanishes at the special points and satisfies $d\pi(v)=0$.
Stable maps were defined by Kontsevich in \cite{K1,KMa}. A
construction of the moduli space can be found in \cite{FP}.

An irreducible component $E\subset C$ is $\pi$-{\em collapsed} if the
image $\pi(E)$ is a point.
Property (iv) is equivalent to a geometric
condition on each $\pi$-collapsed component:
$\pi$ has no infinitesimal automorphisms if and only if
the normalization $$\tilde{E}\rarr E$$
of each $\pi$-collapsed component $E$
contains the preimages of at least $3-2g(\tilde{E})$ special points of $C$.
As $3-2g(\tilde{E})>0$ only if $g(\tilde{E})=0$ or $1$, this condition
only constrains rational and elliptic components.
If the entire domain $C$
is $\pi$-collapsed, property (iv) is equivalent to the Deligne-Mumford
stability
condition for pointed curves $(C,p_1,\ldots,p_n)$.
The moduli space $\overline{M}_{g,n}(X,0)$
is therefore isomorphic to $X\times \overline{M}_{g,n}$.
In particular, $\overline{M}_{g,n}$ is recovered 
as the space of stable maps to a point.

The moduli space $\overline{M}_{g,n}(X,\beta)$ is not always a
{\em nonsingular} Deligne-Mumford stack --- in fact, the space may be singular, non-reduced,
reducible, and of impure dimension. 
While ${M}_{g,n}(X,\beta)\subset \overline{M}_{g,n}(X,\beta)$ 
is an open subset, the inclusion is not necessarily dense.
The space of stable maps may
be quite complicated even when $M_{g,n}(X,\beta)$ is empty.

Most pathologies occur even in case $X=\proj^1$. 
Consider, for example, $\overline{M}_{2}(\proj^1,2)$.
The closure of the locus of hyperelliptic maps $M_{2}(\proj^1,2)$
yields an irreducible component of $\overline{M}_{2}(\proj^1,2)$
of dimension 6. However, 
the set of maps
obtained by attaching a $\pi$-collapsed genus 2 curve 
to a rational double cover of $\proj^1$ 
forms another component of dimension 7. 
In fact, $\overline{M}_{2}(\proj^1,2)$ contains 7 irreducible
components in all.
One of the few global geometric properties always satisfied 
by $\overline{M}_{g,n}(\proj^1,d)$ is connectedness \cite{KiP}.

\subsection{Branch morphisms}
\label{brrrch}
Let $g\geq 0$ and $d>0$.
The moduli space $\overline{M}_{g}(\proj^1,d)$ supports
a natural branch morphism $\br$ which will play
a basic role in the study of the Hurwitz numbers.

The branch morphism is first constructed for the open
moduli space
 $M_{g}(\proj^1,d)$. 
Let $\pi:C \rarr \proj^1$ be a degree $d$ map with a nonsingular
domain $C$. A branch divisor on $\proj^1$
is obtained from the ramifications of $\pi$.
More precisely, the branch divisor $\br(\pi)$ is the $\pi$ push-forward
of the degeneracy locus of the
differential map on $C$:
\begin{equation}
\label{bboo}
\pi^*\omega_{\proj^1} \rarr \omega_C,
\end{equation}
where $\omega_{\proj^1}$ and $\omega_C$ denote the canonical bundles
of $\proj^1$ and $C$ respectively.
By the Riemann-Hurwitz formula, $\br(\pi)$ has degree 
$$r=2g(C)-2+2d.$$
A branch morphism from  $M_{g}(\proj^1,d)$ 
to the space of divisors,
\begin{equation}
\label{oooo}
\br: M_{g}(\proj^1,d) \rarr \text{Sym}^{r}(\proj^1),
\end{equation}
is defined algebraically by the universal degeneracy locus (\ref{bboo}).

A branch divisor $\br(\pi)$ is 
constructed for stable maps $\pi:C \rarr \proj^1$
by the following definition.
Let $N\subset C$ be the
cycle of nodes of $C$.
Let $\nu: \tilde{C} \rarr C$
be the normalization of $C$.
Let $A_1, \ldots, A_a$ be the components of $\tilde{C}$
which dominate $\proj^1$, and let $\{a_i: A_i \rarr \proj^1\}$
denote the natural maps. 
As $a_i$ is a surjective map between nonsingular
curves, the branch divisor $\br(a_i)$ is defined by (\ref{bboo}). 
Let $B_1, \ldots, B_b$ be
the components of $\tilde{C}$ contracted over $\proj^1$, and
let $f(B_j)=q_j\in \proj^1$.
Define $\br(\pi)$ by:
\begin{equation}
\label{ptwise}
\br(\pi)= \sum_{i} \br(a_i) +
\sum_j (2g(B_j)-2)[q_j] + 2\pi_{*}(N).
\end{equation}
Formula (\ref{ptwise}) associates an
effective divisor of degree $r$ on $\proj^1$ to every
moduli point $[\pi]\in \overline{M}_g(\proj^1,d)$.

The branch divisor $\br(\pi)$ for stable maps may be constructed
canonically from the complex:
\begin{equation}
\label{cmplx}
R\pi_*[\pi*\omega_{\proj^1} \rarr \omega_{C}],
\end{equation}
well-defined in the derived category. 
An effective divisor on $\proj^1$ 
is extracted from (\ref{cmplx}) via a determinant construction.  
An algebraic branch morphism
\begin{equation}
\label{ooooo}
\br: \overline{M}_{g}(\proj^1,d) \rarr \text{Sym}^{r}(\proj^1)
\end{equation}
is then obtained from the universal complex (\ref{cmplx}).
The required derived category arguments can be found in \cite{FanP}.

\subsection{Virtual classes}
\subsubsection{Perfect obstruction theories}
\label{bab1}
Let $X$ be a nonsingular projective variety.
The {\em expected} or {\em virtual}  dimension  of the moduli space
$\overline{M}_{g,n}(X, \beta)$ is:
$$\int_\beta c_1(X) + \text{dim}(X)(1-g) + 3g-3+n.$$
$\overline{M}_{g,n}(X,\beta)$ carries a canonical obstruction
theory which yields a {\em virtual class}
$$[\overline{M}_{g,n}(X,\beta)]^{vir} \in A_{\text{exp}}
(\overline{M}_{g,n}(X,\beta), {\mathbb{Q}})$$
in the expected rational Chow group. The virtual class of $\overline{M}_{g,n}(X,\beta)$
was
first constructed in \cite{LiT,B,BFan}.
The virtual class plays a fundamental role in Gromov-Witten
theory --- all cohomology evaluations in the theory are
taken against the virtual class.

The virtual class of $\overline{M}_{g,n}(X,\beta)$ is 
constructed via a canonical {perfect obstruction theory}
carried by the moduli of maps.
A perfect obstruction theory on scheme (or Deligne-Mumford
stack) $V$ consists of the following data:
\begin{enumerate}
\item[(i)] A two term complex of vector bundles
$E^\bullet=[E^{-1} \rarr E^0]$ on $V$.
\item[(ii)] A morphism $\phi: E^\bullet \rarr L^\bullet_{V}$
in the derived category $D_{qcoh}^-(V)$ 
to the cotangent complex $L^\bullet_{V}$
satisfying two properties:
\begin{enumerate}
\item[(a)] $\phi$ induces an isomorphism in cohomology in degree 0.
\item[(b)] $\phi$ induces a surjection in cohomology in degree -1.
\end{enumerate}
\end{enumerate}
A virtual fundamental class of dimension
$dim(E^0) -dim(E^{-1})$ is canonically associated to the data (i) and (ii).

\subsubsection{Categories of complexes}
Let 
$C_{qcoh}^-(V)$ be the category of
complexes of quasi-coherent sheaves bounded from above on $V$. The
objects of $C_{qcoh}^-(V)$ are complexes,
$$F^\bullet=[\ldots \rarr F^{-1} \rarr F^0 \rarr F^1 \rarr \ldots],$$
satisfying $F_i=0$ for $i$ sufficiently large.
The morphisms of $C_{qcoh}^-(V)$ are chain maps of complexes.

A chain map $\sigma: F^\bullet \rarr \tilde{F}^\bullet$ is
a {\em quasi-isomorphism} if $\sigma$ induces an
isomorphism on cohomology: $H^*(\sigma):H^*(F^\bullet) \rarr H^*(\tilde{F}^\bullet)$.

The objects of derived category $D_{qcoh}^-(V)$ are also
complexes of quasi-coherent sheaves bounded from above on $V$. However,
the morphisms of $D_{qcoh}^-(V)$ are obtained by {\em inverting}
all quasi-isomorphisms in $C_{qcoh}^-(V)$.  
A basic result is  
a {\em morphism} 
$F^\bullet\rarr G^\bullet$ in $D_{qcoh}^-(V)$ may be represented by a
diagram:
\begin{equation*}
\begin{CD}
\tilde{F}^\bullet
  @>{\tau} >> G^\bullet \\
@V{\sigma}VV     \\
 F^\bullet,
\end{CD}
\end{equation*}
where $\sigma$ is a quasi-isomorphism and $\tau$ is
map of complexes.

An excellent reference for the derived category 
is \cite{GMa}. A more informal introduction may be found in
\cite{Th}.

\subsubsection{Cotangent complexes}
The cotangent complex $L^\bullet_{V}$ is a canonical object
(up to equivalence) of
$D_{qcoh}^-(V)$.
While the full complex $L^\bullet_{V}$ is constructed abstractly,
we will see the essential
properties which
are required here can be described concretely. 

If $V$ is nonsingular, $L^\bullet_{V}$ is defined by 
the 1 term complex $[\Omega_{V}]$ in degree 0 determined by
the cotangent bundle. 
A nonsingular space $V$ carries a canonical {\em trivial}
perfect obstruction theory:
$$\phi: [0\rarr \Omega_{V}] \stackrel{\sim}{\rarr} L^\bullet_{V}.$$
We will see the virtual fundamental class of this trivial
theory is the ordinary fundamental class of $V$.
For arbitrary $V$, the
cotangent complex may be viewed as a generalized cotangent
bundle.

We first note the $k$ cut-off functor is well-defined in
$D_{qcoh}^-(V)$:
$$F^{\geq k} = 
[ \frac{F^k}{Im(F^{k-1})} \rarr F^{k+1} \rarr F^{k+2} \rarr \ldots],$$
for any complex $F^\bullet$.

The cut-off $L^{\geq -1}_{V}$ of the
cotangent complex for singular $V$ may be geometrically identified by the
following construction. Let
\begin{equation}
\label{ddg}
\overline{M} \subset Y
\end{equation} be
an embedding in a nonsingular scheme (or Deligne-Mumford stack) $Y$.
The cut-off of $L^\bullet_{V}$ is represented by:
\begin{equation}
\label{nnh}
L^{\geq -1}_{\overline{M}}=[I/I^2 \rarr \Omega_Y\otimes \oh_{\overline{M}}].
\end{equation}
Here, $I$ is the ideal sheaf of $V\subset Y$.
The complex (\ref{nnh}) is independent 
(up to equivalence in the derived
category) of the embedding (\ref{ddg}). 

The representation (\ref{ddg}) easily implies the
cohomology of $L^\bullet_{V}$
in degree 0 is the sheaf of differentials $\Omega_{\overline{M}}$. 
The cohomology of $L^\bullet_{V}$ is degree -1 is also
determined by (\ref{ddg}): $H^{-1}(L^\bullet_V)$ 
encodes singularity data of $\overline{M}$.

For the study of perfect obstruction theories and virtual classes, it
will suffice to restrict the cotangent complex to the 
cut-off $L^{\geq -1}_{\overline{M}}$. 

Stack quotient constructions
of $\overline{M}_{g,n}(X,\beta)$ prove the existence of 
nonsingular embeddings (\ref{ddg}) for the moduli space of maps \cite{GrP}.
The quotient constructions also show
the abundance of locally free sheaves on $\overline{M}_{g,n}(X,\beta)$ --- a valuable
property for the derived category.

\subsubsection{Distinguished triangles}
Before proceeding, we include here 
a short review of mapping cones and distinguished
triangles in the derived category.

Let $A$ be a complex in $C_{qcoh}^-(V)$.
Let $A[1]$ denote the shifted complex with negative differential:
$$A[1]^i = A^{i+1}, \ \ d_A[1]=-d_A.$$

Let $\gamma: A^\bullet \rarr B^\bullet$
be a morphism of complexes. The {\em mapping cone} $M[\gamma]$
is the complex with terms and differentials:
$$M[\gamma]^{i}= A[1]^i \oplus B^i, \ \ (d_A[1], \gamma+d_B).$$
The mapping cone may be canonically placed in a triangle of morphisms:
\begin{equation}
\label{wertt}
A^\bullet \stackrel{\gamma}{\rarr} B^\bullet \rarr M[\gamma]^\bullet
\rarr A[1]^\bullet.
\end{equation}

A triangle of morphisms in the derived category,
\begin{equation}
\label{wwertt}
X^\bullet {\rarr} Y^\bullet \rarr Z^\bullet
\rarr X[1]^\bullet,
\end{equation}
is a {\em distinguished triangle} if there exist:
\begin{enumerate}
\item[(i)] a map of complexes $\gamma: A^\bullet \rarr B^\bullet$, 
\item[(ii)] a
triple of isomorphisms in $D_{qcoh}^-(V)$, 
$$A^\bullet \stackrel{\sim}{\rarr} X^\bullet, \
B^\bullet \stackrel{\sim}{\rarr} Y^\bullet, \
M[\gamma]^\bullet \stackrel{\sim}{\rarr} Z^\bullet,$$
\end{enumerate}
which together yield an isomorphism of the triangles (\ref{wertt}) and
(\ref{wwertt}) in the derived category. 

If the triangle (\ref{wwertt}) is distinguished,
it is easily proven that 
$$Y^\bullet {\rarr} Z^\bullet \rarr X[1]^\bullet
\rarr Y[1]^\bullet,$$
$$Z^\bullet {\rarr} X[1]^\bullet \rarr Y[1]^\bullet
\rarr Z[1]^\bullet,$$
are distinguished triangles as well. 
In this sense, the notion of a distinguished triangle has a cyclic
triangular symmetry.

Finally, we note that
a distinguished triangle yields a long exact sequence in
cohomology by a standard result in homological algebra.

\subsubsection{The perfect obstruction theory of the moduli of maps}
Let $$\overline{M}=\overline{M}_{g,n}(X,\beta)$$
The perfect obstruction theory of $\overline{M}$ is obtained
from the deformation theory of maps. The main step is
a construction of a perfect obstruction theory $\tilde{E}^\bullet$ {\em relative}
to the morphism
$$\tau: \overline{M} \rarr {\mathfrak M}$$
where ${\mathfrak M}$ is the Artin stack of quasi-stable curves.
The deformation theory of maps $$\pi:C \rarr X$$ from a fixed domain
curve $C$ is well-known: 
the tangent and obstruction spaces 
are 
$H^0(C,\pi^*T_X)$ and $H^1(C,\pi^*T_X)$ respectively.
A canonical relative perfect obstruction theory is then 
defined by:
\begin{equation}
\label{hjjh}
\tilde{E}^\bullet = [R^\bullet \rho_*(\pi^* T_X)] ^\vee \rarr L^\bullet_\tau,
\end{equation}
where $\rho: U \rarr \overline{M}$ is the universal curve
and $L^\bullet_\tau$ is the relative cotangent complex of the morphism
$\tau$ (see \cite{B}). The relative theory satisfies
conditions (a) and (b) for the morphism (\ref{hjjh}).

The relative cotangent complex $L^\bullet_\tau$ is determined
by a distinguished triangle:
\begin{equation}
\label{rrreee}
\tau^*L^\bullet_{\mathfrak{M}} \rarr L^\bullet_{\overline{M}}
\rarr L^\bullet_\tau \rarr \tau^* L^\bullet_{\mathfrak{M}}[1],
\end{equation}
which generalizes the sequence of relative differentials of a morphism.
The pull-back $\tau^*L^\bullet_{\mathfrak{M}}$ is canonically
identified on $\overline{M}$:
$$\tau^*L^\bullet_{\mathfrak{M}} \stackrel{\sim}{=} [R^\bullet
\underline {Hom}_{\oh_
{\overline{M}}}(-, \oh_U)(\Omega_\rho(P))]^\vee[-1].$$
Here, $\Omega_\rho$ is the sheaf of relative differentials on 
$U$, and $P$ is the divisor of marked points.

The absolute theory  $E^\bullet$ for $\overline{M}$
is then constructed by including the deformations of the domain curve
via a canonical distinguished triangle. 
\begin{equation}
\label{hdm}
\begin{CD}
\tau^*L^\bullet_{\mathfrak{M}}
  @>>> E^\bullet @>>> [R^\bullet\rho_*(\pi^*T_X)]^\vee  @>>> \tau^*L^\bullet
_{\mathfrak{M}} [1]\\
@VVV   @VV{\phi}V  @VVV @VVV\\
 \tau^*L^\bullet_{\mathfrak{M}}
    @>>>L_{\overline{M}}
@>>> L_{\tau}^\bullet @>>> \tau^*L^\bullet_{\mathfrak{M}}
 [1].
\end{CD}
\end{equation}
The right arrow on the top line of (\ref{hdm}) is obtained
from the canonical morphism,
$$\pi^*\Omega_X \stackrel{d\pi}{\rarr} \Omega_\rho \rarr \Omega_\rho(P),$$
together with the identification
$$ R^*\rho_* (\pi^*T_X) \stackrel{\sim}{=} R^*\underline{Hom}_{\oh_{\overline{M}}}
(-,\oh_U) (\pi^*\Omega_X).$$
The top line is then defined to be the distinguished triangle obtained
from the right arrow.
The bottom line of (\ref{hdm}) is the canonical distinguished triangle of
cotangent complexes obtained from the bottom right arrow.
The construction of the diagram is then formal once the
canonical morphisms in the right square are shown to commute.

The projectivity of $X$ may be used to find a two term
sequences of vector bundles representing both the terms and the
morphism,
$$[R^\bullet\rho_*(\pi^*T_X)]^\vee  \rarr \tau^*L^\bullet
_{\mathfrak{M}} [1],$$
in the derived category (see \cite{B, BFan}).
By the mapping cone construction, 
$E^\bullet$ then admits a three term representation:
\begin{equation}
\label{sdqww}
[E^{-1} \rarr E^0 \rarr E^1].
\end{equation}
The stability condition on the moduli space of maps
implies the cohomology of $E^\bullet$ 
vanishes in degree 1.  Hence, the sequence
(\ref{sdqww}) can be reduced to a two term complex.

The defining conditions (i) and (ii) of a perfect
obstruction theory are easily verified for:
$$\phi:E^\bullet \rarr L_{\overline{M}}^\bullet,$$
by the long exact sequence obtained from diagram (\ref{hdm}).

The diagram (\ref{hdm})  is  the primary method of studying
the obstruction theory $E^\bullet$.
Treatments can be found in \cite{B,GrP,LiT} (the latter
pursues a different perspective). A foundational exposition of
these obstruction theories will
be developed in \cite{GrKP}.

Let $[\pi:(C,p_1,\ldots,p_n) \rarr X]$ be a moduli point of $\overline{M}$.
The cohomologies of the dual complex $[E^\bullet_{[\pi]}]^\vee$ are 
the tangent and obstruction spaces of $\overline{M}$ at
$[\pi]$.
The long exact sequence in
cohomology of (the dual of) the top line of (\ref{hdm}) yields the 
the familiar tangent-obstruction sequence:
\begin{equation}
\label{toto}
0 \rarr \text{Ext}^0(\Omega_C(P), \oh_C) 
\rarr H^0(C,\pi^*T_X) \rarr \text{Tan}(\pi) 
\end{equation}
$$ \ \ \ \ \ \ \rarr \text{Ext}^1(\Omega_C(P), \oh_C) \rarr H^1(C,\pi^*T_X) 
\rarr \text{Obs}(\pi) \rarr 0.$$

The following Lemma provides a basic example of the use of the
perfect obstruction theory.
\begin{lm}
\label{qqvv}
If $H^1(C, \pi^*T_X)=0$, then $[\pi]$ is a nonsingular point
of the Deligne-Mumford stack $\overline{M}_{g,n}(X,\beta)$.
\end{lm}
\bpf
If $H^1(C,\pi^*T_X)=0$, then $\text{Obs}(\pi)=0$. By semicontinuity,
the obstruction space vanishes
for {\em every} moduli point in an open set $M$ containing $[\pi]$.
Therefore, the complex $E^\bullet$ must have 
locally free cohomology in degree 0 and vanishing cohomology in
degree -1
on $M$.
By conditions (a) and (b) of the perfect obstruction theory, the cotangent
complex must also have locally free cohomology in
degree 0 and vanishing cohomology in degree -1 on $M$.

Consider an embedding $M \subset Y$ in a nonsingular Deligne-Mumford
stack. The cut-off of the cotangent complex is
$$[I_M/I^2_M \rarr \Omega_Y \otimes \oh_M].$$
By the cohomology conditions, we conclude $I_M/I_M^2$ is
locally free {\em and} injects into $\Omega_Y\otimes \oh_M$. 
By the local criterion for nonsingularity, $M$ is nonsingular.

We note the restriction of the perfect obstruction theory to
$M$ yields the trivial perfect obstruction theory on 
a nonsingular space --- where 
$\phi$ is an isomorphism.
\epf

\subsubsection{Construction of virtual classes}
\label{bab2}
The perfect obstruction theory yields a map in the {\em derived category}
$$\phi: E^\bullet \rarr L^{\bullet}_{\overline{M}}.$$
After an exchange of representatives and cutting-off, we may assume
\begin{equation}
\label{ffa}
\phi: E^\bullet \rarr 
[I/I^2 \rarr \Omega_Y\otimes \oh_{\overline{M}}]
\end{equation}
is a map of {\em complexes}.
The virtual class is obtained from the geometry
of (\ref{ffa}).

The mapping cone associated to (\ref{ffa}) is the following
complex of sheaves:
\begin{equation}
\label{vbng}
E^{-1} \rarr E^0 \oplus I/I^2 \stackrel{\gamma}{\rarr} 
\Omega_Y \rarr 0 .
\end{equation}
The above complex (\ref{vbng}) is right exact
by conditions (a) and (b) satisfied by $\phi$. 
Let $Q$ denote the kernel of $\gamma$. $Q$ is
naturally a quotient of $E^{-1}$ by right exactness.

Let $S$ be a coherent sheaf on $\overline{M}$. The
symmetric tensors define a sheaf of $\oh_{\overline{M}}$
algebras,
$${\mathcal{S}}= \bigoplus_{k=0}^{\infty} Sym^k(S),$$
on $\overline{M}$.
The {\em abelian cone} $C(S)$ is defined to
be $$Spec({\mathcal{S}}) \rarr \overline{M}.$$
In case $S$ is a vector bundle, $C(S)$ is the
total space of $S^*$. We let $E_0$, $E_1$ denote
$C(E^0)$, $C(E^1)$ respectively.

The sequence (\ref{vbng}) yields an {\em exact sequence of abelian
cones}:
$$0 \rarr TY \rarr E_0 \times_{\overline{M}} C(I/I^2) \rarr C(Q) \rarr 0.$$
Here, the vector bundle $TY$ acts fiberwise and freely on the abelian
cone $E_0 \times_{\overline{M}} C(I/I^2)$ with quotient $C(Q)$.

Recall the normal cone $C_{\overline{M}/Y}$ is defined
by:
$$C_{\overline{M}/Y} =   Spec ( \bigoplus_{k=0}^{\infty} I^k/I^{k+1})
\rarr \overline{M}.$$  
$C_{\overline{M}/Y}$ has pure dimension equal to $dim(Y)$ (see \cite{Fu}).
There is closed embedding of $C_{\overline{M}/Y} \subset
C(I/I^2)$ given by a natural surjection of algebras:
$$\bigoplus_{k=0}^{\infty} 
Sym^k(I/I^2) \rarr \bigoplus_{k=0}^{\infty} I^k/I^{k+1}.$$

The fundamental geometric fact is that the
subcone $$E_0 \times_{\overline{M}} C_{\overline{M}/Y}\subset
E_0 \times_{\overline{M}} C(I/I^2)
$$ is invariant under the $TY$ action \cite{BFan}.
The quotient cone
$$D = \frac{E_0 \times_{\overline{M}} C_{\overline{M}/Y}}{TY}$$
is of pure dimension equal to $dim(E_0)$ and lies in $C(Q)$.
There is an embedding of abelian cones
$$C(Q) \subset E_1$$
obtained from the surjection $E^{-1} \rarr C(Q)$.
Hence $D \subset E_1$.

Let $z:\overline{M} \hookrightarrow E_1$ be the
inclusion of the zero section of the vector bundle $E_1$.
Certainly $z^{-1}(D) = \overline{M}$ as $D$ is a cone.
The {\em refined} intersection product
therefore yields a cycle class,
$$z^{!}[D] \in A_{dim(E_0)-dim (E_1)}(\overline{M},{\mathbb{Q}}).$$
The
virtual fundamental class of the perfect obstruction theory
is defined to equal $z^{!}[D]$.

The trivial perfect obstruction theory on a
nonsingular space is easily seen to yield the ordinary fundamental class
as the virtual class.

\subsubsection{Properties}
The virtual class of $\overline{M}_{g,n}(\proj^1,d)$ satisfies 
several remarkable properties --- only two of which will be
required here. 

Since the inclusion of the moduli of maps with nonsingular domains,
$$M_{g,n}(\proj^1,d) \subset \overline{M}_{g,n}(\proj^1,d),$$
is open, there is a well-defined restriction of the
virtual class.

\begin{pr} 
\label{fppp}
Let $d\geq 1$.
$M_{g,n}(\proj^1,d)$ is a nonsingular Deligne-Mumford stack of expected
dimension $2g-2+2d+n$. 
The restriction of virtual class is the ordinary fundamental class of
$M_{g,n}(\proj^1,d)$.
\end{pr}

\bpf 
Let $[\pi: (C,p_1,\ldots,p_n) \rarr \proj^1]$ determine a moduli point of
$M_{g,n}(\proj^1,d)$.
The nonsingularity, the dimensionality, and the identification of the virtual class
follow directly from the vanishing of 
$\text{Obs}(\pi)$  --- as can be seen by Lemma \ref{qqvv} the definitions
of Sections \ref{bab2}.
The  canonical right  exact sequence:
$$ \text{Ext}^1(\Omega_C(D),\oh_C) \stackrel{i}{\rarr}
 H^1(C, \pi^*T_{\proj^1}) \rarr \text{Obs}(\pi) \rarr 0$$
is obtained from the tangent-obstruction sequence (\ref{toto}).
Since $C$ is nonsingular, $\text{Ext}^1(\Omega_C(D),\oh_C) =
H^1(C, T_C(-D))$. Moreover, the map $i$ factors by:
\begin{equation}
\label{reedd} 
H^1(C, T_C(-D)) \rarr H^1(C, T_C) \rarr H^1(C, \pi^* T_{\proj^1}).
\end{equation}
The first map in (\ref{reedd}) is certainly surjective.
Since $d>0$, the sheaf map
$T_C \rarr \pi^*T_{\proj^1}$ has a torsion quotient and the
second map in (\ref{reedd}) is also surjective.
Hence, $i$ is surjective and $\text{Obs}(\pi)=0$.
\epf

The second required property of the virtual class is the $\com^*$-localization
formula discussed in Section \ref{vll}.

\section{\bf Virtual localization}
\label{vll}
\subsection{Atiyah-Bott localization}
Let $V$ be a {\em nonsingular} algebraic variety (or Deligne-Mumford stack)
equipped with an
algebraic $\com^*$-action. The Atiyah-Bott localization formula expresses
equivariant integrals over $V$ as a sum of contributions over
the $\com^*$-fixed subloci. 

Let $H^*_{\com^*}(V)$ denote the equivariant
cohomology of $V$ with ${\mathbb{Q}}$-coefficients.
Let $H^*_{\com^*}(B\com^*)= {\mathbb{Q}}[t]$
be the standard presentation of the equivariant cohomology ring
of $\com^*$.
The equivariant cohomology ring $H_{\com^*}^*(V)$
is canonically a $H^*_{\com^*}(B\com^*)$-module.
Let
$$H^*_{\com^*}(V)_{[\frac{1}{t}]}= H^*_{\com^*}(V) \otimes {\mathbb{Q}}
[t,\frac{1}{t}]$$
denote the $H^*_{\com^*}(B\com^*)$-module localization
at the element $t\in H^*_{\com^*}(B\com^*)$.

Let $A^{\com^*}_*(V)$ denote the closely related equivariant
Chow ring of $V$ with ${\mathbb{Q}}$-coefficients (defined
in \cite{EdGra,To} via homotopy quotients in the algebraic category). 
$A^{\com^*}_*(V)$ is
a module over $A^*_{\com^*}(B\com^*)= {\mathbb{Q}}[t]$.

Let $\{ V^f_i\}$ be the connected
components of the $\com^*$-fixed locus, and   
let $$\iota: \cup_i V^f_i \rarr V$$ denote the
inclusion morphism.
The nonsingularity of $V$ implies each
$V^f_i$ is also
nonsingular \cite{Iv}. Let $N_i$ denote the normal bundle of
$V^f_i$ in $V$, and let $e(N_i)$ denote the equivariant Euler class
(top Chern class) of $N_i$.

The Atiyah-Bott localization formula \cite{AtBo} is:
\begin{equation}
\label{ssddff}
[V] = \iota_* \sum_{i} \frac{[V^f_i]}{e(N_i)} \ \ \in  H^*_{\com^*}(V)_{[\frac{1}{t}]}
\end{equation}
The formula is well-defined as 
the Euler classes $e(N_i)$ are invertible in the localized equivariant
cohomology ring.

By a result of Edidin-Graham, formula (\ref{ssddff})
holds also in the localized equivariant Chow ring 
$A^{\com^*}_*(V)_{[\frac{1}{t}]}$.

Let $\xi \in H^*_{\com^*}(V)$ be a class of 
degree equal to (twice) the dimension of $V$.
The Bott residue formula \cite{B} expresses integrals
over $V$ in terms of fixed point data:
$$\int_{V} \xi = \sum_{i} \int_{V^f_i} \frac{ \iota^*(\xi)}{e(N_i)}.$$
The Bott residue formula is an immediate consequence of (\ref{ssddff}).
Localization therefore provides an effective method
of computing integrals over $V$ when the fixed loci $V^f_i$
are well-understood.

\subsection{Localization of virtual classes}
Let $V$ be an algebraic variety (or Deligne-Mumford stack) equipped with
a $\com^*$-action. Let $V$ carry a
perfect obstruction theory $\phi: E^\bullet \rarr L^\bullet_V$ equipped
with an equivariant $\com^*$-action. 
While $V$ may be arbitrarily singular, a localization formula
for the virtual class holds.

Let $\{ V^f_i\}$ be the connected
components of the scheme theoretic 
$\com^*$-fixed locus as before.
Since $V$ may be singular, the components $V^f_i$ may be singular
as well. However, each
$V^f_i$ is equipped with a canonical perfect obstruction theory
\cite{GrP}. Moreover, a normal complex 
can be found for each $V^f_i$ (replacing the normal bundle in
the nonsingular case). Together, these constructions 
yield a natural extension of the Atiyah-Bott localization
formula  to virtual classes.

Let $E_i^\bullet$ denote the restriction of the
complex $E^\bullet$ to $V^f_i$. 
The complex $E^\bullet_i$ may be decomposed
by $\com^*$-characters:
$$E^\bullet_i = E^{\bullet,f}_i \oplus E^{\bullet,m}_i,$$
where the first summand corresponds to the
trivial character (the $\com^*$-fixed part) and the
second summand corresponds to all the non-trivial characters
(the $\com^*$-moving part). A canonical morphism
\begin{equation}
\label{ffvff}
\phi_i: E^{\bullet,f}_i \rarr L^\bullet_{V^f_i}
\end{equation}
is obtained from the $\com^*$-fixed part of $\phi$.
It is shown in \cite{GrP} that (\ref{ffvff}) is a perfect
obstruction theory on $V^f_i$.
The $\com^*$-moving part $E^{\bullet,m}_i$ is the
defined to be the virtual (co)normal complex $[N^{vir}_i]^\vee$.

The virtual localization formula \cite{GrP} is:
\begin{equation}
\label{ssddfff}
[V]^{vir} = \iota_* \sum_{i} \frac{[V^f_i]^{vir}}{e(N^{vir}_i)} \ \ \in
A^*_{\com^*}(V)_{[\frac{1}{t}]}.
\end{equation}
The Euler class of $N^{vir}_i= [ E_{0,i}^{m} \rarr E_{1,i}^m]$ is defined
to be:
$$e(N^{vir}_i) = \frac{e(E_{0,i}^m)}{ e(E_{1,i}^m)}.$$
The virtual localization  formula is well-defined  
since the Euler classes of the moving parts of the bundles
$E_{0,i}$ and $E_{1,i}$ are invertible after localization.
The proof of (\ref{ssddfff}) in \cite{GrP} requires the
existence of a $\com^*$-equivariant embedding $V\rarr Y$
in a nonsingular variety (or Deligne-Mumford stack) $Y$.

In case $V$ is nonsingular,
the Atiyah-Bott localization formula is recovered from
(\ref{ssddfff}) via the trivial $\com^*$-equivariant perfect obstruction
theory on $V$.

If the nonsingular target $X$ admits a $\com^*$-action, a canonical
$\com^*$-action by translation is induced on
$\overline{M}_{g,n}(X,\beta)$. 
Stack quotient constructions
prove the existence of 
$\com^*$-equivariant 
nonsingular embeddings for $\overline{M}_{g,n}(X,\beta)$ in this case \cite{GrP}.
The virtual localization formula then
provides an effective tool in the study of 
integrals in Gromov-Witten theory of $X$.

\subsection{Virtual localization for $\overline{M}_{g}(\proj^1,d)$}
\subsubsection{The $\com^*$-action on $\proj^1$}
\label{tacp}
We first establish our $\com^*$-action conventions on $\proj^1$. 
Let $V=\com^2$.
Let $\com^*$ act on $V$ with weights $0,1$: 
\begin{equation}
\label{akkkk} 
t\cdot [v_0,v_1] = [v_0, t v_1].
\end{equation}
The action (\ref{akkkk}) canonically induces a $\com^*$-action
on $\proj^1=\proj(V)$. This action will be fixed throughout the
paper.

We identify
$0,\infty \in \proj^1$ with the $\com^*$-fixed points of $\proj(V)$:
 $$p_0=[1,0],\ p_1=[0,1].$$
The canonical $\com^*$-actions on the tangent spaces to  $\proj(V)$
at $p_0$, $p_1$ have weights $+1$, $-1$ respectively.

\subsubsection{The $\com^*$-action on $\overline{M}_{g}(\proj^1,d)$}
The $\com^*$-action on $\proj^1$ canonically induces
a $\com^*$-action on $\overline{M}_g(\proj^1,d)$ by translation
of maps:
$$t\cdot[\pi] = [t\cdot\pi].$$
As the perfect obstruction theory
of $\overline{M}_{g,n}(\proj^1,d)$ is constructed canonically, 
$\com^*$-equivariance is immediate.

The virtual localization formula is studied here for the
translation
action on $\overline{M}_{g,n}(\proj^1,d)$ following
\cite{GrP}. Four properties of the geometry allow for a
complete analysis of the virtual localization formula:
\begin{enumerate}
\item[(1)] The $\com^*$-fixed locus in $\overline{M}_{g,n}(\proj^1,d)$
is a disjoint union of nonsingular (Deligne-Mumford stack) components.
\item[(2)] Each $\com^*$-fixed component is 
 isomorphic a quotient of products of moduli stacks of
pointed curves $\overline{M}_{\gamma,l}$.
\item[(3)] The virtual structure on the $\com^*$-fixed components
is the canonical trivial structure on a nonsingular space.
\item[(4)] The Euler class of the normal complex is identified
in terms of tautological $\psi$ and $\lambda$ classes on the
fixed components.
\end{enumerate}

\subsubsection{The $\com^*$-fixed components}
\label{ccfixx}
Following \cite{K2}, we can identify the components of the $\com^*$-fixed 
locus of $\overline{M}_{g,n}(\proj^1,d)$ with a set of 
graphs.  We will always assume $d>0$. 

A graph $\Gamma\in G_{g,n}(\proj^1,d)$ consists of the data
$(V,E,N,\gamma,j,\delta)$ where:
\begin{enumerate}
\item[(i)] $V$ is the vertex set,
\item[(ii)] $\gamma: V \rarr {\mathbb{Z}}^{\geq 0}$ is a genus assignment,
\item[(iii)] $j: V \rarr \{0, 1\}$ is a bipartite structure,
\item[(iv)] $E$ is the edge set, 
\begin{enumerate}
\item[(a)]
If the edge $e$ connects $v,v'\in V$, then
$j(v)\neq j(v')$ (in particular, there are no self edges), 
\item[(b)]
$\Gamma$ is connected,
\end{enumerate}
\item[(v)] $\delta: E \rarr {\mathbb{Z}}^{>0}$ is a degree
assignment,
\item[(vi)] $N=\{1, \ldots,n\}$ is a set of markings incident to vertices,
\item[(vii)] $g= \sum_{v\in V} \gamma(v) + h^1(\Gamma)$,
\item[(viii)] $d= \sum_{e\in E} \delta(e).$
\end{enumerate}
The $\com^*$-fixed components of $\overline{M}_{g,n}(\proj^1,d)$
are in bijective correspondence with the graph set $G_{g,n}(\proj^1,d)$.

Let $\pi:(C,p_1,\ldots,p_n) \rarr \proj^1$ be a $\com^*$-fixed stable map.
The images
of all marked points, nodes, contracted components, and ramification
points must lie in the $\com^*$-fixed point set $\{ p_0, p_1\}$ of $\proj^1$.   
In particular, 
each non-contracted irreducible component $D\subset C$  
is ramified only over the two fixed points $\{p_0,p_1\}$.  
Therefore $D$ must be nonsingular and rational. 
Moreover,
the restriction $\pi|_D$ is uniquely determined by 
the degree $deg(\pi|_D)$,  $\pi|_D$ must be the rational Galois cover 
  with full ramification over $p_0$ and $p_1$.

To an invariant stable map $\pi:(C,p_1,\ldots,p_n)\rarr \proj^1$, 
we associate a 
graph $\Gamma \in G_{g,n}(\proj^1,d)$
as follows:
\begin{enumerate}
\item[(i)]
$V$ is the set of connected components of 
$\pi^{-1}(\{ p_0, p_1\})$, 
\item[(ii)] 
$\gamma(v)$ is the arithmetic genus of the component corresponding to
$v$ (taken to be 0 if the component is an isolated point),
\item[(iii)] $j(v)$ is defined by $\pi(v)= p_{j(v)}$,
\item[(iv)] $E$ is the set of non-contracted irreducible components $D\subset C$,
\item[(v)] $\delta(D)= deg(\pi|_D)$,
\item[(vi)] $N$ is the marking set.
\end{enumerate}
Conditions (vii-viii) hold by definition.

The set of $\com^*$-fixed 
stable maps with given graph $\Gamma$ is naturally identified
with a finite quotient of a product of moduli spaces of 
pointed curves.  
Define:  
\begin{equation*}
\M_\Gamma = \prod_{v\in V} \M_{\gamma(v),val (v)}.
\end{equation*}
The valence $val(v)$ is the number of incident edges and markings.
$\M_{0,1}$ and
$\M_{0,2}$ are interpreted as points in this product.
Over $\M_\Gamma$, there is
a canonical universal family of
$\com^*$-fixed stable maps, $$\rho: U \rarr \M_\Gamma,$$
$$\pi: U \rarr \proj^1,$$
yielding a morphism
of stacks $\tau_\Gamma: \M_\Gamma \rarr \overline{M}_{g,n}(\proj^1,d).$

There is a natural automorphism group $\A$ acting 
equivariantly on $U$ and $\M_\Gamma$ with respect to the morphisms $\rho$ and
$\pi$.
$\A$ acts via
automorphisms of the Galois covers (corresponding to the
edges) and the symmetries of the graph $\Gamma$.
$\A$ is filtered by an exact sequence of groups:
$$ 1 \rarr \prod_{e\in E} {\Z}/{\delta(e)} \rarr
\A \rarr \text{Aut}(\Gamma) \rarr 1$$
where $\text{Aut} (\Gamma)$ is the 
automorphism group of $\Gamma$: $\text{Aut}(\Gamma)$ is the
subgroup of the permutation group of the vertices and edges which
respects all the structures of $\Gamma$.
$\text{Aut}(\Gamma)$ naturally acts on $\prod_{\rm edges} \Z/ \delta(e)$
and $\A$ is the semidirect product. 

Let $Q_\Gamma$ denote the quotient stack $\M_\Gamma/\A$.
The induced 
map:
$$\tau_\Gamma/ \A : Q_\Gamma \rarr \overline{M}_{g,n}(\proj^1,d)$$
is a closed immersion of Deligne-Mumford stacks.
It should be noted that the subgroup $\prod_{\rm edges} \Z/ \delta(e)$
acts trivially on $\M_\Gamma$. $\Q_\Gamma $
is a nonsingular Deligne-Mumford stack.

The above {\em set-theoretic} analysis proves 
a component of the $\com^*$-fixed stack of $\overline{M}_{g,n}(\proj^1,d)$
is {\em supported} on the substack $Q_\Gamma$.

\subsubsection{The $\com^*$-fixed perfect obstruction theory}
Let $\phi: E^\bullet \rarr L^\bullet_{\overline{M}_{g,n}(\proj^1,d)}$
denote the $\com^*$-equivariant perfect obstruction theory of the moduli
of maps.
Let $E_{\bullet,\Gamma}$ denote the restriction of $E_\bullet$ to
$Q_\Gamma$. 
Denote the cohomology of $E_{\bullet,\Gamma}$ by:
\begin{equation}
\label{plml}
0 \rarr {\text{Tan}} \rarr E_{0,\Gamma} \rarr E_{1,\Gamma} \rarr {\text{Obs}} 
\rarr 0.
\end{equation}
The tangent-obstruction
sequence may be studied on $Q_\Gamma$
--- the sequence is obtained from the cohomology of the
(dual of) the restriction to $Q_\Gamma$  of the top distinguished triangle of
(\ref{hdm}). The fiber of the tangent-obstruction sequence over 
$[\pi]\in Q_\Gamma$ is:
\begin{equation}
\label{totom}
0 \rarr \text{Ext}^0(\Omega_C(P), \oh_C) 
\rarr H^0(C,\pi^*T_{\proj^1}) \rarr {\text{Tan}}
\end{equation}
$$ \ \ \ \rarr \text{Ext}^1(\Omega_C(P), \oh_C) \rarr H^1(C,\pi^*T_{\proj^1}) 
\rarr \text{Obs} \rarr 0.$$
The elements of (\ref{totom}) are 
 {\em vector bundles} on $Q_\Gamma$ (instead of
possibly singular sheaves) as theirs ranks are constant on $[\pi]\in Q_\Gamma$.

The scheme structure of the $\com^*$-fixed stack supported
on $Q_\Gamma$ may be determined from the perfect obstruction theory.
The Zariski tangent space at $[\pi]$ to the $\com^*$-fixed stack is 
$\text{Tan}^f_{[\pi]}$. A direct study of the $\com^*$-fixed part of
(\ref{totom}) in \cite{GrP} shows this Zariski tangent space to
be isomorphic to the tangent space of $Q_\Gamma$.
As $Q_\Gamma$ is a nonsingular stack, we may conclude
the $Q_\Gamma$ {\em is} a component of the
$\com^*$-fixed stack.

The second use of (\ref{totom}) is to determine the
perfect obstruction theory of the $\com^*$-fixed component $Q_\Gamma$
induced by $\phi$. An analysis of the $\com^*$-fixed
part of (\ref{totom}) immediately implies the induced
perfect obstruction theory is trivial \cite{GrP}. It is quite easy to
analyze the sequence (\ref{totom}) as the stable maps
parameterized by $Q_\Gamma$ are of a uniformly simple character.

The virtual localization formula for $\overline{M}_{g,n}(\proj^1,d)$
may now be written as:
\begin{equation}
\label{ssddffff}
[\overline{M}_{g,n}(\proj^1,d)]^{vir} = \sum_{\Gamma\in G_{g,n}(\proj^1,d)}
\frac{1}{|A_\Gamma|} \frac{\tau_{\Gamma *}[\M_\Gamma]}{e(N_\Gamma^{vir})} 
\end{equation}
in
$A_*^{\com^*}(\overline{M}_{g,n}(\proj^1,d))_{[\frac{1}{t}]}$.
The $\com^*$-fixed loci $Q_\Gamma$ enter (\ref{ssddffff}) as
push-forwards of $\M_\Gamma$ via $\tau_\Gamma$.

\subsubsection{The normal complex}
The tangent-obstruction sequence
(\ref{totom}) also determines the
Euler class of the normal complex of $\com^*$-fixed loci
induced by $\phi$. 
The moving parts of the vector bundle
sequences (\ref{plml}-\ref{totom}) imply:
\begin{equation}
\label{weddd}
\frac{1}{e(N^{vir})} =  \frac{e( \text{Ext}^0(\Omega_C(P), \oh_C)^m)}
{e(\text{Ext}^1(\Omega_C(P), \oh_C)^m)} \cdot \frac{ 
e(H^1(C,\pi^*T_{\proj^1})^m )} 
{e(H^0(C,\pi^*T_{\proj^1})^m)}.
\end{equation}

Let $\Gamma\in G_{g,n}(\proj^1,d)$. 
The above identification (\ref{weddd}) precisely specifies 
the $\tau_\Gamma$ pull-back of $1/e(N^{vir})$ 
to $\M_\Gamma$,
\begin{equation}
\label{xccv}
\M_\Gamma = \prod_{v\in V} \overline{M}_{\gamma(v), val(v)}.
\end{equation}
The pull-backs of the Euler classes of the vector bundles
on the right of (\ref{weddd}) naturally split
over the vertex factors of $\Gamma$. We will find:
\begin{equation}
\label{ssggv}
\tau^*_\Gamma( \frac{1}{e(N^{vir})} ) =
(-1)^d \prod_{v\in V} \frac{1}{N(v)},
\end{equation}
where the {\em vertex contributions} $1/N(v)$ lie in localized equivariant
cohomology,
$$\frac{1}{N(v)}\in A_*^{\com^*} (\overline{M}_{\gamma(v),val(v)})_{[\frac{1}{t}]}.$$

\subsubsection{Vertex contributions}
Intermediate vertex {\em and} edge contributions $1/\tilde{N}(v)$ and $1/\tilde{N}(e)$
naturally arise
in the geometric analysis of the
(\ref{weddd}). The intermediate contributions will be joined to
yield the single vertex contribution $1/N(v)$.

There are four types of vertices which we treat independently here. In 
integration formulas, a uniform treatment of the four types is often found.

 A vertex $v$ is {\em stable} if
$2\gamma(v)-2+ val(v) >0.$
If $v$ is stable, the factor $\overline{M}_{\gamma(v), val(v)}$
is a factor of $\overline{M}_\Gamma$ by (\ref{xccv}). The intermediate
contribution
$1/\tilde{N}(v)$ will be a equivariant cohomology class on the factor $\overline{M}_{\gamma(v),val(v)}$ in this case.

\vspace{10pt}
$\bullet$ Let $v$ be a stable vertex.
Let $e_1, \ldots, e_l$ denote the distinct edges incident to $v$ (in bijective
correspondence to a subset of the (local) markings of the moduli
space $\overline{M}_{\gamma(v), val(v)}$).
Let $\psi_i$ denote the cotangent line of the marking corresponding to $e_i$.
\begin{eqnarray*} 
\frac{1}{\tilde{N}(v)} & = & \prod_{i=1}^l \frac{1}
{\frac{(-1)^{j(v)} t}{\delta(e_i)} - \psi_i} 
\cdot \\
 & & ((-1)^{j(v)} t)^{l-1} \cdot \\
& & \sum_{i=0}^{\gamma(v)} (-1)^i ((-1)^{j(v)} t)^{\gamma(v)-i} \lambda_i
.
\end{eqnarray*}
The three factors in $1/\tilde{N}(v)$ 
are the contributions of $\text{Ext}^1(\Omega_C(P), \oh_C)^m$, $H^0(C,\pi^*T_X)^m$, and
$H^1(C,\pi^*T_X)^m$ respectively.  $\text{Ext}^0(\Omega_C(P), \oh_C)^m$  
does not contribute to stable vertices.

We note both the tautological $\psi$ and $\lambda$ classes enter in 
$1/\tilde{N}(v)$.
The Gromov-Witten theory of $\proj^1$
is therefore 
fundamentally related to the intersection theory of the moduli space of
curves.

If $v$ is an unstable vertex, then $\gamma(v)=0$ and $val(v) \leq 2$.
There are three unstable cases: 
two with valence 2 and one with valence 1.

\vspace{10pt}
\noindent $\bullet$ Let $v$ be an {\em unmarked} vertex with
$\gamma(v)=0$ and $val(v)=2$. Let $e_1$ and $e_2$ be the two
incident edges. Then:
$$\frac{1}{\tilde{N}(v)} 
= \frac{1}{\frac{(-1)^{j(v)}t}{\delta(e_1)} + \frac{ (-1)^{j(v)} t}{\delta(e_2)}}
\cdot (-1)^{j(v)}t = \frac{1}{\frac{1}{\delta(e_1)} + \frac{1}{\delta(e_2)}}.$$
The factors are obtained from  
$\text{Ext}^1(\Omega_C(P), \oh_C)^m$ and
$H^0(C,\pi^*T_{\proj^1})^m$ respectively.

\noindent $\bullet$ Let $v$ be a $1$-marked vertex with
$\gamma(v)=0$ and $val(v)=2$. Let $e$ be the unique
incident edge. Then:
$$\frac{1}{\tilde{N}(v)} = 1,$$
there are no contributing factors.

\noindent $\bullet$ Let $v$ be an {\em unmarked} vertex with
$\gamma(v)=0$ and $val(v)=1$. Let $e$ be the unique
incident edge. Then:
$$\frac{1}{\tilde{N}(v)} = \frac{ (-1)^{j(v)} t}{\delta(e)},$$
where $\text{Ext}^0(\Omega_C(P), \oh_C)^m$ is the only contributing factor.

All of these contributions are easily extracted from an analysis of
(\ref{weddd})  \cite{GrP}.

\subsubsection{Edge contributions}
Let $e\in E$ be an edge corresponding to the non-contracted irreducible
component $D\subset C$ (where $$[\pi: (C,p_1,\ldots,p_n)\rarr \proj^1]$$
is a moduli point parameterized by $\overline{M}_\Gamma$). The edge contribution,
$$\frac{1}{\tilde{N}(e)} \in A_*^{\com^*}(B\com^*)_{ 
[\frac{1}{t}]},$$ is the inverse Euler class of 
the $\com^*$-representation
$H^0(D, \pi^* T_{\proj^1})^m.$  The contribution is
obtained from $H^0(C,\pi^*T_{\proj^1})^m$.

Consider the $\com^*$-equivariant 
Euler sequence
on $\proj^1$:
$$0 \rarr \oh \rarr \oh(1) \otimes V \rarr T\proj^1 \rarr 0. $$
After pulling back to $D$ and taking cohomology, we find:
\begin{equation}
\label{mnbv}
0 \rarr \com \rarr H^0(D,\oh(\delta(e)))\otimes V \rarr H^0(D,\pi^*T\proj^1) \rarr 0. 
\end{equation}
The $\com^*$-weight on $\com$ is trivial, 
and the weights of $H^0(D, \oh(\delta(e)))$
are: $$-\frac{i t }{\delta(e)}, \ \ \ 0\leq i \leq \delta(e).$$
The weights of $V$ are $0, 1$.  The
weights of the of the 
middle term in (\ref{mnbv}) are therefore the pairwise sums: 
$$-\frac{it }{\delta(e)}, \ 1-\frac{it }{\delta(e)}, \ \ \ 0\leq i \leq \delta(e).$$  
As only the moving weights concern us, we find:
$$\frac{1}{\tilde{N}(e)}= \frac{1}{ (-1)^{\delta(e)} \frac{{\delta(e)!}^2 }{ \delta(e)^{2\delta(e)}}
t^{2\delta(e)}}.$$

By the analysis of \cite{GrP},
the contributions $1/\tilde{N}(v)$ and $1/\tilde{N}(e)$
together account for the entire right side of (\ref{weddd}). We find:
\begin{equation}
\label{xbh}
\tau^*_\Gamma( \frac{1}{e(N^{vir})} ) =
\prod_{v\in V} \frac{1}{\tilde{N}(v)} \cdot 
\prod_{e\in E} \frac{1}{\tilde{N}(e)}.
\end{equation}

\subsubsection{$1/N(v)$}
\label{convertl}
Since
the intermediate edge contribution $(-1)^{\delta(e)}\tilde{N}(e)$ admits a square root,
$$\sqrt{ \frac{(-1)^{\delta(e)}} { 
\tilde{N}(e)}} = \frac{\delta(e)^{\delta(e)}}{\delta(e)!} t^{-\delta(e)},$$
the edge contributions may be distributed to the incident vertices.
Let $v$ be a vertex with  incident edges $e_1, \ldots, e_l$.
Define $1/N(v)$ by:
$$\frac{1}{N(v)} = 
\frac{1}{\tilde{N}(v)} 
\cdot \prod_{i=1}^l \frac{\delta(e_i)^{\delta(e_i)}}{\delta(e_i)!} t^{-\delta(e_i)}.$$
Equation (\ref{xbh}) then immediately implies (\ref{ssggv}).

\subsubsection{Integration}
Virtual localization yields an integration formula for the
Gromov-Witten theory of $\proj^1$.
The expected dimension of the moduli space $\overline{M}_{g,n}(\proj^1,d)$ is
$2g-2+2d+n$.
Let $\xi$ be an equivariant class
 $$\xi \in H^{2(2g-2+2d+n)}_{\com^*}(\overline{M}_{g,n}(\proj^1,d), {\mathbb{Q}}).$$
Via the canonical morphism,
$$H^*_{\com^*} (\overline{M}_{g,n}(\proj^1,d), {\mathbb{Q}})
\rarr H^*(\overline{M}_{g,n}(\proj^1,d), {\mathbb{Q}}),$$
The class $\xi$ may be viewed as an equivariant lift of an ordinary 
cohomology class on $\overline{M}_{g,n}(\proj^1,d)$ --- called the
{\em non-equivariant limit} of $\xi$.

The virtual residue formula 
for the integral of $\xi$ obtained from virtual localization is:
\begin{equation}
\label{bigint}
\int_{[\overline{M}_{g,n}(\proj^1,d)]^{vir}} \xi =
\sum_{\Gamma \in G_{g,n}(\proj^1,d)} \frac{(-1)^d}{|\A|}
\int_{\overline{M}_\Gamma} \frac{\tau_\Gamma^*(\xi)}{ \prod_{v\in V} N(v)}.
\end{equation}
The left side of (\ref{bigint}) is equal to the integral of 
the non-equivariant limit of $\xi$. Only the $t^0$ terms contribute to
the right side after integration.

Formula (\ref{bigint}) effectively relates integrals in the
Gromov-Witten theory of $\proj^1$ to tautological integrals over
the moduli space of curves.

\subsection{Gravitational descendents}
We explain here an application of the virtual localization formula to
the descendent invariants of $\proj^1$:
\begin{equation}
\label{neddd}
\langle  \prod_{i=1}^r \tau_{a_i} \cdot \prod_{j=r+1}^{r+s}
\tau_{b_j}(\omega) \rangle_{g,d}^{\proj^1} =
\int_{[\overline{M}_{g,n}(\proj^1,d)]^{vir}} \prod_{i=1}^r \psi_{i}^{a_i} \cdot
\prod_{j={r+1}}^{r+s} \psi_j^{b_j} \text{ev}_j^*(\omega).
\end{equation}

All terms of integrand of (\ref{neddd}) are equipped with canonical 
$\com^*$-equivariant lifts. First, the $\com^*$-action is canonically
lifted to the cotangent classes $\psi_i$ of $\overline{M}_{g,n}(\proj^1,d)$. 
Second,
the class $\omega=c_1(\oh(1))$
 is canonically lifted to
$H_{\com^*}^2(\proj^1, {\mathbb{Q}})$ via the canonical $\com^*$-action
on $\oh(1)$ --- the $\com^*$-action on $V$ induces an
action on the tautological line $\oh(-1)$ and (by dualizing)
an action on $\oh(1)$. The $\com^*$-action
on $\oh(1)$ has fiber weights $$w_0=0, \ w_1=-1$$ over the points $p_0,p_1 \in \proj^1$
respectively. Finally,
the class $\text{ev}_j^*(\omega)$ may 
be canonically lifted from the lift of $\omega$.

Let $\xi$ denote the canonical lift of the integrand of (\ref{neddd}).
The virtual localization formula applied to 
$\xi$ determines the
descendent invariant in terms of tautological integrals over the
moduli spaces of curves:
\begin{equation*} 
\langle  \prod_{i=1}^r \tau_{a_i} \cdot \prod_{j=r+1}^{r+s}
\tau_{b_j}(\omega) \rangle_{g,d}^{\proj^1}
=
\sum_{\Gamma \in G_{g,n}(\proj^1,d)} \frac{(-1)^d}{|\A|}
\int_{\overline{M}_\Gamma} \frac{ \tau_\Gamma^*(\xi)}{ \prod_{v\in V} N(v)}.
\end{equation*}

The pull-back of $\xi$ to $\overline{M}_\Gamma$ factorizes over the
vertices of $\Gamma$:
$$\tau_\Gamma^*(\xi) = \prod_{v\in V} \xi(v).$$
There are four types of vertex contributions $\xi(v)$.

\vspace{10pt}
\noindent $\bullet$  Let $v$ be a stable vertex. Let $\{1,\ldots,r+s\}$ denote
the (global) marking set of $\overline{M}_{g,r+s}(\proj^1,,d)$.
Let $$R\subset \{ 1,\ldots,r\}, \ \
S\subset \{r+1, \ldots r+s\}$$
denote the subsets of the global markings lying on $v$.
Then,
$$\xi(v)= \prod_{i\in R} \psi_i^{a_i} \cdot \prod_{i\in S} \psi_i^{b_i} w_{j(v)} t \ \
\in H^*_{\com^*} (\overline{M}_{\gamma(v), val(v)}).$$
Note this contribution vanishes if $j(v)=0$ and $S$ is non-empty.

\noindent $\bullet$ Let $v$ be an unmarked vertex with $\gamma(v)=0$ and $val(v)=2$.
Then,
$$\xi(v)=1.$$

\noindent $\bullet$ Let $v$ be a $1$-marked vertex with $\gamma(v)=0$ and $val(v)=2$.
Let $e$ denote the unique edge incident to $e$.
If the marking $i$ of $v$ satisfies $1\leq i \leq r$, then
$$\xi(v) = \Big(  -\frac{(-1)^{j(v)}t}{\delta(e)} \Big)^{a_i}.$$
If the marking $i$ of $v$ satisfies $r+1 \leq i \leq r+s$, then
$$\xi(v)= \Big(  -\frac{(-1)^{j(v)}t}{\delta(e)} \Big)^{b_i} w_{j(v)} t.$$
Note the second contribution vanishes if $j(v)=0$.

 \noindent $\bullet$ Let $v$ be an unmarked vertex with $\gamma(v)=0$ and $val(v)=1$.
Then,
$$\xi(v)=1.$$

We find an explicit formula for the gravitational descendent invariants
of $\proj^1$ in terms of tautological integrals over the
moduli space of curves.

\begin{pr} The gravitational descendents of $\proj^1$ are determined
by graph sums of Hodge integrals:
$$\langle  \prod_{i=1}^r \tau_{a_i} \cdot \prod_{j=r+1}^{r+s}
\tau_{b_j}(\omega) \rangle_{g,d}^{\proj^1}
=
\sum_{\Gamma \in G_{g,n}(\proj^1,d)} \frac{(-1)^d}{|\A|}
\int_{\overline{M}_\Gamma}  \prod_{v\in V} \frac{\xi(v)}{ N(v)}.$$
\end{pr}

\section{\bf From Hurwitz numbers to Hodge integrals}
\label{fin1}
\subsection{The proof of Theorem \ref{rrrr}}
The Hurwitz numbers $H_{g,\mu}$ count genus $g$ covers of
$\proj^1$ with profile $\mu$ over $\infty$ and simple ramification
over a fixed set of finite points.
The relationship between 
Hurwitz numbers and Hodge integrals is proven here
via the Gromov-Witten theory of $\proj^1$.

The proof of Theorem \ref{rrrr} is immediate in case
$\mu$ is trivial, the case of the Hurwitz numbers $H_{g,d}$.
The Hurwitz numbers $H_{g,d}$ arise as integrals against
$[\overline{M}_{g}(\proj^1,d)]^{vir}$ via the branch morphism. 
The Hodge integral relationship is then
a direct consequence of the virtual residue formula.
The argument for $H_{g,d}$ is explained
first in Section \ref{hurr1}.

Theorem \ref{rrrr} is proven 
for arbitrary profile $\mu$ in Section \ref{hurr2}.
The Hurwitz numbers $H_{g,\mu}$ arise as 
integrals over natural {\em components} 
of $\overline{M}_{g}(\proj^1,d)$. 
A detailed analysis 
is required to extract the relevant component contributions from
the virtual residue formula \cite{GrV}.
Our presentation in Section \ref{hurr2} follows \cite{GrV}.

\subsection{The Hurwitz number $H_{g,d}$}
\label{hurr1}
\subsubsection{Integrals}
The Hurwitz number $H_{g,d}=H_{g,1^d}$ counts genus $g$ covers of $\proj^1$
\'etale over $\infty$ with  
$$r=2g-2+2d$$ fixed finite simple ramification points.
The branch morphism $br$ constructed in Section \ref{brrrch} is:
\begin{equation*}
\br: \overline{M}_{g}(\proj^1,d) \rarr \text{Sym}^{r}(\proj^1).
\end{equation*}
Let ${\xi_p}$ denote (the Poincar\'e dual of) the 
point class of $\text{Sym}^r(\proj^1)$.

\begin{pr} The Hurwitz number $H_{g,d}$ is an integral in Gromov-Witten
theory:
\begin{equation*}
H_{g,d}= \int_{ [\overline{M}_{g}(\proj^1,d)]
^{vir}} \br^*(\xi_p).
\end{equation*}
\end{pr}

\bpf The locus $M_g(\proj^1,d) 
\subset \overline{M}_g(\proj^1,d)$
is nonsingular (of the expected dimension) by Proposition \ref{fppp}.

Let $z_1, \ldots, z_r\in \proj(V)$ be distinct points.
If $[\pi:C\rarr \proj^1]$ is a stable map with a singular domain
curve, then the divisor $\br(\pi)$ must contain a double point. 
Therefore,
$\br^{-1}(\sum_{i=1}^r [z_i]) \subset M_g(\proj^1,d)$.
By Bertini's Theorem applied to the morphism
$$\br:M_g(\proj^1,d) \rarr \text{Sym}^r(\proj^1)=\proj^r,$$
a general divisor $\sum_{i=1}^r [z_i]$
intersects the stack $M_{g}(\proj^1,d)$ transversely via $\br$ in
a finite number of points. These intersections are exactly the finitely many
Hurwitz covers $H_{g,d}$ ramified over $\{z_i\}$ (weighted by 
$1/ |\text{Aut}|$ in the intersection product).
\epf

\subsubsection{Localization}
We follow the conventions set in Section \ref{tacp} regarding
the $\com^*$-action on $\proj^1=\proj(V)$.

The canonical $\com^*$-actions  on the spaces 
$\overline{M}_g(\proj^1,d)$ and 
$\text{Sym}^r(\proj^1)$ are $\br$-equivariant by the 
canonical construction of the branch morphism \cite{FanP}.

Let $\xi$ be the $\com^*$-equivariant lift of the point class $\xi_p$
corresponding to the $\com^*$-fixed divisor $r[p_0] \in \text{Sym}^r(\proj(V))$.
The integral,
\begin{equation*}
H_{g,d}= \int_{ [\overline{M}_{g}(\proj^1,d)]^{vir}} \br^*(\xi),
\end{equation*}
may then be evaluated via 
the virtual residue formula:
\begin{equation}
\label{fforr}
H_{g,d}= 
\sum_{\Gamma\in G_{g}(\proj^1,d)} 
\frac{(-1)^d}{|\A|} \int_{\overline{M}_\Gamma} 
\frac{ \br^*(\xi)} {\prod_{v\in V} N(v)}.
\end{equation}

$\text{Sym}^r(\proj^1)$ has $r+1$ isolated fixed points:
$(r-a)[p_0] + a[p_1]$, for $0\leq a \leq r$.
For each graph $\Gamma$, 
the morphism $\br$ contracts $\overline{M}_\Gamma$
to a fixed point of $\text{Sym}^r(\proj^1)$. 
Therefore, $\br^*(\xi)|_{\overline{M}_\Gamma}=0$ unless
$\br(\overline{M}_\Gamma)= r[p_0].$

Let $[\pi:C \rarr \proj^1]$ be a stable map such that $\br(\pi)= r[p_0]$. 
All nodes, collapsed
components, and ramifications of $\pi$ must lie over $p_0$.
Hence, if $\br(\overline{M}_\Gamma)=r[p_0]$, the graph $\Gamma$ may not have
any vertices of positive genus or valence greater than 1 lying
over $p_1$. Moreover, the degrees of the edges of $\Gamma$
must all be 1.

Exactly one
graph $\Gamma_0$ satisfies $\br(\overline{M}_\Gamma)=r[p_0]$.
$\Gamma_0$ is determined by the following construction.
$\Gamma_0$ has a unique genus $g$ vertex $v_0$ lying over $p_0$ 
which
is incident to exactly
$d$ degree 1 edges. The edges connect $v_0$ to
$d$ 
unstable, unmarked vertices $v_1^1,\ldots,v_1^d$
of valence 1 and genus 0 lying over $p_1$. 

By definition,
$\overline{M}_{\Gamma_0}= \overline{M}_{g,d}$.
Since the automorphism group of $\Gamma_0$
is the full permutation group of the edges, $|A_{\Gamma_0}|= d!$.
The vertex contributions of the Euler class of the
normal complex were found in Section \ref{convertl}:
$$\frac{1}{N(v_0)}= \frac{t^g-t^{g-1}\lambda_1+t^{g-2}\lambda_2 
- t^{g-3}\lambda_3 + \ldots + (-1)^g \lambda_g}
{\prod_{i=1}^d (t-\psi_i)}  t^{-1},$$
for the unique vertex over $p_0$ and
$$\frac{1}{N(v_1^i)}= -1,$$
for each of the $d$ unstable vertices over $p_1$.

By the excess intersection formula, the class
$\br^*(\xi)|_{\overline{M}_\Gamma}$ is the $\com^*$-equivariant Euler class of the
normal bundle of the point $r[p_0]$ in $\text{Sym}^r(\proj(V))$:
$$\br^*(\xi)|_{\overline{M}_{\Gamma_0}} =r!\ t^r,$$
easily computed, for example, via the canonical isomorphism
$$\text{Sym}^r(\proj^1) = \proj(\text{Sym}^r V^*).$$

The sum (\ref{fforr}) contains only one term:
\begin{equation*} 
H_{g,d}= 
\frac{(-1)^d}{|A_{\Gamma_0}|} \int_{\overline{M}_{\Gamma_0}}
\frac{ 
\br^*(\xi)}{ \prod_{v\in V} N(v)}.
\end{equation*}
After substitution of the identified factors, we find:

\vspace{10pt}
\noindent{\bf Theorem \ref{rrrr}}. (For $H_{g,d}$).
$$H_{g,d}= \frac{(2g-2+2d)!}{d!} \int_{\overline{M}_{g,d}}
\frac{1-\lambda_1+\lambda_2 -\lambda_3 + \ldots + (-1)^g \lambda_g}
{\prod_{i=1}^d (1-\psi_i)},$$
for $(g,d) \neq (0,1), (0,2)$.
\vspace{10pt}

The genus 0 formula,
\begin{equation}
\label{gg00}
H_{0,d}= \frac{(2d-2)!}{d!} d^{d-3},
\end{equation}
immediately follows from Theorem \ref{rrrr} together with the
evaluations:
$$\int_{\overline{M}_{0,n}} \psi_1^{a_1} \cdots \psi_n^{a_n} =
\binom{n-3}{a_1,\ldots, a_n}.$$
Equation (\ref{gg00}) was first found by Hurwitz.

\subsection{The Hurwitz number $H_{g,\mu}$}
\label{hurr2}
\subsubsection{Overview}
The proof of Theorem \ref{rrrr} for $H_{g,\mu}$ requires 
a study of maps with fixed profile over $\infty$.
However,
the strategy of Section \ref{hurr1} is maintained. The Hurwitz number $H_{g,\mu}$
is first
identified as an integral over a restricted moduli space of maps.
Then, Theorem \ref{rrrr} is deduced from a vertex contribution via
the virtual residue formula.
The presentation here follows \cite{GrV}.
\subsubsection{Moduli spaces and integrals}
Let $C$ be a nonsingular genus $g$ curve.
Let $\pi:C \rarr \proj^1$ be a map with  
profile $\mu=(m_1,\ldots,m_l)$ over $p_1=\infty$.
Let $d=|\mu|$ be the degree of $\pi$.
Let $r= 2g-2+d+l$ be the number of simple ramifications of $\pi$
over finite points. 
Let $k=\sum_{i} (m_i-1) = d-l$.
The branch morphism is:
$$\br: \overline{M}_g(\proj^1,d) \rarr \text{Sym}^{r+k}(\proj^1) =
\proj^{r+k}.$$

Let $L_k$ denote the {\em linear} subspace of $\text{Sym}^{r+k}(\proj^1)$
defined by:
$$L_k = \{\ D+ k[p_1] \ | \ D \in \text{Sym}^r(\proj^1) \}.$$ 
As $\pi$ has profile $\mu$ over $p_1$, the branch divisor
satisfies $\br(\pi)\in L_k$.
Define $\overline{M}_g(L_k)$ by the $\com^*$-equivariant fiber square:
\begin{equation}
\label{hdmm}
\begin{CD}
\overline{M}_g(L_k)
  @>>> \overline{M}_{g}(\proj^1,d) \\
@V{\br_k}VV   @V{\br}VV  \\
 L_k
    @>{\iota}>> \text{Sym}^{r+k}(\proj^1).
\end{CD}
\end{equation}
A virtual class of dimension $r$ is induced on $\overline{M}_{g}(L_k)$
by the Gysin map:
$$[\overline{M}_{g}(L_k)]^{vir} = \iota^![\overline{M}_g(\proj^1,d)]^{vir}.$$
Theorem \ref{rrrr} is proven by virtual localization on
$\overline{M}_g(L_k)$.

As before, let $M_g(\proj^1,d)$ be the
open moduli space of maps with nonsingular domains.
By Proposition \ref{fppp}, $M_g(\proj^1,d)$ is
a nonsingular Deligne-Mumford stack of pure dimension
$r+k$. 
Let $M_{g}(\mu) \subset M_{g}(\proj^1,d)$ denote the
(reduced) substack of maps with profile
$\mu$ over $p_1$.
$M_g(\mu)$ is of pure dimension $r$.
Let $$M_g(\mu) \subset \overline{M}_g(\mu)$$ denote the
closure.

$\overline{M}_{g}(\mu)$ is a substack of $\overline{M}_g(L_k)$ equal
to a union of irreducible components.
The restricted branch divisor is well-defined:
$$\br_k = \overline{M}_g(\mu) \rarr L_k.$$
Let $\xi_p$ denote (the Poincar\'e dual of) the point class of
$L_k$.

\begin{pr} 
\label{pjjl}
The Hurwitz number $H_{g,\mu}$ is an integral:
$$H_{g,\mu} = \int_{[\overline{M}_g(\mu)]} \br_k^*(\xi_p).$$
\end{pr}

\bpf The integral is well-defined as $\overline{M}_g(\mu)$
is of pure dimension $r$. By Bertini's Theorem,
a general point $\sum_{i=1}^r [z_i] +k[p_1]$ of $L_k$
intersects the stack $\overline{M}_{g}(\mu)$ transversely via $\br_k$ in
a finite number of nonsingular points of $M_{g}(\mu)$.
These intersections are exactly the finitely many
Hurwitz covers $H_{g,\mu}$ simply ramified over $\{z_i\}$ (weighted by 
$1/ |\text{Aut}|$ in the intersection product).
\epf

\subsubsection{Multiplicity}
\label{kkrr}
The moduli space $M_{g}(\mu)\subset M_{g}(\proj^1, d)$
occurs as an open {\em set} of the intersection 
$\br^{-1}(L_k) \cap M_g(\proj^1,d)$.
The multiplicity of $\br^{-1}(L_k)\cap M_g(\proj^1,d)$ along $M_g(\mu)$
will be required in the proof of Theorem \ref{rrrr}.

\begin{lm}
\label{mmull}
The intersection $\br^{-1}(L_k) \cap M_g(\proj^1,d)$ is of uniform
multiplicity 
$$\text{mult}(\mu)= k! \prod_{i=1}^l \frac{m_i^{m_i-1}}{m_i!}$$
along $M_g(\mu)$.
\end{lm}

\bpf
Let $m\leq r+k$.
Let $x_1,\ldots,x_m$ be distinct points of $\proj^1$.
Define the linear space $L(x_1,\ldots,x_m)\subset \text{Sym}^{r+k}(\proj^1)$ by:
$$L(x_1,\ldots,x_m) = 
\{\ D+ \sum_{i=1}^m[x_i] \ | \ D \in \text{Sym}^{r+k-m}(\proj^1) \}.$$

Let $[\pi]\in M_g(\mu)$ be a map with simple ramification over
the points $z_1,\ldots z_r \in \proj^1$.
Assume 
the linear space $L(z_1,\ldots,z_r)$  
intersects $M_g(\mu)$ transversely via $\br$ at nonsingular
reduced points (including $[\pi]$). The assumption holds for all
$[\pi]$ in a dense open subset of $M_g(\mu)$
by Bertini's Theorem.

Let $\{z'_j(s)\}_{j=1}^k$ be holomorphic paths in $\proj^1$ satisfying:
\begin{enumerate}
\item[(i)] $z'_{j_1}(s)\neq z'_{j_2}(s)$, for all $j_1\neq j_2$ and $0\neq s\in \com$,
\item[(ii)] $z'_j(0)=p_1$, for all $j$.
\end{enumerate}
The  substacks $\br^{-1}(L(z'_1(s),\ldots,z'_k(s)))\cap M_g(\proj^1,d)$ 
form a flat family specializing to $\br^{-1}(L_k) \cap M_g(\proj^1,d)$
at $s=0$.

For all except finitely many special values of $s$,
$z_i \neq z'_j(s).$ At nonspecial values,
 $L(z_1,\ldots,z_r)$  
intersects $\br^{-1}(L(z'_1(s),\ldots,z'_k(s)))$ transversely via $\br$ at nonsingular
reduced points corresponding to $H_{g,d}$ Hurwitz covers
with simple ramification over $\{z_i\}\cup \{z'_j(s)\}$.

Let $H(s)$ denote the set of the Hurwitz covers specified by $s$.
Let $H(\pi) \subset H(s)$ be the subset of Hurwitz covers
which specialize to $[\pi]$ as $s\rarr 0$.
The multiplicity of $\br^{-1}(L_k)$ at $[\pi]$ is 
equal to $|H(\pi)|$.

$H(s)$ is equal to the set of $(r+k)$-tuples
of 2-cycles $$(\gamma_1,\ldots,\gamma_r,\gamma'_1,\ldots\gamma'_k)$$
modulo $S_d$-conjugation satisfying:
\begin{enumerate}
\item[(a)] $\gamma_1, \ldots, \gamma_r, \gamma'_1,\ldots,\gamma'_k$ generate 
a transitive subgroup of $S_{d}$,
\item[(b)] $\prod_{i=1}^r \gamma_i\prod_{j=1}^k \gamma'_j =1$.
\end{enumerate}

Let $c_{m_1}\cdots c_{m_l}\in S_d$ be a fixed element with
cycle decomposition $\mu$. 
The elements $H(\pi)\subset H(s)$ 
bijectively correspond
to solutions of the equation:
\begin{equation}
\label{mull}
\prod_{j=1}^k\gamma'_j=c_{m_1} \cdots c_{m_l}.
\end{equation}
The number of solutions of (\ref{mull}) 
is proven to equal
$$k! \prod_{i=1}^l \frac{m_i^{m_i-1}}{m_i!}$$
in Lemma \ref{solomon} below.
\epf

\begin{lm}
\label{solomon}
The equation $\prod_{j=1}^k\gamma'_j=c_{m_1} \cdots c_{m_l} \in S_d$
has 
$$k! \prod_{i=1}^l \frac{m_i^{m_i-1}}{m_i!}$$
solutions for $k$-tuples $(\gamma_1',\ldots, \gamma_k')$.
\end{lm}

\bpf
A $2$-cycle $(x_1x_2)$ {\em lies in the span} of a
cycle $c=(y_1\cdots y_m)$ if $$\{x_1,x_2\}\subset \{ y_1,\ldots,y_m\}.$$
Each solution of 
\begin{equation}
\label{notri} 
\prod_{j=1}^k\gamma'_j=c_{m_1} \cdots c_{m_l}
\end{equation}
has the following
property: for each $i$, exactly $m_i-1$ of the $2$-cycles
$\gamma'_j$ lie in the span of $c_{m_i}$.

An elegant proof of the above property is given in \cite{GrV}.
A solution of (\ref{mull}) defines a degree $d$ cover $D\rarr \proj^1$ 
with simple ramifications determined
by $\{\gamma_j'\}$ at fixed finite points $q_1,\ldots q_k$
and profile $\mu$ over $p_1$.
The arithmetic genus of $D$ is $1-l$ by the Riemann-Hurwitz formula.
As the preimage of $p_1$ contains $l$ nonsingular points of $D$,
$D$ has at most $l$ components. Hence, $D$ must consist of exactly $l$
disconnected genus 0 components $\cup_{i=1}^l D_i$. Each $D_i$ is
fully ramified over $p_1$ with profile $m_i$. Therefore,
$D_i$ must be simply ramified over exactly $m_i-1$ finite points.
The proof of the property is complete.

As the number of factorizations of an $m$-cycle into
$m-1$ transpositions in $S_m$ is well-known to be $m^{m-2}$,
the solutions of (\ref{notri}) are now easily counted:
$$\binom{k}{m_1-1,\ldots,m_l-1} \prod_{i=1}^l {m_i^{m_i-2}} =
k!\prod_{i=1}^l \frac{m_i^{m_i-1}}{m_i!}.$$
\epf

\subsubsection{Localization}
The virtual localization formula for $\overline{M}_g(\proj^1,d)$
yields:
\begin{equation}
\label{ppwwe}
[\overline{M}_g(L_k)]^{vir} = \sum_{\Gamma\in G_g(\proj^1,d)}
\frac{1}{|\A|} \frac{ \tau_{\Gamma*}[\overline{M}_\Gamma] \cap \br^*[L_k]}
{e(N_\Gamma^{vir})}
\end{equation}
in $A_*^{\com^*}(\overline{M}_{g}(L_k))_{[\frac{1}{t}]}$ via the
Gysin map.

Let $\xi$ be the $\com^*$-equivariant lift of the point class $\xi_p$ of
$L_k$
corresponding to the $\com^*$-fixed point $r[p_0]+ k[p_1] \in L_k$.
The integral
\begin{equation}
\label{pjl}
\int_{[\overline{M}_g(L_k)]^{vir}} \br_k^*(\xi)
\end{equation}
is determined by the localization formula (\ref{ppwwe}).

However, the Hurwitz number $H_{g,\mu}$ is
not equal to (\ref{pjl}), but rather to 
the corresponding integral over $\overline{M}_g(\mu)$ by
Proposition \ref{pjjl}.
The central result is the identification of the contribution of
$\overline{M}_g(\mu)$ to the integral (\ref{pjl}).

Let
$\Gamma_\mu$ be the following distinguished graph.
$\Gamma_\mu$ has a unique genus $g$ vertex $v_0$ lying over $p_0$ 
which
is incident to exactly
$l$ edges of degrees $m_1,\ldots,m_l$. The edges connect $v_0$ to
$l$ 
unstable, unmarked vertices $v_1^1,\ldots,v_1^l$
of valence 1 and genus 0 lying over $p_1$. 

\begin{pr}
\label{sawwas}
$$mult(\mu) \int_{\overline{M}_g(\mu)} \br_k^*(\xi) =\frac{1}{|A_{\Gamma_\mu}|}
\int_{\overline{M}_{\Gamma_\mu}}
\frac{ \br^*[L_k] \cup \br_k^*(\xi)}
{e(N_{\Gamma_\mu}^{vir})}.$$
\end{pr}

\noindent
Proposition \ref{sawwas} is proven in Section \ref{huhh} below.

By definition,
$\overline{M}_{\Gamma_\mu}= \overline{M}_{g,l}$.
The order of the automorphism group is easily determined:
$$|A_{\Gamma_\mu}|= |\text{Aut}(\mu)| \prod_{i=1}^l m_i.$$
The vertex contributions of the Euler class of the
normal complex,
$$\frac{1}{e(N_{\Gamma_\mu}^{vir})} = \frac{1}{\prod_{v\in V} N(v)},$$
 were found in Section \ref{convertl}:
$$\frac{1}{N(v_0)}= \frac{t^g-t^{g-1}\lambda_1+t^{g-2}\lambda_2 
- t^{g-3}\lambda_3 + \ldots + (-1)^g \lambda_g}
{\prod_{i=1}^d (\frac{t}{m_i}-\psi_i)} t^{l-1-d} 
\prod_{i=1}^l\frac{m_i^{m_i}}{m_i!} ,$$
for the unique vertex over $p_0$ and
$$\frac{1}{N(v_1^i)}= -\frac{t^{1-m_i}}{m_i} \frac{m_i^{m_i}}{m_i!},$$
for the $i^{th}$ unstable vertex over $p_1$.

By the excess intersection formula, the class
$$\br^*(L_k)\cup \br^*_k(\xi)|_{\overline{M}_\Gamma} =(-1)^k r! k! t^{r+k}$$
 is the $\com^*$-equivariant Euler class of the
normal bundle of the point $r[p_0]+k[p_1]$ in $\text{Sym}^{r+k}(\proj^1)$.

After substitution of these identified factors,
Propositions \ref{pjjl} - \ref{sawwas} and Lemma \ref{mmull}
yield Theorem \ref{rrrr}.
The Hurwitz number $H_{g,\mu}$ equals
$$ 
\frac{(2g-2+|\mu|+l)!}{|\text{\em Aut}(\mu)|} 
\prod_{i=1}^l \frac{m_i^{m_i}}{m_i!}
\int_{\overline{M}_{g,l}}
\frac{1-\lambda_1+\lambda_2 -\lambda_3 + \ldots + (-1)^g \lambda_g}
{\prod_{i=1}^l (1-m_i \psi_i)},
$$
in the stable range $2g-2+\ell(\mu)>0$.

\subsubsection{Localization isomorphisms}
The following result proven in \cite{EdGra,Kr} will be used several times in
the proof of Proposition \ref{sawwas}.

\begin{lm}
\label{egkr}
Let $V$ be an algebraic variety (or Deligne-Mumford stack) equipped
with a $\com^*$-action. Let 
$$\iota: \cup_i V^f_i \rarr V$$ be the inclusion of the connected
components of the $\com^*$-fixed locus of $V$. Then $\iota_*$ is
an isomorphism after localization:
\begin{equation}
\iota_*:\bigoplus_i A^{\com^*}_*(V^f_i)_{[\frac{1}{t}]} \stackrel{\sim}{\rarr}
A^{\com^*}_*(V)_{[\frac{1}{t}]}.
\end{equation}
\end{lm}
\bpf 
We prove the result in case $V$ admits a nonsingular $\com^*$-equivariant
embedding $V\rarr Y$. The full result in proven in \cite{EdGra,Kr}.

The surjectivity of $\iota_*$ after localization
follows from the right exact sequence of equivariant
Chow groups of a closed inclusion:
$$A^{\com^*}_*(\cup_i V_i^f) \rarr A^{\com^*}_*(V) \rarr A^{\com^*}_*(U) \rarr 0.$$
Since $U$ admits a fixed point free $\com^*$-action, there is an
isomorphism
$$A^{\com^*}_*(U) \stackrel{\sim}{=} A_*(U/\com^*).$$
The right Chow group has finite grading (as $U/\com^*$ is
a finite dimensional algebraic variety (or Deligne-Mumford stack)) . 
Therefore, $A^{\com^*}_*(U)$ is $t$-torsion and vanishes after
localization. 

Injectivity is easily proven in 
case $V$ admits an equivariant nonsingular embedding $V\rarr Y$.
Let $$j:Y^f\rarr Y$$ denote the inclusion of the $\com^*$-fixed
locus. $Y^f$ is nonsingular (but possibly disconnected)
Let $N$ denote the normal bundle of $Y^f\subset Y$. 
Since $Y^f \cap V = \cup_i V_i^f$,
there is a Gysin map obtained by intersection with $Y_f$ in $Y$:
$$j^!:A^{\com^*}(V) \rarr \bigoplus_i A^{\com^*}(V^f_i).$$
The composition $ j^! \circ \iota_* $ is equal to multiplication by $e(N)$.
As 
$e(N)$ is invertible after localization, $\iota$
is injective after localization.
\epf

As the moduli space $\overline{M}_{g}(\proj^1,d)$ admits $\com^*$-equivariant
nonsingular embeddings, Lemma \ref{egkr} will only be used in the
restricted case considered in the proof.

\subsubsection{Proof of Proposition \ref{sawwas}}
\label{huhh}
Let $X_0= \overline{M}_g(\mu)$.   
$X_0$ is a union of irreducible components of $\overline{M}_g(\mu)$ of
pure dimension $r$.
Let
$\cup_{j\geq 0} X_j = \overline{M}_g(L_k)$ where 
$\{X_j\}_{j\geq1}$ are the remaining irreducible
components of $\overline{M}_g(\mu)$.
The virtual class admits a (non-canonical)
decomposition:
$$[\overline{M}_g(L_k)]^{vir} = \iota_*  
\sum_{j\geq 0} R_j \ \ \in A^{\com^*}_{r}(\overline{M}_g(L_k)),$$
where $R_j \in A_{r}^{\com^*}(X_j)$. 
Since
$X_0$ is of multiplicity $mult(\mu)$
in $\br^{-1}(L_k) \cap M_g(\proj^1,d)$, 
$$R_0= mult(\mu) \ [X_0].$$ 
The classes $\{R_j\}_{j\geq 1}$ are quite difficult to describe.

The $\com^*$-fixed loci of $\overline{M}_g(L_k)$ correspond to the set of graphs
$$G_g(\mu) \subset G_g(\proj^1,d)$$
satisfying $\br(\overline{M}_\Gamma) \in L_k$.
Let $Q_{\Gamma}= \overline{M}_\Gamma/ A_\Gamma$
denote the $\com^*$-fixed locus corresponding to $\Gamma$ (as in Section \ref{ccfixx}).
The localization formula (\ref{ppwwe}) may be written as:
$$[\overline{M}_g(L_k)]^{vir} = \iota_* \sum_{\Gamma\in G_g(\proj^1,d)} C_\Gamma \ \ \in
A^{\com^*}_*(\overline{M}_g(L_k))_{[\frac{1}{t}]},$$
where $C_\Gamma \in A^{\com^*}_*(Q_\Gamma)_{[\frac{1}{t}]}$.

To prove Proposition \ref{sawwas}, we must show:
\begin{equation}
\label{bbaaq}
\int_{\overline{M}_g(L_k)} \iota_*R_0 \cap br_k^*(\xi) = \int_{\overline{M}_{g}(L_k)} 
\iota_* C_{\Gamma_\mu} 
\cap \br_k^*(\xi).
\end{equation}

For each $j\geq 0$, we may use the localization isomorphism of
Lemma \ref{egkr} to uniquely 
determine classes $R_{j,\Gamma}\in A^{\com^*}_*(Q_\Gamma)
_{[\frac{1}{t}]}$ satisfying:
$$\iota_*R_j = \iota_* \sum_{\Gamma\in G_g(\mu)} R_{j,\Gamma}.$$
The localization isomorphism implies:
$$\sum_{j\geq 0} R_{j,\Gamma} = C_\Gamma$$
for all $\Gamma \in G_g(\mu)$.

We may rewrite the desired equation (\ref{bbaaq}) in the following form:
\begin{equation}
\label{bbaaqq}
\int_{\overline{M}_g(L_k)} \iota_* \sum_{\Gamma\in G_{g}(\mu)} 
R_{0,\Gamma} \cap br_k^*(\xi) = \int_{\overline{M}_{g}(L_k)} 
\iota_* \sum_{j\geq 0} R_{j,\Gamma_\mu} 
\cap \br_k^*(\xi).
\end{equation}
It will therefore suffice to prove the following vanishing
results:
\begin{enumerate}
\item[(i)]$\int_{\overline{M}_g(L_k)}
\iota_* R_{0,\Gamma} \cap br_k^*(\xi)= 0$ 
for $\Gamma\neq \Gamma_\mu$,\\
\item[(ii)] $\int_{\overline{M}_g(L_k)}
\iota_* R_{j,\Gamma_\mu} 
\cap \br_k^*(\xi)=0$
for $j\neq 0$.
\end{enumerate}
A study of 
the component geometry of $\overline{M}(L_k)$ will
be required to prove (i) and (ii).

Let $\Gamma\in G_g(\mu)$ and let $\br(\overline{M}_\Gamma)=
a_{\Gamma,0}[p_0]+a_{\Gamma,1}[p_1]$. The inequality
$$a_{\Gamma,1} \geq k$$
holds since
$\br(\overline{M}_\Gamma)\in L_k$.
If $a_{\Gamma,1}>k$, then $[Q_\Gamma] \cap \br_k^*(\xi)=0$ as
$\xi$ is the class corresponding to the $\com^*$-fixed point $r[p_0]+k[p_1]$.
Since $R_{j,\Gamma} \in A^{\com^*}_r(Q_\Gamma)_{[\frac{1}{t}]}$, we find
the trivial vanishing:
\begin{enumerate}
\item[(a)]
$R_{j,\Gamma}\cap \br^*_k(\xi)=0$
if $a_{\Gamma,1}>k$.
\end{enumerate}

As $X_j$ is a $\com^*$-equivariant locus, we may apply the
localization isomorphism to decompose $R_j$ {\em on the
$\com^*$-fixed locus of $X_j$.} By uniqueness, we conclude
another trivial vanishing:
\begin{enumerate}
\item[(b)]
$R_{j,\Gamma} = 0$ if $X_j \cap Q_\Gamma = \emptyset$.
\end{enumerate}

\vspace{10pt}
\noindent{\em Proof of (i).}
The fixed locus
$Q_{\Gamma_\mu}$ meets $X_0=\overline{M}_g(\mu)$. In fact
the limit 
$$\text{Lim}_{t\rarr 0} \ \ t\cdot [\pi]$$ of
every element of $M_{g}(\mu)$ lies in $Q_{\Gamma_\mu}$ (see \cite{KiP}). 

Define a stable map $\pi$ to have {nonsingular profile $\mu$ over $p_1$}
if $\pi^{-1}(\mu)$ is a divisor of shape $\mu$ lying in the
{\em nonsingular locus} of the domain.
A limit $[\pi]$ of elements in
$M_g(\mu)$ must either have nonsingular profile $\mu$ {\em or}
degenerate over $p_1$. Any degeneration is easily seen
to {\em increase} the branching order of $\pi$ over $p_1$.
As $Q_{\Gamma_\mu}$ is the unique fixed locus with nonsingular profile
$\mu$ over $p_1$, 
$Q_{\Gamma_\mu}$ is the unique fixed locus meeting
$\overline{M}_g(\mu)$  with branching order
exactly $k$ over $p_1$.  
Vanishings (a) and (b) then imply (i).

\vspace{10pt}
\noindent{\em Proof of (ii).}
Let $j\neq 0$. By vanishing (b), we may assume $X_j\cap Q_{\Gamma_\mu} \neq \emptyset$.
Since $X_j\subset \overline{M}_g(L_k)$, every element
$[\pi]\in X_j$ corresponds to a map with branching order at least $k$ over $p_1$.
Since $X_j\cap Q_{\Gamma_\mu} \neq \emptyset$, the general
map $[\pi]\in X_j$ must have nonsingular
profile $\mu$ over $p_1$.

As in the proof of (i), $Q_{\Gamma_\mu}$ must be
the unique fixed locus meeting
$\overline{M}_g(\mu)$  with branching order
exactly $k$ over $p_1$.

Maps $\pi$ with no contracted components and nonsingular profile
$\mu$ over $p_1$ are easily shown to be limits of $M_g(\mu)$.
As $X_j\neq X_0$,
the general map $[\pi]\in X_j$ must contain a domain component collapsed
away from $p_1$. By the definition of the branch morphism, 
$\br(\pi)$ then lies in the singular
sublocus $L_k^{sing} \subset L_k$:
$$L_k^{sing}= \{\ D+ k[p_1] 
\ |\ D=2[x_1]+[x_2]+\cdots +[x_{r-1}] \in \text{Sym}^{r}(\proj^1) \}.$$
As $L_k^{sing}$ is a proper subvariety, the following integral vanishes:
\begin{equation}
\label{ddfccs} 
\int_{\overline{M}_g(L_k)} \iota_* R_j \cap \br_k^*(\xi) =0.
\end{equation}

By the localization isomorphism, the integral (\ref{ddfccs}) may be rewritten as
$$\int_{\overline{M}_{g}(L_k)} \iota_*\sum_{a_{\Gamma,1}=k} 
R_{j,\Gamma} \cap \br_k^*(\xi) 
+ \iota_* \sum_{a_{\Gamma,1} > k} R_{j,\Gamma} \cap
\br^*_k(\xi)=0.$$
The second sum vanishes completely by (a).
As $\Gamma={\Gamma_\mu}$ is
the unique graph
satisfying $a_{\Gamma,1}=k$ and $X_j\cap Q_{\Gamma} \neq \emptyset$,
all other term in the first sum vanish by (b). We conclude:
$$\int_{\overline{M}_{g}(L_k)} \iota_* R_{j,\Gamma_\mu}=0$$
for $j\neq 0$.

\vspace{10pt}
\noindent
The proof of 
Proposition \ref{sawwas} is complete.
\epf

%
%
%


\part{Asymptotics of Hurwitz numbers}

\section{Random trees}
\label{st2}
\subsection{Overview}
The analysis of the $N\to\infty$ asymptotics of the
Hurwitz numbers 
$H_{g, N\mu}$ via the asymptotic enumeration of  
branching graphs 
will require a study of trees.
Trees naturally arise via edge
terms in the homotopy classification of branching graphs.
The enumerative and probabilistic
results for trees which will be required 
are discussed here.
The asymptotic analysis of $H_{g,N\mu}$ is undertaken
in Section \ref{fin2}.

Section \ref{revprob} contains a minimal discussion of
probabilistic terminology. Section \ref{tr1} is a review  of the basic
enumeration formulas for trees. 
The required properties of random edge trees are discussed
in Sections \ref{hrt}-\ref{semiper}.

The literature on trees and random trees is very large. 
An excellent place to start is Chapter 5 of \cite{Sta}.
An introduction from a more probabilistic perspective can be
found in \cite{Pit}. Many asymptotic properties of random
trees find a unified treatment in the theory of continuous
random trees due to Aldous \cite{Al1, Al2}.
Fortunately, all the properties of 
random trees that we shall need are quite basic. Instead
of locating them in the literature, we will 
prove these properties from first principles.

The trees that we will consider
naturally come with a choice of two distinguished vertices
(a root and a top). Random trees are more often studied with
one special vertex (rooted trees) or with no special
vertices (plain trees). The properties we will require are
simpler in the presence of a root and a top. The
analogous results for
rooted trees are less elementary both to 
state and to prove. 

 
\subsection{Review of probabilistic terminology}\label{revprob}

A {\em probability space} is a triple $(\Omega,\mathfrak{B},P)$,
where $\Omega$ is any nonempty set, $\mathfrak{B}$ is a $\sigma$-algebra
of subsets of $\Omega$ called the algebra of events, and 
$$
P: \mathfrak{B} \to \R_{\ge 0}
$$
is a measure such that $P(\Omega)=1$. We will primarily be concerned
with finite
sets $\Omega$:  $\mathfrak{B}$ will then include
 all subsets of $\Omega$,
and 
$P$  will typically be the uniform probability measure.
 When the probability measure 
is understood, it will be denoted by the symbol $\Prob$. 

Any measurable function
$$
X : \Omega \to \R^k
$$
is called a  vector-valued \emph{random variable}. The push-forward measure $X_*P$
on $\R^k$ is called the distribution of $X$. The integral
$$ 
\la X \ra= \int_\Omega X(\omega) \, P(d\omega) = \int_{\R^k} x \, X_*P(dx) 
$$
is called the {\it expectation} of $X$.   Two random 
variables $X$ and $Y$ are said to be {\it independent} if the distribution
of their direct sum $X\oplus Y$, also known as the
joint distribution of $X$ and $Y$, is a product-measure. 

A sequence $\{\mm_n\}$ of measures  on $\R^k$ is said to {\it converge
weakly} to a measure $\mm$ if
$$
\int_{\R^k} f(x) \, \mm_n(dx) \to \int_{\R^k} f(x) \, \mm(dx)
$$
for any bounded continuous function $f$.  A sequence of random variables
$X_n$ on a sequence of probability spaces $(\Omega_n,\mathfrak{B}_n,P_n)$
is  said to {\it converge in distribution} to a random variable $X_\infty$
 if the measures
$\mm_n=(X_n)_*P_n$ converge weakly to the distribution of $X_\infty$, or,
equivalently, if
$$
\la f(X_n) \ra \to \la f(X_\infty) \ra  \,, \quad n\to\infty\,, 
$$
for any bounded continuous function $f$.

In particular, $X_n$
converge in distribution to the variable identically equal to $0$ if
$$
\Prob(X_n \in U) \to 1
$$
for any neighborhood of $U$ of $0$.  A basic continuity 
property of convergence in distribution, which we will use often,
is 
the following standard result (see, for example, \cite{Bil}).

\begin{lm}\label{limdist} 
Let $X_n$ and $Y_n$, $n=1,2,\dots,\infty$, be vector-valued 
random variables on a sequence 
$(\Omega_n,\mathfrak{B}_n,P_n)$   of probability spaces. If, as $n\to
\infty$,  we have
$$
X_n \to X_\infty, \quad X_n - Y_n \to  0\,,
$$
in distribution, then also 
$$
Y_n \to X_\infty
$$
in distribution. 
\end{lm}

\begin{proof}
Let $\|X\|$ denote a vector norm of $X$ and consider the function
$$
g_{AB} : \R_{\ge 0} \to [0,1]
$$
such that $g_{AB}(x)=0$ for $x>B$, $g_{AB}(x)=1$ for
$x<A$ and $g_{AB}$ linearly interpolates between $1$ and $0$
on $[A,B]$. Clearly, for any $X$, 
$$
\Prob\{\|X\| \le A\}\le \la g_{AB}(\|X\|) \ra \le \Prob\{\|X\| \le B\} \,.
$$
For any $\epsilon>0$, we can find $A$ such that
$$
\Prob\{\|X_\infty\|\le A\} > 1- \epsilon  \,.
$$
By hypothesis, for any $A$ and $B$ we have 
$$
 \la g_{AB}(\|X_n\|) \ra \to \la g_{AB}(\|X_\infty\|) \ra\,, \quad  n\to\infty \,.
$$
Therefore, for any $B>A$,
$$
\Prob\{\|X_n\| \le B \} > 1- 2\epsilon 
$$
for all sufficiently large $n$. For any $C>B$, we have
$$
\Prob\{\|Y_n - X_n\| > C-B\} \to 0 \,,
$$
and, therefore, for all sufficiently large $n$ we have
$$
\Prob\{X_n\in K\} > 1- 2\epsilon \,, 
\quad \Prob\{Y_n \in K\} > 1- 3\epsilon \,, 
$$
where $K$ denotes the compact set 
$$
K=\{X, \|X\| \le C\} \,.
$$

Since $f$ is continuous and $K$ is compact, $f$ is uniformly continuous on $K$
and hence there exists $\delta>0$ such that
$$
|f(X)-f(Y)|<\epsilon \,, 
$$
whenever $X,Y\in K$ and $\|X-Y\| < \delta$.  We can choose 
$n$ large enough so that
$$
\Prob \{ \|X_n-Y_n\| < \delta\} > 1- \epsilon \,.
$$
Collecting all estimates, we obtain
$$
\left| \la f(X_n)  - f(Y_n) \ra \right| \le 
\epsilon (1+ 12 \max |f|) \,.
$$
for all sufficiently large $n$. Since $\epsilon$ is
arbitrary, the Lemma follows.
\end{proof}

The above Lemma will often be used in  the following situation.
Suppose there exists a ``good'' subset 
$$
\Omega'_n \subset \Omega_n
$$
such that on this good subset we have
$$
\sup_{\omega\in \Omega'_n} |X_n(\omega)-Y_n(\omega)| \to
0\,, \quad n\to\infty \,,
$$
and also such that, asymptotically, most $\omega$ are good, that is,
$$
\Prob  \Omega'_n \to 1 \,, \quad n\to\infty \,.
$$
In this case, Lemma \ref{limdist} implies that if $X_n$ has a limit in
distribution then $Y_n$ converges to the same limit.

\subsection{Enumeration of trees}
\label{tr1}
\subsubsection{Definitions}
A {\em tree} $(V,E)$ is a connected graph with no circuits.
Let $\T(n)$ denote the set of trees with $n$ vertices.
We will consider trees $T$ with additional structures: 
vertex and edge labels, and distinguished vertices. 

Labelings of vertices and edges are bijections
$$
\phi_V: V \to \{1,2,3,\dots,|V|\}\, ,
$$
$$
\phi_E: E \to \{1,2,3,\dots,|E|\}\, .
$$
We will denote the set of vertex marked trees with
$n$ vertices by
$\V(n)$. Let $\E(n)$ denote the set of edge marked trees with
$n$ vertices.

One of the vertices of a tree $T$ may be designated as
a distinguished vertex, called the \emph{root} of $T$.
The tree $T$ is this case is called a \emph{rooted tree}.
Let $\T^1(n)$ denote the set of rooted trees with $n$ vertices.
Similarly, let 
$\V^1(n)$ and $\E^1(n)$ denote the sets $n$ vertex rooted trees with
marked vertices and 
marked edges respectively.

In addition to the root vertex, one may choose a \emph{top} vertex of $T$. 
The top vertex may or may not be allowed to
coincide with the root vertex. Let $\V^{11}$
denote vertex marked trees with distinct root
and top vertices. $\V^{2}$ will denote the larger set in which
root and top are allowed to coincide. Let
$\E^{11}$ and $\E^2$ denote the corresponding sets for edge
marked trees.

Of these flavors of trees, two will be particularly
important for us and deserve special names. Edge
marked trees will be also called  \emph{branching
trees}, and the $\E^{11}$-trees will be called  
\emph{edge trees}. Edge trees naturally arise in
the study of the edge contributions 
in the asymptotic analysis of branching graphs (see
Section \ref{fin2}), whence the name. The term branching
tree is justified by the following:

\begin{lm}\label{brtrpl}
A branching tree with $n$ vertices is isomorphic to the data of a
branching graph on the sphere
$\Sigma_0$ with perimeter $(n)$, where $(n)$ denotes the length 1 partition
of $n$.
\end{lm}
 
\bpf
First, we make a general remark. By definition, the
edges of a branching graph are labeled by roots of unity,
whereas the edges of an edge marked tree $T\in\E(n)$
are labeled by
$1,2,3,\dots,n-1$. We will 
identify the two kinds of labeling using the bijection 
$$
\{1,\dots,n-1\} \owns k \mapsto e^{2\pi i k/(n-1)} \in U_{n-1} \,.
$$
A branching tree $T$ can be canonically (up to homeomorphism) embedded
in an oriented sphere $\Sigma_0$.
The embedding is uniquely determined by the following condition:
the cyclic order induced on the
edges 
incident to each vertex by the orientation of $\Sigma_0$ must
agree with the cyclic order of the markings of the edges. The tree   
$T\subset \Sigma_0$ then defines a branching graph on $\Sigma_0$ (see
Section \ref{hur}).

Conversely, every branching graph on $\Sigma_0$ must be a
tree (as the complement determines 1 cell). The edge
markings then determine a branching tree structure.
\epf

\subsubsection{Automorphisms and counting}
Trees $T\in T(n)$ may have non-trivial automorphism
groups. However, labelled trees in the sets
$\V(n)$ and $\E(n)$
admit no non-trivial
automorphisms preserving their markings, the only exception being the unique
element of $\E(2)$. 

We will exclusively count labelled trees (with
distinguished vertices). Therefore,
by the number of trees, we will 
mean the actual number (except in the $\E(2)$ case
where the number is set, by definition, to $1/2$ in order to account
for the order 2 automorphism group). 
Similarly, when considering random labelled trees,  
we will always take the uniform
probability measure on the corresponding set.

\subsubsection{Cayley's formula and its consequences}

We recall the following fundamental result about trees:

\begin{pr}[Cayley] We have 
  \begin{equation}
    \label{Cayley}
\sum_{T\in \V(n)} 
\prod_{i=1}^n z_i^{\val(i)}  = z_1 \cdots z_n (z_1+\dots+z_n)^{n-2}\,
\end{equation}
where the summation is over all trees $T$ with vertex set
$\{1,\dots,n\}$ and $\val(i)$ denotes the valence of the 
vertex $i$ in the tree $T$. 
\end{pr} 

See, for example \cite{Sta}, Theorem 5.3.4,
 for a proof of this formula. The formula
\eqref{Cayley} has a large number of  corollaries.

\begin{cor}\label{treeenum}
We have
\begin{alignat*}{4}
&|\V(n)|&&=n^{n-2}\,, \qquad & &|\E(n)|&&=n^{n-3} \,,  \\
&|\V^1(n)|&&=n^{n-1}\,, \qquad & &|\E^1(n)|&&=n^{n-2} \,,  \\
&|\V^{11}(n)|&&=(n-1)\,n^{n-1}\,, \qquad & &|\E^{11}(n)|&&=(n-1)\,n^{n-2} \,,  \\
&|\V^{2}(n)|&&=n^{n}\,, \qquad & &|\E^{2}(n)|&&=n^{n-1} \,,
\end{alignat*}
\end{cor}

Recall that $|\E(2)|=1/2$, by our convention, reflects the
order 2 automorphism group of the unique element of $\E(2)$.

\begin{proof}
The enumeration of $\V(n)$ is obtained by setting $z_i=1$, $i=1\dots n$, in 
Cayley's formula \eqref{Cayley}. 

Given a vertex marked tree $T$ with $n$ vertices, 
one can mark its edges in $(n-1)!$ ways. The vertex marking
can then be removed by  dividing by $n!$ which gives $|\E(n)|=n^{n-3}$.
The remaining formulas are obvious.
\end{proof}

\begin{cor}\label{c_val} The number of trees in $\V(n)$
such that the valence $\val(1)$ of the vertex $1$ is
$k+1$ equals $(n-1)^{n-k-2}\, \binom{n-2}{k}$.
\end{cor}

This is obtained by setting $z_i=1$, $i=2\dots n$, in \eqref{Cayley}
and extracting the coefficient of $z_1^{k+1}$. 

Consider the probability that in a uniformly random vertex
marked tree $T\in \V(n)$ the valence $\val(1)$ of the vertex marked by $1$
equals $k+1$. We have 
$$
\Prob\big\{\val(1)=k+1\big\} = 
\frac{(n-1)^{n-k-2}\, \binom{n-2}{k}}{n^{n-2}} \to \frac{e^{-1}}{k!} \,, \quad
n\to\infty \,.
$$
In other words, the valence distribution of a given vertex in
a large random tree converges in distribution to
one plus a Poisson random
variable with mean 1. 

This observation has an immediate generalization
for the joint distribution of valences of several vertices. Given a
vertex $v\in T$, let us call the number $\val(v)-1$ the \emph{excess valence}
of the vertex $v$. 

\begin{cor}\label{poisslimit}
As $n\to\infty$, the excess valences of vertices 
of a random tree $T\in\V(n)$ converge in distribution
to independent Poisson random variables
with mean $1$.
\end{cor}

Recall that a {\em forest} is graph which is a disjoint union
of trees. A forest is {\em rooted} if a distinguished vertex,
called root, is specified in each connected component. 
 
\begin{cor}\label{numfor} 
The number of rooted forests with  with vertex set
$\{1,\dots,n\}$ and $k$ connected components is equal to
$k\, \binom{n}{k} \, n^{n-k-1}$.
\end{cor}

\begin{proof}
There exits a simple bijection between rooted forests that we want to enumerate
and trees with vertex set $\{0,1,\dots,n\}$ such that the 
vertex $0$ is $k$-valent. We just add new edges which
join $0$ to the roots of the forests. Now we apply Corollary \ref{c_val}. 
\end{proof}

\subsubsection{Factorization into transpositions and trees}
By Definition \ref{Hur2}, the Hurwitz number
$H_{0,(n)}$ equals
the automorphism weighted count of branching graphs on the
sphere $\Sigma_0$ with one cell of perimeter $n$. 
By Corollary \ref{treeenum},
$H_{0,(n)} = n^{n-3}$.

By Definition \ref{Hur3},
$H_{0,(n)}$ is also equal to $(1/n!)$ times
the number of $(n-1)$-tuples of transpositions in $S_n$
with product in the conjugacy class of an $n$-cycle. 
Equivalently,
$n H_{0,(n)}$ equals the number of solution to the equation:
$$\gamma_1\ldots\gamma_{n-1}= (123\ldots n) \in S_n,$$
for $2$-cycles $\gamma_i\in S_n$.

We therefore obtain the following classical result (used in the
proof of Lemma \ref{solomon} in Section \ref{kkrr}):
\begin{cor}
\label{dqqqp}
The number of factorization of an $n$ cycle into $n-1$ transpositions
in $S_n$ is $n^{n-2}$.
\end{cor}

\noindent
Corollary \ref{dqqqp}
is a 
particular case of a formula due to Hurwitz \cite{Hu,Str} and
was also discovered by D\'enes \cite{De}.

\subsection{Trunk of a random edge tree}
\label{hrt}

Given $T\in\E^{11}(n)$, denote by $\tk T$ the {\em trunk} of $T$, that is,
the shortest path from then root to the top in $T$. Let $|\tk T|$ denote 
the number of vertices in the trunk of $T$.  We are interested 
in the distribution of 
this quantity with respect to the uniform probability
measure on $\E^{11}(n)$ as $n\to\infty$.

Recall that an exponential random variable $\xi$ with mean $1$ is,
by definition, the variable with distribution density
$e^{-x} \, dx$ on $[0,+\infty)$. The random variable $\sqrt{2\xi}$, which
 has the density
$x  \, e^{-x^2/2} \, dx$ on the half-line $(0,\infty)$, is called a
\emph{Rayleigh} random variable. 

\begin{pr}\label{ptrh}
As $n\to\infty$, the random variable $\dfrac1{\sqrt{n}}\, |\tk T|$,
where $T$ is a random edge tree with $n$ vertices, 
converges in distribution to a Rayleigh random variable. 
\end{pr}

\begin{proof}
The same distribution of trunk heights is obtained
if, instead of edge trees, we consider random elements of $\V^{11}(n)$. 
The notion of trunk and its height have an
obvious analog for such trees.

Given a tree $T\in \V^{11}(n)$ with $|\tk T|=k$ and
 $n$ vertices, we can associate
to it a forest with $k$ components by deleting the trunk path $\tk T$
from $T$.
This forest
comes with an additional structure, namely, an ordering
on the components of the forest. Since there are $k!$
possible orderings, we conclude using Corollary \ref{numfor}
that the probability to
have $|\tk T|=k$  equals
\begin{equation}\label{prtrh}
\left. k \, k!\binom{n}{k}\,n^{n-k-1} \right/ (n-1)\, n^{n-1} = \frac{k}{n-1}\,
\prod_{i=1}^{k-1} \left(1-\frac i{n}\right) \,. 
\end{equation}
If $k=x\sqrt n$ then, as $n\to\infty$, we have
$$
\ln \, \prod_{i=1}^{k-1} \left(1-\frac in\right) 
\sim - \frac1n \sum_{i=1}^{k-1}
i \to -x^2/2\,,
$$
hence the probability \eqref{prtrh} is asymptotic to 
$
\displaystyle \frac{1}{\sqrt n} \, x e^{-x^2/2}
$, 
which completes the proof. 
\end{proof}

\begin{cor}\label{tk0}
For any $\epsilon>0$, we have
$$
\Prob \left\{ |\tk T|> n^{1/2+\epsilon} \right\} \to 0 \,, \quad
n \to \infty \,,
$$
with respect to the uniform probability measure on $\E^{11}(n)$. 
\end{cor}

The trunk of a tree $T$
appears in the literature under various names.
See, for example, \cite{Jo, Lab, MM0}.
In particular, our trunk is called the spine of $T$ in \cite{AlPit}.

\subsection{Size of the root component of a random tree}
\label{sofrc}

\begin{figure}[!hbt]
\centering
\scalebox{.6}{\includegraphics{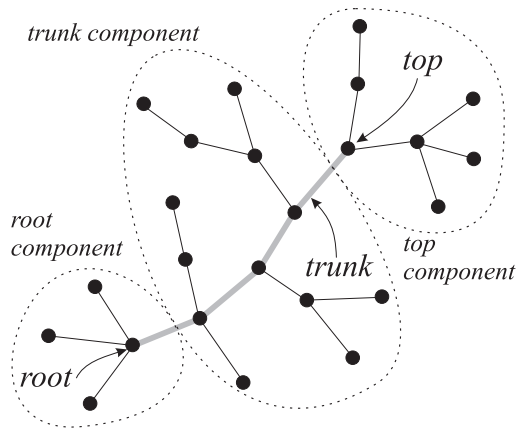}}
\caption{The components of tree $T\in \E^{11}(r)$}
\label{fig7}
\end{figure}

Given  $T\in \E^{11}(r)$, consider the edges incident to the root
vertex. One of these edges belongs to the trunk $\tr T$,
we will call it the \emph{trunk edge}. One of the two components of
$T$ that the trunk edge separates contains the root vertex,
we call this component  the
{\em root component} of $T$. We define the \emph{top
component} of $T$ similarly and call the complement of
the root and top components of $T$ the \emph{trunk component}
of $T$. These notions are illustrated in Figure \ref{fig7}

\begin{pr}\label{rootcomp}
As $n\to\infty$, the probability that the root component of a
random edge tree $T\in \E^{11}(n)$ contains
$k$ vertices has limit $\displaystyle \frac{k^{k-1}}{k!} \, e^{-k}$.
\end{pr} 

\begin{proof}
As in proof of Proposition \ref{ptrh}, we can replace
random edge trees by random elements of $\V^{11}$. 
We can construct elements of $\V^{11}(n)$ with given
root component of size $k$ as follows: partition the
$n$ vertices into sets of order $k$ and $n-k$, take an element
of $\V^{1}(k)$ and an element of $\V^{2}(n-k)$, join
their roots by an edge, and choose the root of the first
tree to be the root of the union. 

It follows that the probability that
the root component has size $k$ equals
$$
\binom{n}{k}\, \frac{k^{k-1} \, (n-k)^{n-k}}{(n-1)\, n^{n-1}} 
\to  \frac{k^{k-1}}{k!} \, e^{-k} \,, \quad k\to\infty \,,
$$   
where the asymptotics follow immediately from the
Stirling formula \eqref{Stirl}. 
\end{proof}

The root component of an element of $\V^{11}(n)$ determines
a rooted tree in $\T^1$ after forgetting the vertex labels.
The argument for Lemma \ref{rootcomp} proves more precise statements.

\begin{cor}\label{roottopcomp1}
The probability that $T\in \T^{1}(k)$
corresponds to
the root component of a random tree $T\in\V^{11}(n)$ is asymptotic
to
$$\frac{e^{-k}}{|\Aut(T)|}$$ 
as $n\to\infty$. 
\end{cor}

\begin{cor}\label{roottopcomp2}
The top component of random tree
$\V^{11}(n)$ has the same distribution as the root component.
Moreover, in the $n\to\infty$ limit,
the root and top component distributions are independent. 
\end{cor}

The asymptotic probabilities of Proposition \ref{rootcomp}
determine a probability measure.

\begin{lm}\label{probdist}
The measure $\Prob(k)=\displaystyle \frac{k^{k-1}}{k!} \, e^{-k}$
is a probability measure on natural numbers. 
\end{lm}
\begin{proof}
This can be seen, for example, from the equation
$$
w(z)=z \, e^{w(z)}\,,
$$
satisfied by the function 
$$
w(z)=\sum_{k=1}^\infty \frac{k^{k-1}}{k!} \, z^k \,,
$$
which is the generating function for $|\V^{1}(k)|$ and
is, essentially, the same as the  Lambert W-function. The equation
implies that
$$
1 = w(1/e)=\sum_{k=1}^\infty \frac{k^{k-1}}{k!} \, e^{-k} \,.
$$ 
\end{proof}

In fact, the measure
$
\Prob(k)=\displaystyle \frac{k^{k-1}}{k!} \, e^{-k}
$
is the {\em Borel distribution} \cite{Bo} and is well known to appear
in the context of branching processes and random trees. See, for
example, Section 7 in \cite{Pit}.

Informally, Proposition \ref{rootcomp} and Lemma \ref{probdist} imply 
 that the size of the root component of a typical tree stays finite as the
size of the tree goes to infinity. A more formal statement
is the following:  

\begin{cor}\label{rootfin}
For any $\epsilon>0$ there exists $M$ such that for all $n$ the
probability that  a random
tree $T\in\E^{11}(n)$ has 
the root component with more than $M$ vertices is less
than $\epsilon$.  
\end{cor}

Similarly, we have: 

\begin{cor}
For any sequence $\{c_n\}$ such that $c_n\to\infty$,
 the probability that a random
tree $T\in\E^{11}(n)$ has 
the trunk component of size $\ge n-c_n$ goes 
to $1$ as $n\to\infty$. In other words, all but
finitely many vertices
of a typical large edge tree $T$ lie in the trunk 
component. 
\end{cor}

\subsection{Semiperimeters}\label{semiper}

\subsubsection{Definitions}

Let $T\in\E^{11}(n)$ be an edge tree. Make $T$ planar as
in Lemma \ref{brtrpl} and let $\lambda$ be a path with follows
the perimeter of $T$ once clockwise. Formally,
$\lambda$ is a function 
$$
\Z\owns k \mapsto \lambda_k \in E
$$
periodic with period $2n-2$, 
which lists the edges in the order of their appearance
along the boundary of $\Sigma_0 \setminus T$. 

Let $\phi:E\to\{1,\dots,n-1\}$ be the marking of the 
edges of $T$ which is, by definition, a part of the structure of an edge tree.
Define the angle between two edges $e,e'\in E$ by
\begin{equation}
  \label{angledef}
 \measuredangle(e,e')=\frac{2\pi(\phi(e')-\phi(e))}{n-1} \bmod 2\pi\,,
\quad  \measuredangle(e,e')\in
(0,2\pi] \,.  
\end{equation}
Consider the \emph{perimeter} of $\lambda$ which, by 
definition, equals 
$$
\per(\lambda)= \frac1{2\pi} \sum_{k=1}^{2n-1} 
\measuredangle(\lambda_k,\lambda_{k+1}) 
\,.
$$
Of course, $\per(\lambda)=n$ 
because every vertex contributes $1$ to the above  sum. 

\begin{figure}[!hbt]
\centering
\scalebox{.6}{\includegraphics{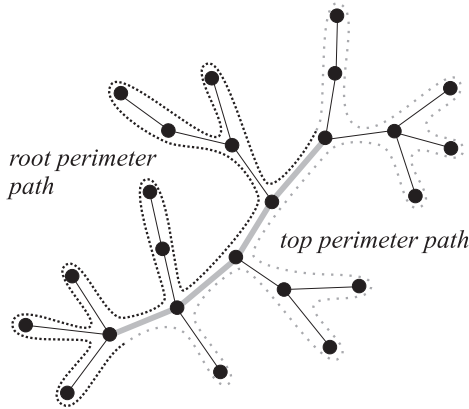}}
\caption{The paths  $\lambda_R$ and $\lambda_T$ for the tree from 
Figure \ref{fig7}}
\label{fig8}
\end{figure}

We now want to split the
path $\lambda$, and its perimeter, into two parts: the 
\emph{root perimeter path}
$\lambda_R$ and the \emph{top perimeter path} $\lambda_T$. 
We proceed as follows. Let $e_r,e_t\in E$ denote the trunk
edges at the root and the top of $T$, respectively. 
As we follow the path $\lambda$, these edges appear in cycles
of the form 
$$
(e_r,\underbrace{\dots,e_r,\dots,e_t}_{\lambda_R},\dots,e_t,\dots)\,,
$$
where the dots stand for other edges of $T$. 
The root part $\lambda_R$ starts after the first
appearance of $e_r$ and ends with the first appearance
of $e_t$ as shown above. 
Similarly, we define the top part $\lambda_T$. We also define
the two perimeters, $P_R$ and $P_T$ as  the perimeters 
of two paths $\lambda_R$ and $\lambda_T$, respectively, and call
them the \emph{semiperimeters} of $T$.  
Since the paths $\lambda_R$ and $\lambda_T$ are not closed,
these semiperimeters may be fractional. The definition
of $\lambda_R$ and $\lambda_T$ is illustrated in Figure \ref{fig8}

We also define the canonical marking 
$$
\psi: V\to \{1,\dots,n\}
$$
of the vertices $V$ by the order of their appearance in
the concatenated path $\lambda_R+\lambda_T$.

\subsubsection{Perimeter estimates}
Let us denote the root and top vertices by $v_r$ and $v_t$,
respectively. 

A basic consequence of the
definitions is:
\begin{lm}
\label{jjjj} For $T\in \E^{11}(n)$,
$| P_R + P_T  - n | \leq 2\,.$
\end{lm}
\bpf
The difference between $P_R+P_T$ and the vertex number $n$ occurs from losses at $v_r$ and $v_t$.
\epf

\begin{lm}
\label{kksk}
$|P_R- \psi(v_t)| \le  |\tk T|.$
\end{lm}
\bpf 
As we follow $\lambda_R$, every vertex on the trunk contributes
$1$ to $\psi(v_t)$ and between $0$ and $1$ to $P_R$. Every other vertex
 contributes $1$ to both $\psi(v_t)$ and $P_R$. 
\epf

\subsubsection{Semiperimeter distribution}
\label{semdis}
\begin{pr}
\label{xmml}
As $n\to\infty$, 
the normalized semiperimeter $\dfrac{P_R}{n}$ converges
in distribution to the uniform distribution 
 on $[0,1]$. 
\end{pr}

\begin{proof}
Since, by Lemma \ref{kksk},
$$
\left|\frac{P_R}{n}- \frac{\psi(v_t)}{n}\right| \le  \frac{|\tk T|}n
$$ 
and the right-hand side converges to $0$ in distribution by Corollary \ref{tk0},
it suffices to prove that $\dfrac{\psi(v_t)}{n}$ converges to
the $[0,1]$-uniform random variable. 

Consider the subset of $\E^{11}(n)$ formed by trees with the root
component of fixed cardinality $k\in \{1,2,\dots\}$. Clearly, on this
subset, $\psi(v_t)$ is uniformly distributed on the interval 
$\{k+1,\dots,n\}$
and, hence, on this subset, $\dfrac{\psi(v_t)}{n}$ converges, in
distribution,  to
the $[0,1]$-uniform random variable. Now Corollary 
\ref{rootfin} concludes the proof. 
\end{proof}

\subsubsection{Perimeter measure}\label{sperim}

Let $A\subset {\mathbb R}^2_{\geq 0}$ be a compact polygonal
region. 
Define the perimeter measure $\mm^P(A)$ by:
\begin{equation}
\label{ddeffi}
\mm^P(A) = \lim_{N\rarr \infty} \ \ \frac{1}{\sqrt{N}}\, \sum_{n\ge 1}
\sum_{\substack{T\in\E^{11}(n)\\ \\ (P_R(T),P_T(T))\in N A}}
\frac{1}{e^n\, (n-1)!} \,,
\end{equation}
where $P_R(T)$ and $P_T(T)$ denote the root and top perimeters
of an edge tree $T$ and $N A$ denotes the region $A$ scaled by 
a factor of $N$. 

\begin{pr}
\label{ddl}
We have:
$$\mm^P(A) =  \frac{1}{\sqrt{2\pi}}  \int_A
\frac{dx \, dy}
{(x+y)^{3/2}} $$
\end{pr}

\begin{proof}
It suffices to prove the Proposition for the sets $A$ of the
form
$$
A_{c,d}=\left\{(x,y), (x+y)\le c , \frac{x}{x+y} \le d \right\}\,,
\quad d\in [0,1] \,.
$$    
It is clear, that in \eqref{ddeffi} we can replace the
summation over $n\ge 1$ by the summation over $n \ge M$ for
any $M$, because the contribution of any particular value
of $n$ is suppressed by the factor $ \dfrac{1}{\sqrt{N}}$.
By Proposition \ref{xmml},  
choosing $M$ sufficiently large, we can make the distribution
of $\dfrac{P_R}{P_R+P_T}$ be arbitrarily close to the uniform distribution
on $[0,1]$, whence
$$
\mm^P(A_{c,d}) = d \, \mm^P(A_{c,1}) \,.
$$
The measure $\mm^P(A_{c,1})$, by Lemma \ref{jjjj}, just counts the
number of trees of size $\le cN$, or, more concretely, 
\begin{multline*}
\mm^P(A_{c,1}) = \lim_{N\to\infty}  
\frac{1}{\sqrt{N}} \sum_{n=1}^{cN} \frac{(n-1)\,n^{n-2}}{e^{n}\, (n-1)!} =\\
 \lim_{N\to\infty}  
\frac{1}{\sqrt{N}} \sum_{n=1}^{cN} \frac{1}{\sqrt{2\pi n}} =
\lim_{N\to\infty}  
\frac{1}{N} \sum_{n=1}^{cN} \frac{1}{\sqrt{2\pi (n/N)}} 
= \frac1{\sqrt{2\pi}}\, \int_0^c \frac{dt}{\sqrt{t}} \,, 
\end{multline*}
where the second equality uses the Stirling formula
(see e.g.\ \cite{Ap})
\begin{equation}
  \label{Stirl}
n! = 
\sqrt{2\pi n} \, \frac{n^n}{e^n} \,
\left(1+O\left(\tfrac{1}{n}\right)\right)  
\end{equation}
and the last equality is by the
 definition of the integral. 
This determines the measure $\mm^P$ uniquely and concludes
the proof. 
\end{proof}

\subsubsection{Independence of semiperimeters and root/top
components}

The analysis in Section \ref{semdis} can be repeated exactly
to obtain the
asymptotics of the normalized perimeter $\dfrac{P_R}{n}$
for trees with fixed root and top components. Concretely,
suppose the fixed root and top components have 
 $k$ and $l$ vertices respectively.
Then, by moving the point where the top component is attached to
the trunk component,  one sees that on this set $\psi(v_t)$
is uniformly distributed on the interval $\{k+1,\dots,n-l+1\}$. Hence,
$\dfrac{\psi(v_t)}{n}$ converges to the uniform distribution
on $[0,1]$. Using Corollary \ref{rootfin} we conclude that: 

\begin{pr}\label{indperroot}
In the limit $n\to\infty$, the normalized perimeter $\dfrac{P_R}{n}$ of
a random tree $T\in\E^{11}(n)$ is
independent of the root and top components of $T$ and,
in particular, independent of the valences of the root and top
vertices of $T$.   
\end{pr}

\subsubsection{Effect of relabeling the edges}\label{relab}

Recall that, by definition, the edges of an edge tree $T\in \E^{11}$ 
come with a bijective labeling
$$
\phi: E \to \{1,2,\dots,n-1\} \,.
$$
This labeling goes into the definition of the angle \eqref{angledef}
between two adjacent edges of $T$. 

Suppose we have a monotone injective map
$$
\sigma: \{1,2,\dots,n-1\} \to \{1,2,\dots,N\} \,, \quad N \ge n-1 \,,
$$
using which we modify the definition \eqref{angledef} as follows:
$$
\widetilde\measuredangle(e,e')=\frac{2\pi}{N}
(\sigma(\phi(e'))-
\sigma(\phi(e)))  \bmod 2\pi\,,
\quad  \widetilde\measuredangle(e,e')\in
(0,2\pi]\,,
$$
and, accordingly, we introduce modified semiperimeters $\widetilde P_R$ and
$\widetilde P_T$. Because the vertices which do not belong to
the trunk still contribute $1$ to $\widetilde P_R$ or
$\widetilde P_T$, respectively, 
we have
$$
\left|
\frac{\widetilde P_R}{n}  - \frac{P_R}{n} 
\right| \le \frac{|\tk T|}{n} \,. 
$$
Recall that by Corollary \ref{tk0} the right-hand side
converges to $0$ in distribution as $n\to\infty$.

\section{\bf Asymptotics of the Hurwitz numbers}
\label{fin2}
\subsection{Overview}
Let $\mu$ be a partition with $l$ distinct parts.
By Definition \ref{Hur2}, the Hurwitz number 
$H_{g,\mu}$ is a weighted  count of the branching graphs
on  $\Sigma_g$ with perimeter $\mu$.
By Proposition \ref{lhur}, the asymptotics of $H_{g,N\mu}$
as $N\to\infty$ recover the $l$-point function
$P_g(\mu_1,\dots,\mu_l)$ defined in \eqref{P_g}. 
In the present section, we compute these
asymptotics using results about random trees
obtained in Section \ref{st2}.

Instead of analyzing the asymptotics of 
$H_{g,N\mu}$ for particular partitions $\mu$, 
we will study the Hurwitz
asymptotics averaged over a neighborhood $U$ of
$\mu$ values, that is, the asymptotics of the number
\begin{equation}
\label{pmdp}
\sum_{\nu\in NU} H_{g,\nu}
\end{equation}
as $N\to\infty$. Here, the sum is over 
integral points $\nu$ of
$N U$ (the set $U$ scaled by a factor of
$N$). After weighting by functions of $\nu$, these
asymptotics will define a Hurwitz measure
$\mm_g(U)$ of $U$ (see Section \ref{shurm}). 
The Hurwitz measure $\mm_g$ will uniquely determine the
$l$-point function $P_g$.
The averaged data (\ref{pmdp}) 
arises naturally in the Laplace transform of the Hurwitz asymptotics required to recover
Kontsevich's series $K_g$ \eqref{kgen}.
We will find that the averaging leads to simplifications in the
asymptotic analysis.

The Hurwitz measure $\mm_g$ is analyzed by studying the distribution of
cell perimeters of a random branching graph with a large total
perimeter. The following
strategy will be used. 
The branching graphs of genus $g$ with $l$ cells 
are partitioned into finitely many homotopy classes
indexed by maps $G\in\bGth_{g,l}$ with vertices of
valence $3$ and higher. Accordingly, the Hurwitz
measure $\mm_g$ is decomposed into a finite sum of contributions
$\mm_G$ of homotopy classes.

In Section \ref{q11}, $\mm_G$ is proven to vanish unless 
$G$ is trivalent. For a trivalent $G$, 
$\mm_G$ is shown to be a push-forward under a linear map of a
product measures which is the product of the perimeter
measures $\mm^P$ (see Section \ref{sperim}) over the edges of $G$.
After the Laplace transform, we recover precisely
the contribution of $G$ to Kontsevich's combinatorial
model \eqref{dfgg}. This establishes Theorem \ref{grasy} and
completes the proof of Theorem \ref{kon}.

\subsection{Hurwitz measure}\label{shurm}
Let $A\subset {\mathbb R}^l_{\geq 0}$ be a compact polygonal
region.
Define the genus $g$ Hurwitz measure $\mm_g(A)$ by:
\begin{equation}
\label{dipo}
\mm_g(A) = \lim_{N\rarr \infty} \, \frac{1}{N^{3g-3+3l/2}} 
\sum_{\mu \in N A}  \frac{H_{g,\mu}}{e^{|\mu|}\ r(g,\mu)!}\,.
\end{equation}
Here, the sum is over integral points $\mu$ in $NA$ 
(the set $A$ scaled by the factor of $N$),
and
$$
r(g,\mu)=2g-2+|\mu|+\ell(\mu)
$$
is a number of simple ramifications of a genus $g$ Hurwitz cover
corresponding to the partition $\mu$. 

\begin{pr}
\label{ddll}
The Hurwitz measure is determined by:
$$\mm_g(A) = \int_A  H_g(x) \ dx_1 \cdots dx_l.$$
\end{pr}

\bpf
$H_g(x)$ is defined on rational points $x\in {\mathbb Q}^l_{>0}$
satisfying $x_i \neq x_j$ by:
$$H_{g}(x_1,\ldots,x_l) =
\lim_
{N\rarr \infty} \ \  \frac{1}{N^{3g-3+l/2}}\frac{H_{g,Nx}} 
{e^{d(Nx)} \
r(g,Nx)!}.$$
$H_g(x)$ is a polynomial function by Theorem \ref{rrrr}.
The Proposition then follows directly from Stirling's formula
\eqref{Stirl} and the definition of the Riemann integral.
\epf

The Laplace transformed measure $L\mm_g$ is then
determined as a function of $y_1,\ldots, y_l$ by the equivalent formulas:
\begin{eqnarray*}
L\mm_g(A) & = & \lim_{N\rarr \infty} \, \frac{1}{N^{3g-3+3l/2}} 
\sum_{\mu \in NA}  e^{-y\cdot \mu/N} \, \frac{
 H_{g,\mu}}{e^{|\mu|} r(g,\mu)!}  \\
& = & \int_A e^{-y\cdot x} H_g(x) \ dx_1 \cdots dx_l
\end{eqnarray*}
By construction, 
\begin{eqnarray*}
LH_g(y_1,\ldots,y_l) & = 
& \int_{{\mathbb R}^l_{\geq 0}} e^{-y\cdot x} \, H_g(x) \, dx_1 \cdots dx_l
\\  
& =& L\mm_g( {\mathbb R}^l_{\geq 0}).
\end{eqnarray*}

Our strategy, introduced in Section \ref{grasym}, is to express 
$L\mm_g(A)$ as a sum of contributions of the possible
homotopy types $G\in\bGth_{g,l}$. We define
\begin{equation}\label{LmG}
 \mm_G(A)  =  \lim_{N\rarr \infty} \, \frac{1}{N^{3g-3+3l/2}} 
\sum_{\mu \in NA} \frac{
 H_{G,\mu}}{e^{|\mu|}\, r(g,\mu)!} \,,  
\end{equation}
where the number $H_{G,\mu}$ is counting branching graphs
with homotopy type $G$, see Section \ref{grasym}.

The limit $\mm_G(A)$ will be shown to exist for all
$G\in\bGth_{g,l}$ and, in fact, shown to vanish
unless $G$ is trivalent. 

\subsection{Assembling branching graphs from edge trees}

Fix a homotopy type $G\in\bGth_{g,l}$ with $|E|$ edges.
Denote by $\bH_{G,\mu}$ the set of all branching graphs
with homotopy type $G$ and perimeter $\mu$. We will use  
the following procedure for enumerating all
elements of $\bH_{G,\mu}$. 

We will use the symbol
$\displaystyle\binom{r}{r_1,\dots,r_k}$ to denote {\em both}
the multinomial coefficient {\em and} the corresponding
set of $k$-tuples of subsets of an $r$-element set. 
Fix an arbitrary orientation of the edges of $G$. 
There exists a natural \emph{assembly map}
\begin{equation}\label{Asm}
\Asm_G:\bigsqcup_{r_1,\dots,r_{|E|}} 
\binom{r}{r_1,\dots,r_{|E|}} \times \prod_{i=1}^{|E|}
\E^{11}(r_i+1) \to  \bigsqcup_{|\mu|=r+2-2g-l} \bH_{G,\mu} 
\cup \{\bem\} \,,  
\end{equation}
which is defined as follows. 

Let $e_1,\dots,e_{|E|}$ be the edges of $G$ and let
$$
T_1,\dots,T_{|E|}
$$
be an $|E|$-tuple of edge trees. First, we 
replace, preserving the order, the edge markings of each tree by
a subset of the set $U_r$ of $r$th roots of unity  according to the given element
of $\displaystyle \binom{r}{r_1,\dots,r_{|E|}}$.

After that, we 
replace each oriented edge $e_i$ of $G$ by the corresponding edge
tree $T_i\in \E^{11}(r_i+1)$ 
in such a way that the root vertex of $T_i$ is identified
with the initial
vertex of $e_i$ and the top vertex of $T_i$ is identified
with the final vertex of $e_i$.

This replacement is done so that
 at the vertices $v$ of $G$
the edges coming from the same tree $T_i$ have consecutive places in
the clockwise
order around $v$ with the trunk edge being the
last one. 

\begin{figure}[!hbt]
\centering
\scalebox{.6}{\includegraphics{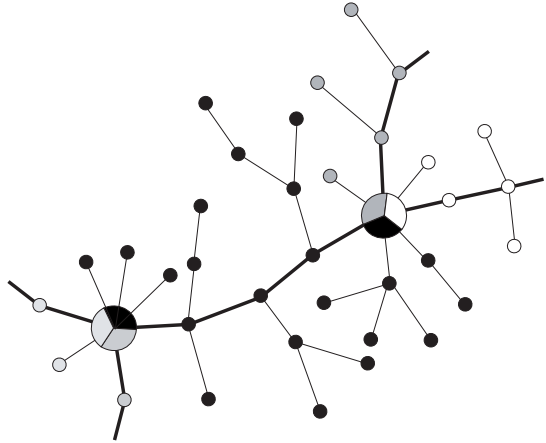}}
\caption{The tree from Figures \ref{fig7} and \ref{fig8} as part of an assembly}
\label{fig9}
\end{figure}

This procedure is illustrated in Figure \ref{fig9}
where it is shown how the tree from Figures \ref{fig7} and \ref{fig8}
may be used in the assembly of a branching graph. 
In Figure \ref{fig9}, the vertices of different trees are shown
in different color and those vertices which are shared by several
trees (the ones which correspond to the vertices of the graph $G$)
are painted accordingly. Also observe how the trunks of the trees in
Figure \ref{fig9} form the edges of the homotopy type graph $G$.

The resulting graph $H$ is a branching graph if 
the edge labels of $H$ at each vertex $v$ respect the cyclic order
of the roots of unity.
If $H$ is a branching graph, then define $\Asm_G=H$.
Otherwise,
the assembly result is declared a failure, indicated by the formal symbol 
$\Asm_G=\bem$.  

The group $\Aut(G)$ acts\footnote{
This action
is always free, whereas the action on just edges of $G$ may not be as
the $(g,n)=(1,1)$ example shows.}
on the oriented edges of $G$. This action makes $\Aut(G)$ act
on the domain of the assembly map by permutation of factors in
the Cartesian product and reversal of the root/top choice 
in individual factors. The following property of the 
assembly map is obvious from its
construction: 

\begin{pr}\label{prasm}
The assembly map is surjective and for any $H\ne\bem$
the preimage $\Asm_G^{-1}(H)$ is a single $\Aut(G)$-orbit. 
\end{pr}

We will call the edge trees in $\Asm_G^{-1}(H)$ the \emph{edge parts}
of a branching graph $H$, ignoring the minor ambiguity coming
from the action of $\Aut(G)$.

\subsection{Vanishing for non-trivalent graphs}
\label{q11}

Fix a homotopy type $G$ and let $|E|$ be the number
of edges in $G$. From the equation
$$
r=2g-2+|\mu|+l
$$
and Corollary \ref{treeenum} we have
\begin{multline}\label{r!r!}
 \frac{1}{e^{|\mu|} \, r!}
\binom{r}{r_1,\dots,r_{|E|}} \times \prod_{i=1}^{|E|} 
\left|\E^{11}(r_i+1)\right| = \\
e^{|E|+l+2g-2} \, \prod_{i=1}^{|E|} 
\frac{\left|\E^{11}(r_i+1)\right|}{e^{r_i+1}\, r_i!} 
\sim 
\frac{e^{|E|+l+2g-2}}{(2\pi)^{|E|}}\,  \prod \frac1{\sqrt{r_i}}\,, 
\end{multline}
as $r_i\to\infty$. By Proposition \ref{prasm}, this implies
that
$$
\sum_{|\mu|\le N} \frac{1}{e^{|\mu|}\, r!} \, H_{G,\mu}
= O\left(\prod_{i=1}^{|E|}\sum_{r_i=1}^N \frac1{\sqrt{r_i}}\right)=
O\left(N^{|E|/2}\right)\,, \quad N\to\infty \,.
$$
It follows that the limit \eqref{LmG} vanishes unless 
$$
|E|=6g-6+3l\,,
$$
which is equivalent to $G$ being trivalent. Thus, we 
have established the following:

\begin{pr}\label{prvanntrv}
If the homotopy type $G$ is not trivalent, then $\mm_G=0$. 
\end{pr}

Further, for a trivalent graph $G$ we can assume the edge trees
which participate in the assembly map to be arbitrarily 
large. This is because the contribution of edge trees of any
fixed size to $\mm_G$ obviously vanishes  for the same reason
as above.

This puts us
in the asymptotic regime of large random edge trees, which was
considered in Section \ref{st2}. In particular, in this regime
$$
\left|\Asm_G^{-1}(H)\right| = |\Aut(G)|\,,
$$
for a typical branching graph $H$. This is because the probability of 
having two isomorphic edge parts or an edge part which has an
automorphism permuting root and top clearly goes to zero
as the size of the edge parts goes to infinity. 

\subsection{Probability of assembly failure}

In particular, let us compute the probability that 
for large random edge trees the assembly \eqref{Asm} 
will end in failure $\bem$. In other words, we want
to compute the probability that the cyclic order condition will be violated
at one of the vertices $v$ of a trivalent graph $G$. 

Suppose we have three random disjoint sequences $X_i$, $i=1,2,3$ of 
elements
of a some cyclically ordered set. The conditional
probability that the concatenation
$$
X_1,X_2,X_3
$$
is cyclically ordered, given that each $X_i$ is cyclically ordered and
$|X_i|=k_i$, $i=1,2,3$, 
is easily seen to be equal to 
$$
\frac {(k_1-1)! \,(k_2-1)!  \,(k_3-1)!} 
{(k_1+k_2+k_3-1)!} \,.
$$
In our case, the cyclically ordered set is the 
 set $U_r$ of $r$th  roots 
of unity and the sequences $X_i$ are the labels of the
edges incident to the 3 root/top vertices $v_1$, $v_2$, $v_3$ that
are glued together at $v$.

The probability that the valence if $v_i$ in the corresponding
tree equals $k_i$ approaches, by Corollary \ref{poisslimit}, the limit
$e^{-1}/(k_i-1)!$ as the size of tree goes to infinity. Therefore,
the probability of failure at a given vertex $v$ converges to
\begin{multline}
e^{-3}\sum_{k_1,k_2,k_3=1}^\infty  \frac{1}{(k_1+k_2+k_3-1)!} = \\
e^{-3} \sum_{k=3}^{\infty}\, \sum_{k_1+k_2+k_3=k} \frac{1}{(k-1)!} =
\frac{e^{-3}}2 \sum_{k=3}^\infty \frac{(k-1)(k-2)}{(k-1)!} = 
\frac{e^{-2}}2 \,.
\end{multline}

Also, we see that at each vertex $v$, the probability
to fail depends only on the three valences $k_1$, $k_2$, $k_3$ 
involved, and hence, in
the limit of large random edge trees, failures at vertices of $G$
become independent events by Corollary \ref{roottopcomp1}.  

Since a trivalent graph $G$ has $4g-4+2l$ vertices, we obtain
the following conclusion (the second assertion follows from
Proposition \ref{indperroot}) 

\begin{pr} \label{asmbfail}
For any trivalent homotopy type $G$, the probability
of assembly failure goes to 
$$
e^{-8g+8-4l} \, 2^{-4g+4-2l}
$$ 
as
the sizes of all edge trees go to infinity. Further, assembly
failure is asymptotically independent of the normalized 
semiperimeters of the edge trees involved. 
\end{pr}

\subsection{Computation of the Hurwitz measure}

By definition \eqref{LmG},
the Hurwitz measure $\mm_G(A)$ involves the asymptotics
of the weighted number of branching
graphs $H$ of homotopy type $G$ such that
$$
\frac{\mu}{N} \in A\,,
$$
where $\mu$ is the perimeter of $H$.

Let $D$ be a cell of $H$. The boundary $\partial D$,
followed in the clockwise direction, 
is a sequences of edges
$$
e_1,e_2,\dots,e_s\,, \quad e_i\in E(H) \,.
$$
The perimeter of $D$ is, by definition,
the following sum
\begin{equation}\label{perD}
\per(D) = 
\frac1{2\pi}\sum_{k=1}^s 
\arg\left( \frac{\gamma(e_{k})}{\gamma(e_{k+1})} \right).
\end{equation}
where $\gamma:E(H)\to U_r$ is the labeling of the edges
of $H$ by roots of unity, the argument takes values in
$(0,2\pi]$, and $e_{s+1}=e_1$. 

For most terms in \eqref{perD} both
$e_k$ and $e_{k+1}$ belong to the same edge part of $H$.
The only exception are the terms corresponding to
the vertices of $G$ on the
boundary of $D$.  The contribution of these exceptional terms
to $\per(D)$ is bounded by $4g-4+2l$ because the contribution
of each of the  $4g-4+2l$ vertices of $G$ to $\per(D)$ is bounded
by $1$. Since we are interested
in the distribution of $\dfrac{\per(D)}{N}$ as $N\to\infty$,
this contribution may be ignored. 

Further, we can substitute the contribution of each edge part
of $H$ by the semiperimeter of the corresponding edge tree $T$. The
difference between these two numbers is that the first is computed using the
labeling $\gamma$ of $E(T)$ by the $r$th roots of unity,
whereas the second uses the labeling
$$
\phi: E(T) \to \{1,2,\dots,|E(T)|\} 
$$
which is a part of the structure of an edge tree. By the
results of Section \ref{relab} the difference between these
two numbers is of 
the order of magnitude $\sqrt{N}$. In particular, the
probability of having a difference of size $N^{1/2+\epsilon}$ 
goes to zero for any $\epsilon>0$. This means that the effect
of this difference on $\dfrac{\per(D)}{N}$ 
is negligible in the $N\to\infty$ limit. 

It follows that, asymptotically, $\dfrac{\per(D)}{N}$ is a 
sum of independent random variables which are the
normalized semiperimeters of the edge parts of $H$ along
the boundary of $D$. Recall that the distribution of 
the normalized semiperimeters of a random edge tree is governed
by the perimeter measure $\mm^P$, which was studied in
Section \ref{sperim}. 
Together with \eqref{r!r!} and Proposition
\ref{asmbfail}, we conclude the following: 

\begin{pr}
For any $G\in\bGt_{g,l}$, we have
$$
\mm_G = \frac{2^{-4g+4-2l}}{|\Aut(G)|}\, 
\left(\asm_G\right)_*\left(\bigotimes_{E(G)} \mm^P\right)
$$  
where the product of perimeter measures is over all edges
of $G$ and $\asm_G$ is the linear map which takes the
normalized semiperimeters to their sums along the boundaries of the
cells of $G$.  
\end{pr}

This result can be more conveniently stated in terms of the 
Laplace transform 
$$
L\mm_G (y_1,\dots,y_l) = \int_{\R^l_{>0}}
e^{-y\cdot x} \, \mm_G(dx) \,,
$$ 
for which it implies the following factorization
$$
L\mm_G (y) = \frac{2^{-4g+4-2l}}{|\Aut(G)|}\, \prod_{e\in E(G)}
L\mm^P (y_{i(e)},y_{j(e)}) \,,
$$
where $i(e)$ and $j(e)$ are the numbers of the cells of $G$ that 
the edge $e$ 
separates and $L\mm^P$ is the Laplace transform of the
perimeter measure $\mm^P$. 

It remains, therefore, to compute $L\mm^P$ which by Proposition
\ref{ddl} equals the following integral
$$
L\mm^P (y_1,y_2)=
\frac1{\sqrt{2\pi}} \iint_0^\infty \frac{e^{-x_1 y_1 - x_2 y_2}}
{(x_1+x_2)^{3/2}}\, dx_1 dx_2 \,,
$$ 
where $y_1,y_2 > 0$. 

Making a change of variables
$$
x_1+x_2 = u\,, \quad x_1-x_2 = v \,,
$$
and integrating out $v$, we obtain
$$
L\mm^P (y_1,y_2) = - \frac1{\sqrt{2\pi}} \, \frac1{y_1-y_2}
\, \int_0^\infty \left(e^{-y_1 u} - e^{-y_2 u}\right) \, \frac{du}{u^{3/2}}
\,.
$$
For $\Re \alpha > -1$ we have
$$
\int_0^\infty \left(e^{-y_1 u} - e^{-y_2 u}\right) \, u^{\alpha-1} \,
du = \Gamma(\alpha) \, \left(\frac1{y_1^{\alpha}} -  \frac1{y_2^{\alpha}}
\right)
$$
which for $\Re \alpha > 0$ follows from the definition of the 
$\Gamma$-function and can be extended to $\Re \alpha > -1$ by
analytic continuation because the integral remains 
absolutely convergent. Since
$$
\Gamma(-1/2) = - 2 \sqrt\pi\,,
$$
plugging in $\alpha=-1/2$, we obtain
$$
L\mm^P (y_1,y_2) = \frac{\sqrt{2}}{\sqrt{y_1}+\sqrt{y_2}} \,.
$$ 
Since a trivalent graph $G$ has $6g-6+3l$ vertices, it follows
that 
$$
L\mm_G (y) = \frac{2^{2g-2+l}}{|\Aut(G)|}\, \prod_{e\in E(G)}
\frac{1}{\sqrt{2y_{i(e)}}+\sqrt{2y_{j(e)}}} \,.
$$

Finally, summing over all homotopy types $G$ and using the
vanishing for nontrivalent homotopy types established in
Proposition \ref{prvanntrv}, we obtain the following:

\begin{pr}
We have
$$
LH_g(y_1,\dots,y_l) = \sum_{G\in\bGt_{g,l}} \, 
\frac{2^{2g-2+l}}{|\Aut(G)|}\, \prod_{e\in E(G)}
\frac{1}{\sqrt{2y_{i(e)}}+\sqrt{2y_{j(e)}}} \,,
$$
where the product is over all edges of a trivalent 
map $G$ and $i(e)$ and $j(e)$ are the numbers of the edge $e$ 
separates. 
\end{pr}

This Proposition completes the proof of Theorem \ref{grasy}.
Theorems \ref{connn} and \ref{grasy}  imply Theorem \ref{kon}
which, therefore,  is established.

\subsection{Connection with the edge-of-the-spectrum matrix model}

After studying the asymptotic enumeration of branching
graphs on $\Sigma_g$ with $l$ cells, 
we see that the problem is exactly 
parallel  to the
asymptotic enumeration of simple maps on $\Sigma_g$ with $l$ cells
carried out in \cite{O2}. 

As in the case of  branching graphs, there
exist only finitely many homotopy types of simple maps, of which only
the trivalent homotopy types make a nonvanishing contribution 
to the asymptotics. The cell perimeters of a map are now
exactly equal to the sums of semiperimeters of  edge parts
along the cell boundaries. An edge part of map is an unmarked
 planar tree
with a choice of a root and top vertex (the semiperimeter
distribution is easily seen to be asymptotically uniform). 

This complete parallelism explains the equality \eqref{Pmap} and therefore
explains the connection between intersection theory of $\mgn$
and the matrix model \eqref{ESMM}.

\subsection{Lower order asymptotics}
The lower order terms (in $N$) of the
asymptotics of $H_{g,N\mu}$ 
govern Hodge integrals on $\overline{M}_{g,l}$ with integrand linear
in the $\lambda$ classes.
It appears quite difficult to extract lower order
asymptotics from the random tree analysis. However, the
{\em lowest} order term, related to the $\lambda_g$ integrals
$$\langle \tau_{k_1} \cdots \tau_{k_l} \lambda_g\rangle_g =
\int_{\overline{M}_{g,l}} \psi_1^{k_1} \cdots \psi_l^{k_l} \lambda_g,$$
 is well-understood from a different perspective.

The $\lambda_g$ integrals arise in the degree 0 sector of the
Virasoro conjecture for an elliptic target curve. In \cite{GeP}, the Virasoro
conjecture for this degree 0 sector was shown to be equivalent to the following
equation:
\begin{equation}
\label{sdflkj}
\langle \tau_{k_1} \cdots \tau_{k_l} \lambda_g\rangle_g  = \binom{2g-3+l}{
k_1, \ldots, k_l} \langle \tau_{2g-2} \lambda_g\rangle_g,
\end{equation}
where $\langle \tau_{-2} \lambda_0\rangle_0=1$.
The $\lambda_g$ conjecture (\ref{sdflkj}) was later proven in \cite{FaP2} via
virtual localization techniques (independent of the Hurwitz connection
developed here). The integrals $\langle \tau_{2g-2} \lambda_g\rangle_g$ are
determined by:
\begin{equation*}
\label{qqq}
\sum_{g\geq 0}t^{2g} \langle \tau_{2g-2} \lambda_g \rangle_g
= \Big( \frac{t/2}{\sin(t/2)} \Big),
\end{equation*}
proven in \cite{FaP1}.

Hodge integrals over the moduli space of curves are intimately related
to Gromov-Witten theory via virtual localization, Virasoro constraints, 
Toda equations, and Mirror symmetry. Additional Hodge integral formulas
and predictions may be found in \cite{BP,Fa,FaP1,FaP2,FaP3,GeP,MSS,O1,P1,P2}.

\appendix
  
\section{Degeneration formulas for Hurwitz numbers}
\label{degdeg}
Classical recursive formulas for $H_{g,\mu}$ are obtained
by studying the degenerations of covers as a finite branch point is
moved to $\infty$. The recursions provide
an elementary (though combinatorially complex) 
method of calculating $H_{g,\mu}$. 
We derive the degeneration formulas here
from Definition 2 of the Hurwitz numbers following a suggestion
of R. Vakil. There are very many different proofs of these
formulas (see, for example, \cite{GoJVa,IP,LR,LZZ}).

A Hurwitz cover $\pi:C \rarr \proj^1$ together with a marking of the
fiber $\pi^{-1}(\infty)$ is a {\em marked} Hurwitz cover.
Let $H^*_{g,\mu}$ denote the automorphism
weighted count of marked Hurwitz covers with ramification $m_i$
at the 
$i^{th}$ marked point. We find:
$$H^*_{g,\mu} = |\Aut(\mu)| \cdot H_{g,\mu}.$$
By Definition 2, $H^*_{g,\mu}$ equals a count of
distinct $\mu$-graphs $H^*$ with {\em marked cells} 
on $\Sigma_g$ (weighted by $1/|\Aut(H^*)|$).
The Hurwitz numbers $H^*_{g,\mu}$ are more
convenient for the degeneration formulas.

Let $\mu=(m_1,\ldots,m_l)$ be a partition with positive parts.
The following partition terminology will be needed:
\begin{enumerate}
\item[$\bullet$]
$\mu-m_i$ equals the partition (possibly empty) of length $l-1$
 obtained by deleting $m_i$. 
\item[$\bullet$]
$\mu(m_i+m_j)$ equals
the partition of length $l-1$ obtained by combining $m_i$ and $m_j$.
\item[$\bullet$]
$\mu(a_1+a_2=m_i)$ equals the partition of length $l+1$ obtained
by splitting $m_i$ into positive parts $a_1$ and $a_2$.
\item[$\bullet$]
$\mu+a$ equals the partition of length $l+1$ obtained adding a positive part $a$.
\end{enumerate}
Finally, let $\mu_1 + \mu_2$ denote the union of the
partitions $\mu_1$ and $\mu_2$.

As in Section \ref{hur},
let $r(g,\mu)=2g-2+|\mu|+\ell(\mu)$ be the number of finite 
branch points of the Hurwitz covers counted by $H^*_{g,\mu}$.
If $r(g,\mu)$ vanishes, then $g=0$ and $\mu=(1)$. In this
case, $H^*_{0,(1)}=1$.  The Hurwitz numbers $H^*_{g,\mu}$
are determined recursively by the following Theorem.

\begin{tm}
\label{dddd}
Let $r(g,\mu)>0$. The degeneration relation holds:
\begin{eqnarray*}
H^*_{g,\mu} & = & \ \ \
\sum_{i\neq j} \ \ \ \  \ \frac{m_i+m_j}{2} H^*_{g,\mu(m_i+m_j)} \\
& & + 
   \sum_i  \sum_{a_1+a_2=m_i} \frac{a_1a_2}{2} H^*_{g-1, \mu(a_1+a_2=m_i)} \\
& &
   +      \sum_i \sum_{a_1+a_2=m_i}\ \sum_{ g_1+g_2=g}\ 
\sum_{ \mu_1+\mu_2=\mu-m_i} 
\epsilon \frac{a_1a_2}{2} 
  H^*_{g_1,\mu_1+a_1} H^*_{g_2,\mu_2+a_2},
\end{eqnarray*}
where $\epsilon$ denotes a binomial coefficient in the last sum:
$$\epsilon= \binom{r(g,\mu)-1}{r(g_1,\mu_1+a_1)} .$$
\end{tm}

\bpf The
degeneration of a Hurwitz cover as a branch point 
is moved to $\infty$ corresponds
simply to edge removal for the associated $\mu$-graphs.

Let $H^*$ be a $\mu$-graph with marked cells on $\Sigma_g$.
Let $r=r(g,\mu)$. Let $U_r$ be the set of $r^{th}$ roots of unity marking
the edges.
There are three possibilities for the graph $X$ obtained 
after removal of the edge $e$ marked by the unit $1\in U_r$.

\vspace{+5pt}
\noindent {\em Case I.} 
The edge $e$ separates two distinct cells of $H^*$ with markings $i\neq j$.
Then,
$X$ is canonically a $\mu(m_i+m_j)$-graph with marked cells on $\Sigma_g$.
The edge markings of $X$ lie in $U_r \setminus\{1\}$. 

Conversely, let $X$ be a $\mu(m_i+m_j)$-graph with marked cells on $\Sigma_g$.
Let the edge markings of $X$ lie in set $U_r \setminus \{1\}$.  
Let $D$ be the cell corresponding to the part $(m_i+m_j)$.
There are
$m_i+m_j$ distinct ways an edge $e$ with marking $1$ 
may be added  which separates $D$ into two cells of perimeters $m_i$ and $m_j$ {\em and}
respects the edge orientation
conditions.

\vspace{+5pt}
\noindent {\em Case II.} The two sides of $e$ bound the same
cell of $H^*$ {\em and} $e$ is not a  disconnecting edge.
Then, $X$ is canonically a $\mu(a_1+a_2=m_i)$-graph with marked
cells on $\Sigma_{g-1}$. Conversely, there are $a_1a_2$ ways to
add $e$ to $X$ to recover a $\mu$-graph with marked cells on
$\Sigma_g$. 

\vspace{+5pt}
\noindent {\em Case III.} The two sides of $e$ bound the
same cell $H^*$ {\em and} $e$ is a disconnecting edge.
Then, $X=X_1 \cup X_2$ is the union where $X_i$ is
a $\mu_i+a_i$-graph with marked cells on $\Sigma_{g_i}$.
Conversely, there are $a_1a_2$ ways to
add $e$ to $X$ to recover a $\mu$-graph with marked cells on
$\Sigma_g$. 

\vspace{+5pt}
The degeneration formula follows from counting these three
cases (weighted by the possible locations of $e$).
\epf

The degeneration formulas may be viewed as a first
geometric approach to the Hurwitz numbers.
Unfortunately, a direct analysis  of $H_{g,\mu}$ 
via Theorem \ref{dddd} appears combinatorially
difficult. More efficient recursive strategies for the Hurwitz
have been found (see \cite{FanP,GoJV}), but these formulas
are genus dependent.

\section{Integral tables}
\label{tabtab}
Hodge integrals on $\overline{M}_{g,n}$ are {\em primitive}
if neither the string or dilaton equation may be applied.
With the exception of $\langle \tau_0^3\rangle_0$ and
$\langle \tau_1 \rangle_1$, the primitive condition is equivalent to the absence of 
$\tau_0$ and $\tau_1$ factors in the integrand.
The first table contains all primitive Hodge integrals with a single
$\lambda$ class for $g\leq 2$.

$$
\begin{array}{||c||c||}
\hline
g=0 & \langle \tau_0^3\rangle_0 =1 \\ \hline
g=1 & \langle \tau_1 \rangle_1=1/24, \ \langle \lambda_1 \rangle_1=1/24 \\ \hline
g=2 & \langle\tau_4\rangle_2=1/1152, \ \langle \tau_3 \tau_2 \rangle_2=29/5760, \ \langle \tau_2^3\rangle_2=7/240 \\ \hline
&  \langle \tau_3 \lambda_1\rangle_2=1/480, \ \langle \tau_2^2 \lambda_1 \rangle_2=5/576 \\
\hline 
&  \langle \tau_2 \lambda_2\rangle_2 =7/5760 \\ \hline
\end{array}
$$

The second table contains Hurwitz numbers $H_{g,\mu}$ for 
$g\leq 2$ and partitions $\mu$ satisfying $|\mu|\leq 4$.

$$
\begin{array}{||c||c|c|c|c|c|c|c||}
\hline 
H_{g,\mu}  &(1)& (2)&(1,1)&(3)&(2,1)&(1,1,1) \\ \hline 
g=0 & 1&1/2 &1/2 &1 & 4 & 4    \\ \hline
g=1 &0 &1/2 &1/2 & 9 &  40 & 40    \\ \hline
g=2 &0 &1/2 &1/2 &  81 & 364 & 364   \\ \hline
\end{array}
$$

$$
\begin{array}{||c||c|c|c|c|c|c|c||}
\hline 
H_{g,\mu}  &(4)& (3,1)&(2,2)&(2,1,1)&(1,1,1,1) \\ \hline 
g=0 &4 & 27& 12 & 120 & 120   \\ \hline
g=1 & 160& 1215 & 480 & 5460&5460     \\ \hline
g=2 &5824 &45927 &17472 &206640  &206640     \\ \hline
\end{array}
$$

\vspace{+10 pt}
\noindent
Department of Mathematics \\
UC Berkeley \\
Berkeley, CA 94720\\
okounkov@math.berkeley.edu \\

\vspace{+10 pt}
\noindent
Department of Mathematics \\
\noindent California Institute of Technology \\
\noindent Pasadena, CA 91125 \\
\noindent rahulp@cco.caltech.edu

\begin{thebibliography}{99}

\bibitem{ASV}
M.~Adler, T.~Shiota, P.~van Moerbeke,
\emph{Random matrices, Virasoro algebras, and 
non-commutative KP}, Duke Math.\ J.\ \textbf{94} (1998), 
379--431. 




\bibitem{Al1}
D.~J.~Aldous,
\emph{The continuum random tree I},
Ann.~Prob.\ \textbf{19} (1991), 1--28.


\bibitem{Al2}
D.~J.~Aldous,
\emph{The continuum random tree II: an overview},
in \emph{Stochastic Analysis}, edited by M.~T.~Barlow and
N.~H.~Bingham, Cambridge University Press, 1991, 23--70.


\bibitem{AlPit}
D.~J.~Aldous and J.~Pitman,
\emph{Tree-values Markov vhain derived from Galton-Watson
processes}, Ann.\ Inst.\ Henri Poincar\'e
\textbf{34} (1998), no.\ 5, 637--686.


\bibitem{Ar} V.\ I.\ Arnold, 
{\em Topological classification of complex trigonometric polynomials and the combinatorics of graphs with an identical number of vertices and edges},
 Funct.\ Anal.\ Appl.\ \textbf{30} (1996), no.\ 1, 1--14.  



\bibitem{Ap} T.\ Apostol, {\em 
Mathematical analysis}, Addison-Wesley: 
1974. 



\bibitem{AtBo} M.\ Atiyah and R.\ Bott, {\em The moment map
and equivariant cohomology}, Topology {\bf 23} (1984), 1-28.



\bibitem{B} K.\ Behrend, {\em Gromov-Witten invariants
in algebraic geometry}, Invent.\ Math.\ {\bf 127} (1997), 601-617.

\bibitem{BFan} K.\ Behrend and B.\ Fantechi, {\em The intrinsic normal
cone},  Invent.\ Math.\ {\bf 128} (1997), 45-88.

\bibitem{Bil}
P.~Billingsley, 
\emph{Convergence of probability measures},
John Wiley \& Sons, New York, 1999. 


\bibitem{BI} P.~Bleher and A.~Its, Talk at the Random Matrix conference
at MSRI, May 1999. 

\bibitem{Bo}
E.~Borel, 
\emph{
Sur l'emploi du th\'eor\`eme de Bernoulli pour faciliter le calcul 
d'une infinit\'e de coefficients. Application au probl\`eme de 
l'attente \`a un guichet},
C.~R.~Acad.\ Sci.\ Paris \textbf{214} (1942). 452--456. 


\bibitem{BK} E.~Br\'ezin and V.~Kazakov, \emph{
Exactly solvable field theories of closed strings},
Phys.\ Let.\ \textbf{B236} (1990), 144-150.

\bibitem{BP} J.~Bryan and R.~Pandharipande, {\em  BPS states of curves in
Calabi-Yau 3-folds}, preprint 2000.

\bibitem{De}
J.~D\'enes,
\emph{
The representation of a permutation as the product of a minimal number of transpositions, and its connection with the theory of graphs},
Magyar Tud.\ Akad.\ Mat. Kutat\'o Int.\ K\"ozl.\ \textbf{4} (1959),
63--71. 

\bibitem{dF} P.\ Di Franceso, {\em 2-d quantum gravities and topological
gravities, matrix models, and integrable differential systems}, in
{\em The Painlev\'e property: one century later}, (R.\ Conte, ed.), 229-286,
Springer: New York, 1999.

\bibitem{DFGZ} P.\ Di Franceso, P.~Ginzparg, and J.~Zinn-Justin,
\emph{2D Quantum gravity and random matrix models},
Phys.\ Rep.\ \textbf{254} (1995). 1--131. 

\bibitem{DIZ}  P.\ Di Franceso, C.~Itzykson, and J.~B.~Zuber,
\emph{Polynomial averages in the Kontsevich model}, Comm.\ Math.\ Phys.,
\textbf{151} (1993), 193--219. 




\bibitem{Dou} M.~Douglas, 
\emph{Strings in less than one dimension and generalized KP hierarchies},
Phys.\ Let.\ \textbf{B238} (1990) 

\bibitem{DS}
M.~Douglas and S.~Schenker, 
\emph{Strings in less than one dimension},
Nucl.\ Phys.\ \textbf{B335} (1990).



\bibitem{EdGra} D.\ Edidin and W.\ Graham, {\em Equivariant intersection
       theory}, Invent.\ Math.\ {\bf 131} (1998), 595-634. 



\bibitem{EgHX} T.\ Eguchi, K.\ Hori, and C.-S.\ Xiong, {\em
Quantum cohomology and Virasoro algebra}, Phys.\ Lett.\ {\bf B402} (1997),
71-80.


\bibitem{EgY} T.\ Eguchi and S.-K.\ Yang, {\em
The topological $CP^1$ model and the large-N matrix integral},
Mod.\ Phys.\ Lett.\ {\bf A9} (1994), 2893-2902.


\bibitem{ELSV} T.\ Ekedahl, S.\ Lando, M.\ Shapiro, and A.\ Vainshtein,
{\em On Hurwitz numbers and Hodge integrals}, 
preprint 2000.




\bibitem{Fa}  C.\ Faber, {\em A conjectural description of the 
tautological ring of the moduli space of curves},
in {\em Moduli of Curves and Abelian Varieties
(The Dutch Intercity Seminar on Moduli) } (C.\ Faber and
E.\ Looijenga, eds.),
109-129, Aspects of Mathematics E33, Vieweg: Wiesbaden, 1999.

\bibitem{FaP1}  C.\ Faber and R.\ Pandharipande, {\em
Hodge integrals and Gromov-Witten theory}, Invent.\ Math.\ 
{\bf 139} (2000), 173-199.


\bibitem{FaP2}  C.\ Faber and R.\ Pandharipande, {\em
Hodge integrals, partition matrices, and the $\lambda_g$ conjecture},
preprint 1999.

\bibitem{FaP3}  C.\ Faber and R.\ Pandharipande, {\em
Logarithmic  series and Hodge integrals in the tautological
ring}, preprint 2000.

\bibitem{FanP}  B.\ Fantechi and R.\ Pandharipande, {\em
Stable maps and branch divisors}, preprint 1999.

\bibitem{FGOR} Ph.~Flajolet,
Z.~Gao, Zhi-Cheng, A.~Odlyzko, and B.~Richmond, 
{\em The distribution of heights of binary trees and other simple trees},
 Combin.\ Probab.\ Comput.\ \textbf{2} (1993), no.~2, 145--156. 

\bibitem{Fu} W.\ Fulton, {\em Intersection theory}, Springer-Verlag: Berlin, 1998.

\bibitem{FP}
W.\ Fulton and R.\ Pandharipande,
{\em Notes on stable maps and quantum cohomology}, in
{\em Proceedings of symposia in pure  mathematics: Algebraic Geometry
Santa Cruz 1995} Volume 62, Part 2 (J.\ Koll\'ar, R.\ Lazarsfeld, and D.\ Morrison, eds.),
45-96, AMS: Rhode Island, 1997.

\bibitem{GMa} S.\ I.\ Gelfand and Y.\ Manin, {\em Methods of homological 
algebra}, Springer-Verlag: Berlin,
1996. 

\bibitem{Ge1} E.\ Getzler, in preparation.

\bibitem{Ge2}  E.\ Getzler, {\em The Virasoro conjecture for
Gromov-Witten invariants}, preprint 1998.



\bibitem{GeP} E.\ Getzler and R.\ Pandharipande, {\em Virasoro
constraints and the Chern classes of the Hodge bundle},
Nucl.\ Phys.\ {\bf B530} (1998), 701-714.

\bibitem{GoJVa} I.\ Goulden, D.\ Jackson, and A.\ Vainstein, {\em
The number of ramified coverings of the sphere by the torus and surfaces of
higher genera}, preprint 1999.


\bibitem{GoJV} I.\ Goulden, D.\ Jackson, and R.\ Vakil, {\em
The Gromov-Witten potential of a point, Hurwitz numbers,
and Hodge integrals}, preprint 1999.

\bibitem{GrKP} T.\ Graber, A.\ Kresch,  and R.\ Pandharipande, in preparation.


\bibitem{GrP} T.\ Graber and R.\ Pandharipande, {\em Localization
of virtual classes}, Invent.\ Math.\ {\bf 135} (1999), 487--518.

\bibitem{GrV} T.\ Graber and R.\ Vakil, {\em Hodge integrals and
Hurwitz numbers via virtual localization},
preprint 2000.

\bibitem{Gro} M.~Gromov, {\em
Pseudoholomorphic curves in symplectic manifolds}, Invent.\ Math.\ {\bf 82} (1985), 
307--347. 


\bibitem{GM1} D.~Gross and A.~Migdal,
\emph{Nonperturbative two-dimensional quantum gravity},
Phys.\ Rev.\ Lett.\ \textbf{64} (1990)

\bibitem{GM2} D.~Gross and A.~Migdal,
\emph{A nonperturbative treatment of two-dimensional quantum gravity},
Nucl.\ Phys.\ \textbf{B340} (1990), 333--365.

\bibitem{Hu} A.\ Hurwitz, {\em \"Uber die Anzahl der
Riemann'schen Fl\"achen mit gegebenen Verzweigungspunkten}, Math.
Ann.\ {\bf 55} (1902), 53-66.

\bibitem{IP} E.\ Ionel and T.\ Parker, {\em Relative Gromov-Witten
invariants}, preprint 1999.

\bibitem{IZ} C.~Itzykson, and J.~B.~Zuber,
\emph{Combinatorics of the modular group. 2.\ The Kontsevich
integral}, Int.\ J.\ Mod.\ Phys.\ \textbf{A7} (1992), 1--23. 

\bibitem{Iv} B.\ Iverson, {\em A fixed point formula for
action of tori on algebraic varieties}, Invent.\ Math.\ {\bf 16} (1972),
229-236.


\bibitem{Jo}
A.~Joyal,
\emph{
Une th\'eorie combinatoire des s\'eries formelles},
Adv.\ in Math.\ \textbf{42} (1981), no.\ 1, 1--82.

\bibitem{KS}
V.~Kac and A.~Schwarz,
\emph{Geometric interpretation of the partition function
of 2D gravity}, Physics Letters B \textbf{257} (1991),
329--334.




\bibitem{KiP} B.\ Kim and R.\ Pandharipande, {\em The connectedness of
the moduli space of maps to homogeneous spaces}, preprint 2000.

\bibitem{K1}  M.\ Kontsevich, {\em Intersection theory on the moduli
space of curves and the matrix Airy function}, Comm.\ Math.\ Phys.
{\bf 147} (1992), 1-23.

\bibitem{K2}  M.\ Kontsevich, {\em Enumeration of rational curves via
torus actions}, in {\em The moduli space of curves},
(R.\ Dijkgraaf, C.\ Faber, and G.\ van der Geer, eds.), 335-368, Birkh\"auser: Boston,
1995.

\bibitem{KMa} M.\ Kontsevich and Y.\ Manin, {\em
Gromov-Witten classes, quantum cohomology, and enumerative geometry}, Comm.\ Math.\ Phys.\ 
{\bf 164} (1994),
525-562.\ 

\bibitem{Kr} A.\ Kresch, {\em Cycle groups for Artin stacks},
 Invent.\ Math.\ {\bf 138} (1999),  495-536. 

\bibitem{Lab}
G.~Labelle, 
\emph{Une nouvelle d\'emonstration combinatoire des formules 
d'inversion de Lagrange}, Adv.\ in
Math\ \textbf{42} (1981), no.\ 3, 217--247.


\bibitem{LR} A.-M.\ Li and Y.\ Ruan, {\em Symplectic surgery and
Gromov-Witten invariants of Calabi-Yau 3-folds I}, preprint 1998.

\bibitem{LZZ}  A.-M.\ Li, G.\ Zhao, Q.\ Zheng, {\em The number of ramified
covers of a Riemann surface by a Riemann surface}, preprint 1999.


\bibitem{LiT} J.\ Li and G.\ Tian, {\em Virtual moduli cycles
and Gromov-Witten invariants of algebraic varieties},
Jour.\ AMS {\bf{11}} (1998), no.\ 1, 119-174.



\bibitem{Lo} E.\ Looijenga, {\em 
Intersection theory on Deligne-Mumford compactifications (after Witten and Kontsevich)},
S\'eminaire Bourbaki, Vol.
1992/93.\ Astérisque No.\ 216, (1993), Exp.\ No.\ 768, 4, 187-212. 

\bibitem{Lu} T.~\L uczak, 
\emph{Random trees and random graphs},
 Proceedings of the Eighth International Conference "Random Structures and Algorithms" 
(Poznan, 1997). Random
Structures Algorithms \textbf{13} (1998), no.\ 3-4, 485--500.

\bibitem{Ma} Yu.\ Manin, {\em Generating functions in algebraic
geometry and sums over trees},
in {\em The moduli space of curves}, 
(R.\ Dijkgraaf, C.\ Faber, and G.\ van der Geer, eds.), 401-417,  Birkh\"auser: Boston,
1995.

\bibitem{Me} M.~Mehta, 
\emph{Random matrices},  Academic Press,  1991.


\bibitem{MM0}
A.~Meir and J.~W.~Moon,
\emph{The distance between points in random trees}, 
J.\ Combinatorial Theory \textbf{8} (1970), 99--103. 


\bibitem{MM}
A.~Meir and J.~W.~Moon,
\emph{On major and minor branches of rooted trees},
Canad.\ J.\ Math.\ \textbf{39} (1987), no.\ 3, 673--693. 

\bibitem{MSS} S.~Monni, J.~Song, and Y.~Song, {\em The Hurwitz enumeration problem
of branched covers and Hodge integrals}, preprint 2000.

\bibitem{Mu}  D.\ Mumford, {\em Towards an enumerative geometry of
the moduli space of curves}, in {\em Arithmetic and Geometry, Part II},  
(M.\ Artin and J.\ Tate, eds.), 271-328, Birkh\"auser: Boston 1983.


\bibitem{O2}  A.\ Okounkov, {\em Random matrices and random permutations},
Inter.\ Mat.\ Res.\ Notices \textbf{20} (2000), math.CO/9903176.

\bibitem{O3}  A.\ Okounkov, {\em Toda equations for Hurwitz numbers}, Math.\ Res.\ Letters,
{\bf 7} (2000) 447-453. 

\bibitem{O1}  A.\ Okounkov, {\em Generating functions for the intersection numbers
on moduli spaces of curves}, preprint, 2000. 


\bibitem{OP}  A.\ Okounkov and R.\ Pandharipande, in preparation.

\bibitem{P1}  R.\ Pandharipande, {\em The Toda equations and
the Gromov-Witten theory of the Riemann sphere}, preprint 1999.

\bibitem{P2}  R.\ Pandharipande, {\em Hodge integrals and
degenerate contributions}, 
Comm.\ Math.\ Phys.\ {\bf 208} (1999), 489-506.


\bibitem{Pit}
J.~Pitman,
\emph{Enumeration of trees and forests related to branching
processes and random walks},
DIMACS Series in Discrete Mathematics and Theoretical Computer
Science, vol.\ \textbf{41}, 1998. 

\bibitem{Pit2}
J.~Pitman,
\emph{Coalescent random forests},
J.~Combin.\ Theory Ser.\ A \textbf{85} (1999), no.\ 2, 165--193.

\bibitem{RS}
A.~R\'enyi and G.~Szekeres,  
\emph{On the height of trees}. 
J.\ Austral.\ Math.\ Soc.\ \textbf{7} (1967) 497--507. 


\bibitem{Sosh}
A.~Soshnikov,
\emph{Universality at the edge of the spectrum in Wigner random matrices},
 Comm.\ Math.\ Phys.\ \textbf{207} (1999), no.\ 3, 697--733. 

\bibitem{Sta} R.~Stanley,
\emph{Enumerative combinatorics}, Vol.~2, Cambridge Studies in Advanced Mathematics, 62.
Cambridge University Press, Cambridge, 1999.

\bibitem{Str}
V.~Strehl,
\emph{
Minimal transitive products of transpositions---the reconstruction of a proof of A.~Hurwitz},
S\'em.\ Lothar.\ Combin.\ \textbf{37} (1996), Art.\ S37c, 12 pp.\ (electronic).


\bibitem{Th}  R.\ Thomas, {\em Derived categories for the working mathematician},
preprint 1999.


\bibitem{To} B.\ Totaro, {\em Chow ring of the symmetric group}, preprint 1998.

\bibitem{TW}
C.~A.~Tracy and H.~Widom,
\emph{Level-spacing distributions and the Airy kernel},
Commun.\ Math.\ Phys., \textbf{159}, 1994, 151--174.  


\bibitem{Vi} A.\ Vistoli, {\em Intersection theory on algebraic stacks
and their moduli}, Invent.\ Math.\ {\bf 97} (1989), 613-670.



\bibitem{W} E.\ Witten, {\em Two dimensional gravity and intersection
theory on moduli space}, Surveys in Diff.\ Geom.\ {\bf 1} (1991), 243-310.


\bibitem{Z}
A.~Zvonkin,
\emph{Matrix integrals and map enumeration: an
accessible introduction},
Math.\ Comput.\ Modelling, \textbf{26}, 1997,
no.~8--10, 281--304.

\end{thebibliography}
\end{document}